\documentclass{amsart}
\usepackage{amsmath,amsthm,upref}

\newcommand{\ttt}[1]{\hspace{-#1pt}&\hspace{-#1pt},\hspace{-#1pt}&\hspace{-#1pt}}

\newcommand\N{\mathbb N}
\newcommand\ZZ{{\mathbb Z}/2}
\newcommand\Z{\mathbb Z}

\newcommand\R[1]{{\mathbb R}^{#1}}
\newcommand\C[1]{{\mathbb C}^{#1}}
\newcommand\ra{\rightarrow}

\newcommand{\ha}[1]{\overset{\rightharpoonup}{{#1}}}
\renewcommand{\mod}{\text{mod }}
\newcommand{\card}{\operatorname{card}}

\newcommand\cB{{\mathcal{B}}}
\newcommand\cC{{\mathcal{C}}}
\newcommand\cH{{\mathcal{H}}}
\renewcommand\H{{\mathcal{H}}}
\newcommand\cP{{\mathcal{P}}}
\newcommand\q{{\mathcal{Q}}}

\DeclareMathOperator{\Card}{Card}

\newtheorem{Theorem}{Theorem}
\newtheorem*{Thm}{Theorem}

\newtheorem{Proposition}{Proposition}

\newtheorem{Lemma}{Lemma}
\newtheorem*{Lm}{Lemma}
\newtheorem{Corollary}{Corollary}

\theoremstyle{remark}
\newtheorem{Remark}{Remark}

\theoremstyle{definition}
\newtheorem{Definition}{Definition}
\newtheorem{Example}{Example}

\begin{document}
\title
[Components of the moduli space of Abelian differentials]
{Connected  components   of   the   moduli   spaces   of  Abelian
differentials with prescribed singularities}

\author{Maxim Kontsevich}
\address{Institut des Hautes \'Etudes Scientifiques,
Le Bois-Marie, 35 Route de Chartres,
\linebreak
F-91440 Bures-sur-Yvette, France}
\email{maxim@ihes.fr}

\author{Anton Zorich}
\address{
Institut Math\'ematique de Rennes,
Universit\'e Rennes-1,
Campus de Beaulieu,
\hspace*{20pt}  \linebreak
35042 Rennes, cedex, France \hspace*{230pt}
}
\email{Anton.Zorich@univ-rennes1.fr}




\begin{abstract}
Consider the  moduli space of pairs  $(C,\omega)$ where $C$  is a
smooth compact  complex curve of a given genus  and $\omega$ is a
holomorphic 1-form on $C$ with a given list  of multiplicities of
zeroes. We  describe  connected  components  of  this space. This
classification is important  in the study of dynamics of interval
exchange transformations  and billiards in rational polygons, and
in the study of geometry of translation surfaces.
\end{abstract}
\maketitle
\tableofcontents
%
\section{Introduction}
%
\subsection{Stratification  of  the  moduli   space   of  Abelian
differentials}
For integer $g\ge  2$  we define the space  $\H_g$  as the moduli
space of pairs $(C,\omega)$ where $C$ is a smooth compact complex
curve of genus $g$ and  $\omega$  is a holomorphic 1-form on  $C$
(i.e. an Abelian  differential) which is not equal identically to
zero. Obviously, $\H_g$ is a complex algebraic orbifold (in other
words, a  smooth stack) of dimension  $4g-3$. It is  fibered over
the moduli space $\mathcal{M}_g$ of  curves  with  the fiber over
$[C]\in  \mathcal{M}_g$  equal  to  the  punctured  vector  space
$\Gamma(C,\Omega^1_C)\setminus  \{0\}$  (modulo  the action of  a
finite group $Aut(C)$).

Orbifold $\H_g$ is naturally stratified by  the multiplicities of
zeroes of $\omega$. Let $k_1,\dots,k_n$ be a sequence of positive
integers, $n\ge 1$ with the sum $\sum_i k_i$ equal to  $2g-2$. We
denote by $\H(k_1,\dots,k_n)$ the subspace of  $\H$ consisting of
equivalence classes  of  pairs  $(C,\omega)$  where  $\omega$ has
exactly  $n$  zeroes  and  their  multiplicities   are  equal  to
$k_1,\dots,k_n$  (for  some ordering of zeroes). Our notation  is
symmetric,    $\H(k_1,\dots,k_n)=\H(k_{\pi(1)},\dots,k_{\pi(n)})$
for any permutation $\pi\in \mathfrak{S}_n$. One has then
$$
\H_g=\bigsqcup_{\substack{n,\,(k_1,\dots,k_n)\\
k_1\le \dots \le k_n\\k_1+\dots+k_n=2g-2}}
                     \H(k_1,\dots,k_n)
$$
Thus, we have a stratification of the moduli space $\H_g$.  It is
well-known that  each stratum $\H(k_1,\dots,k_n)$ is an algebraic
orbifold of dimension
\begin{equation}
\label{eq:dim:H}
\dim_{\C{}}\H(k_1,\dots,k_n)=2g+n-1
\end{equation}
(see~\cite{M82},  \cite{V82},
\cite{V90}). Moreover,  it  carries  a natural holomorphic affine
structure. Here is the description of this structure.

With any pair  $(C,\omega)$  we associate an element $[\omega]\in
H^1(C,  Zeroes(\omega);\C{})$,  the  cohomology  class  of   pair
$(C,Zeroes(\omega))$ represented by  closed complex-valued 1-form
$\omega$. Locally near  each  point $x$ of $\H(k_1,\dots,k_n)$ we
can identify cohomology spaces $H^1(C, Zeroes(\omega);\C{})$ with
each  other  using  the  Gauss---Manin  connection.  (For  points
$x=(C,\omega)$ with nontrivial symmetry we  would  need  to  pass
first to a finite covering of the neighborhood of $x$).  Thus, we
obtain (locally)  a {\em period mapping} from $\H(k_1,\dots,k_n)$
to a domain in a complex vector space. It is well known that this
mapping  is  holomorphic and locally one-to-one. The pullback  of
the tautological affine structure on $H^1(C,Zeroes(\omega);\C{})$
gives   an   affine  structure   on   $\H(k_1,\dots,k_n)$.   (See
also~\cite{Kat}  for  a related  construction  concerning  smooth
closed 1-forms.)

In general, the strata $\H(k_1,\dots,k_n)$ are  not fiber bundles
over the moduli space of curves $\mathcal{M}_g$. For example, the
dimension of the stratum $\H(2g-2)$ for  $g\ge  2$  equals  $2g$,
while dimension of the moduli  space  of  curves  $\mathcal{M}_g$
equals $3g-3$ which is strictly larger than $2g$ for $g\ge 4$.

The  goal  of  this paper is  to  describe  the  set of connected
components of  all  strata  $\H(k_1,\dots,k_n)$. Surprisingly, we
found that  the answer is  quite complicated, some strata have up
to 3 connected  components. The full description of the connected
components of strata is given in  the section~\ref{ss:main}. This
result was announced in the paper~\cite{K}.

\begin{Remark}
\label{rm:number}
For any  {\em  sequence}  $(k_1,\dots,k_n)$  of positive integers
$k_i\ge   1$    such    that    $\sum_i   k_i=2g-2$   we   define
$\H^{num}(k_1,\dots,k_n)$   the    moduli   space   of    Abelian
differentials on curves  with {\em numbered} zeroes such that the
first    zero    has    multiplicity    $k_1$    etc.    Orbifold
$\H^{num}(k_1,\dots,k_n)$    is    a     finite    covering    of
$\H(k_1,\dots,k_n)$. One can show that preimage  of any connected
component of $\H(k_1,\dots,k_n)$  in $\H^{num}(k_1,\dots,k_n)$ is
connected, i.e.  the  classification  of  connected components is
essentially identical in both cases, no matter whether the zeroes
are numbered or not.
\end{Remark}

\subsection{Applications to interval exchange transformations}

The motivation for  our study came from dynamical systems, namely
from the theory of so called interval exchange transformations.

First of  all, there is an  alternative description of  $\H_g$ in
terms of differential geometry.  Outside  of zeroes of an Abelian
differential $\omega$ one can chose locally  a complex coordinate
$z$ in  such way that  $\omega=dz$. This coordinate is defined up
to  a  constant, $z=z'+const$, so it determines  a  natural  flat
metric $|dz^2|$ on the Riemann surface $C$ punctured at zeroes of
$\omega$. At zero of $\omega$  of  multiplicity  $k_i$ the metric
has a conical singularity with the cone angle $2\pi(k_i+1)$. This
flat metric has trivial holonomy in the group $SO(2)$: a parallel
transport of  a vector tangent to  the Riemann surface  $C$ along
any closed path avoiding conical singularities  brings the vector
back  to  itself.  Thus,  choosing  a  tangent  direction  at any
nonsingular point we can  extend  it using the parallel transport
to all other nonsingular points, getting a smooth distribution on
the punctured  Riemann  surface. This distribution is integrable:
it defines a  foliation with singularities at the conical points.
The oriented foliation defined by the positive real direction $x$
in coordinate $z=x+iy$  is  called {\it horizontal}; the oriented
foliation defined by the positive purely  imaginary direction $y$
is called {\it vertical}. At  a  conical point with a cone  angle
$2\pi(k_i+1)$ one gets $k_i+1$ horizontal (vertical) directions.

Conversely, a flat structure with trivial $SO(2)$-holonomy having
several cone type singularities plus a choice of, say, horizontal
direction uniquely determines a complex structure on the surface,
and an Abelian differential in this complex structure.

An Abelian differential  $\omega$  defines also two smooth closed
real-valued        1-forms        $\omega_v=Re(\omega)$       and
$\omega_h=Im(\omega)$ on  $C$  considered  as  a  smooth oriented
two-dimensional  surface  $M^2$.  The  vertical  and   horizontal
foliations  described  above  are  the kernel foliations  of  the
1-forms $\omega_v$ and $\omega_h$ correspondingly.

Conversely, let $M^2$ be a compact  smooth  oriented  surface  of
genus $g$ with a pair of closed 1-forms $\omega_v, \omega_h$ such
that $\omega_v\wedge\omega_h>0$  everywhere on $M^2$ outside of a
finite set. Then there  is  a unique point $[(C,\omega)]\in \H_g$
producing such $M^2$ with forms $\omega_v,\omega_h$.

There is a non-holomorphic continuous action  on  $\H_g$  of  the
group $GL(2,\R{})_{+}$  (the  group  of  matrices  with  positive
determinants).     In     terms     of    pairs    of     1-forms
$(\omega_v,\omega_h)=\left(Re(\omega),Im(\omega)\right)$     this
action is given simply by linear transformations
$$
(\omega_v,\omega_h)\mapsto (a\omega_v+b\omega_h,c\omega_v+d\omega_h).
$$

Later in  this text we shall use all  descriptions of $\H_g$: the
algebro-geometric one,  the one in  terms of flat surfaces with a
choice of the horizontal direction, and the one in terms of pairs
of measured oriented foliations.

It was  proved  by  H.~Masur,  see~\cite{M82}  and by W.~A.~Veech
(see~\cite{V82}) that for a generic (with respect to the Lebesgue
measure) point  of any stratum $\H(k_1,\dots,k_n)$ the horizontal
foliation (and also the vertical one) is uniquely ergodic. Let us
take any interval  $I$ on the  surface $M^2$ transversal  to  the
vertical foliation,  with  the  canonical induced length element.
The first  return  map  $T:I\longrightarrow  I$  (defined  almost
everywhere on $I$) is an interval exchange map, i.e. a one-to-one
map  with  finitely  many  discontinuity  points  such  that  the
derivative  of  $T$  is  equal  almost  everywhere  to  $+1$. The
interval  exchange  map  is  parametrized  by  the number $m$  of
maximal open subintervals $(I_i)_{i=1,\dots,m}$ of continuity  of
the  transformation  $T$,  by  the sequence of lengths  of  these
subintervals  $\lambda_1,\dots,\lambda_m$  where   $\lambda_i>0$,
$i=1,\dots,m$,  and  by  a  permutation  $\pi\in  \mathfrak{S}_m$
describing the order in which  intervals  $T(I_i)$  are placed in
$I$:  the  $k$-th interval  is  sent to  the  place $\pi(k)$.  It
follows from the  unique ergodicity that the permutation $\pi$ is
{\it  irreducible},  which  means  in our context  that  $\forall
k=1,\dots,m-1$ we have $\pi(\{1,\dots,k\})\ne\{1,\dots,k\}$.

Conversely, for any interval  exchange  map $T$ one can construct
an Abelian differential $\omega$ and a horizontal interval $I$ on
a complex curve $C$  such that the first return map to  $I$ along
the  vertical  foliation  of  $\omega$  is  the  given  map  $T$,
see~\cite{M82},  \cite{V82}.  Though   the  Abelian  differential
$\omega$ is not uniquely determined by the interval exchange map,
the collection of multiplicities of zeroes  $(k_1,\dots, k_n)$ of
$\omega$ and even the connected component  of  the  moduli  space
$\H(k_1,\dots,k_n)$ containing point  $[(C,\omega)]$ are uniquely
determined by the permutation $\pi$, see~\cite{M82},  \cite{V82}.
Thus one may decompose  the  set of irreducible permutations into
groups called  {\it  extended  Rauzy  classes}  corresponding  to
connected components of the strata $\H(k_1,\dots,k_n)$.

The application of our result to the theory  of interval exchange
maps is based  on  the corollary  of  the fundamental theorem  of
H.Masur~\cite{M82} and W.Veech~\cite{V82} which we present in the
next section. The  corollary  is as follows: dynamical properties
of a  generic interval exchange map  depend only on  the extended
Rauzy class of the permutation  of  subintervals.  Genericity  is
understood here with respect to the Lebesgue measure on the space
${\mathbb R}_+^m$ parameterizing lengths $(\lambda_i)_{1\le i \le
m}$ of subintervals under exchange.

Actually, the extended Rauzy classes can  be  defined  in  purely
combinatorial terms,  see Appendix~\ref{a:rauzy:cl} for  details.
Thus  the  problem  of  the  description  of  the  extended Rauzy
classes, and hence, of the description of connected components of
the strata  of  Abelian  differentials,  is purely combinatorial.
However, it seems to be  very  hard to solve it directly.  Still,
for  small  genera  the  problem  is  tractable.  W.Veech  showed
in~\cite{V90}  that   the   stratum  $\H(4)$  has  two  connected
components.  P.Arnoux  proved that the stratum $\H(6)$ has  three
connected components.

In the  present paper we  give a classification of extended Rauzy
classes using not  only combinatorics but also tools of algebraic
geometry, topology and of dynamical systems.

\subsection{Ergodic  components  of  the  Teichm\"uller  geodesic
flow}

There  is a  natural  immersion of the  moduli  space of  Abelian
differentials into  the  moduli  space  of  holomorphic quadratic
differentials:  we  associate  to  an  Abelian  differential  its
square. With every quadratic differential we  can again associate
a   flat   metric   with   conical   singularities.   The   group
$GL(2,\R{})_+$ acts  naturally  on  this  larger  moduli space as
well; this action leaves the immersed  moduli  space  of  Abelian
differentials invariant,  moreover,  on  the immersed subspace it
coincides with the action  defined  in the previous section. This
action preserves the  natural  stratification of the moduli space
of quadratic differentials by multiplicities of zeroes.

The  action   of  the  diagonal  subgroup  of  $SL(2,\R{})\subset
GL(2,\R{})_{+}$ on  the  moduli  space of quadratic differentials
can be naturally identified with the geodesic flow  on the moduli
space of  curves  for  Teichm\"uller  metric  (which is piecewise
real-analytic   Finsler   metric   on   $\mathcal{M}_g$).   Group
$SL(2,\R{})$  preserves  the hypersurface in the moduli space  of
quadratic  differentials  consisting of those ones for which  the
associated flat metric has the total area equal to $1$.

Numerous important results in the  theory  of  interval  exchange
maps, of measured foliations, of billiards  in rational polygons,
of dynamics on translation surfaces  are  based  on the following
fundamental     observation     by     H.Masur~\cite{M82}     and
W.Veech~\cite{V82}:

\begin{Thm}[H.~Masur; W.~Veech]
The  Teichm\"uller   geodesic  flow  acts  ergodically  on  every
connected  component  of  every  stratum of the moduli  space  of
quadratic  differentials  with  total  area  equal  to  $1$;  the
corresponding  invariant  measure on  the  stratum  is  a  finite
Lebesgue equivalent measure.
\end{Thm}

Thus our classification  of connected components of the strata of
Abelian  differentials  gives  the   classification   of  ergodic
components of the Teichm\"uller geodesic  flow  on  the strata of
squares of Abelian differentials in the moduli space of quadratic
differentials.

The complete classification of connected components  of strata of
quadratic differentials  is  in  progress  (see  an  announcement
in~\cite{L01}). For  example,  the  stratum  of  those  quadratic
differentials  on  a  curve  of  genus  $g=4$,  which  cannot  be
represented as a square of an  Abelian  differential,  and  which
have a single zero of degree $12$, has  two connected components,
but at the moment a topological invariant which would distinguish
representatives  of  these two connected components is not  known
yet.

In general, it seems to be very interesting to describe invariant
submanifolds (closures  of  orbits,  invariant  measures) for the
action  of  $GL(2,\R{})_{+}$  on  the  moduli  spaces.  Connected
components  of  the  strata  are  only   the  simplest  invariant
submanifolds,  there   are   many   others.   For   example   the
Teichm\"uller  disks  of Veech curves form the smallest  possible
invariant submanifolds.

One  can  use   a  submanifold  invariant  under  the  action  of
$GL(2,\R{})_{+}$  to  produce  other  invariant  submanifolds  in
higher genera  applying some fixed ramified covering construction
to all  pairs  $(C,\omega)$  constituting  the  initial invariant
submanifold. In section~\ref{ss:hyp:com} we use a particular case
of this construction to define some  special connected components
of some strata.

\section{Formulation of results}

\subsection{Hyperelliptic components}
\label{ss:hyp:com}

First  of  all,  we  introduce the moduli spaces  of  meromorphic
quadratic differentials.

\begin{Definition}
For integer $g\ge 0$ and  collection  $(l_1,\dots,  l_n),n\ge  1$
such  that  $l_j\ge  -1,  l_j\ne  0$  for  all  $j$  and  $\sum_j
l_j=4g-4$, denote  by  $\q(l_1,\dots,l_n)$  the  moduli  space of
pairs $(C,\phi)$  where $C$ is  a smooth compact complex curve of
genus $g$ and $\phi$  is  a meromorphic quadratic differential on
$C$ with  zeroes of orders  $l_j$ (simple poles if $l_j=-1$) such
that  $\phi$   is  not  equal  to   the  square  of   an  Abelian
differential.
\end{Definition}

It  is  known  (see~\cite{V90})  that  $\q(l_1,\dots,l_n)$  is  a
complex algebraic orbifold of dimension
\begin{equation}
\label{eq:dim:Q}
\dim_{\C{}}\q(l_1,\dots,l_n)=2g+n-2
\end{equation}

Sometimes  we   shall  use  ``exponential''  notation  to  denote
multiple zeroes (simple poles) of  the  same  degree, for example
$\q(-1^5,1):=\q(-1,-1,-1,-1,-1,1)$. The condition  that $\phi$ is
not  a  square is  automatically  satisfied if  at  least one  of
parameters $l_j$ is odd.

One can  canonically  associate  with every meromorphic quadratic
differential $(C,\phi)$  another  connected  curve  $C'$  with an
Abelian differential $\omega$ on  it.  Namely, $C'$ is the unique
double  covering  of $C$  (maybe  ramified  at  singularities  of
$\phi$), such that  the  pullback  of $\phi$ is a  square  of  an
Abelian    differential    $\omega$.   We    have   automatically
$\sigma^\ast(\omega)=-\omega$ where $\sigma$ is the involution on
$C'$ interchanging points in the  generic  fiber  over $C$. Curve
$C'$ is  connected because of the condition that  $\phi$ is not a
square of an Abelian differential.

Thus, we obtain a map from  the  stratum  $\q(l_1,\dots,l_n)$  of
meromorphic    quadratic    differentials    to    the    stratum
$\H(k_1,\dots,k_m)$  of  Abelian   differentials,  where  numbers
$(k_i)$ are obtained from $(l_j)$ by the following  rule: to each
even $l_j>0$ we associate a pair of zeroes of $\omega$  of orders
$(l_j/2,l_j/2)$  in  the  list  $(k_i)$, to each odd  $l_j>0$  we
associate one zero of order $l_j+1$,  and  associate  nothing  to
simple poles (e.g. to $l_j=-1$).

\begin{Lemma}
\label{lm:Q:to:H}
The canonical map
described above
$$
\q(l_1,\dots,l_n)\to\H(k_1,\dots,k_m)
$$
is an immersion.
\end{Lemma}
\begin{proof}
Denote as  above by $C'$ the double covering  of $C$ with Abelian
differential $\omega$ and involution $\sigma$.

Consider the induced involution
$$
\sigma^{\ast}: H^1(C,Zeroes(\omega);\C{})\to
H^1(C,Zeroes(\omega);\C{})
$$
It   defines   decomposition    $H^1(C,Zeroes(\omega);\C{})\simeq
V_1\oplus V_{-1}$ of the first cohomology into the  direct sum of
subspaces invariant  and  anti  invariant  under  the  involution
$\sigma^{\ast}$. By construction $[\omega]\in  V_{-1}$.  Thus, we
obtain (locally) a  mapping  from $\q(l_1,\dots,l_n)$ to a domain
in     the     complex      vector     space     $V_{-1}\subseteq
H^1(C,Zeroes(\omega);\C{})$. It  is  well known that this mapping
is  holomorphic   and   locally   one-to-one.   Since  the  space
$\H(k_1,\dots,k_m)$      is      locally     identified      with
$H^1(C,Zeroes(\omega);\C{})$ by means of the period mapping, this
completes the proof of Lemma. \end{proof}

The  following  two  series of maps  of  this  kind  would play a
special role for us:
\begin{equation}
\label{eq:hyp:series}
\begin{array}{lcl}
\q(-1^{2g'+1},2g'-3) &\ra &  \H(2g'-2)\\
\q(-1^{2g'+2},2g'-2) &\ra &  \H(g'-1,g'-1),
\end{array}
\end{equation}
where  $g'\ge 2$ in  both  cases.  In  both cases  curve  $C$  is
rational (i.e. $g=0$), and  hence  curve $C'$ is hyperelliptic of
genus $g'$. In these two cases the dimension of the image stratum
of Abelian  differentials  coincides  with  the  dimension of the
original stratum of meromorphic quadratic differentials.  Indeed,
formula~\eqref{eq:dim:Q} gives
\begin{align*}
&\dim_{\C{}}\q(-1^{2g'+1},2g'-3)=2\cdot 0+(2g'+2)-2=2g'\\
&\dim_{\C{}}\q(-1^{2g'+2},2g'-2)=2\cdot 0+(2g'+3)-2=2g'+1,
\end{align*}
while formula~\eqref{eq:dim:H} gives the following dimensions  of
the image strata:
\begin{align*}
&\dim_{\C{}}\H(2g'-2)=2g'+1-1=2g'\\
&\dim_{\C{}}\H(g'-1,g'-1)=2g'+2-1=2g'+1.
\end{align*}
\begin{Remark}
We have constructed a map $\q(l_1,\dots,l_n)\to\H(k_1,\dots,k_m)$
using certain canonical double covering $C'\to  C$. Choosing some
other (ramified) covering of some  fixed  type  one can construct
some other (local) maps between moduli  spaces  of  quadratic  or
Abelian differentials. The reader can find a detailed description
of all maps of this kind between moduli spaces of  {\it quadratic
differentials},  which   give   coincidence   of  dimensions,  in
paper~\cite{L01}.
\end{Remark}

Before returning  to  maps~\eqref{eq:hyp:series}  which  are of a
particular  interest  for  us  we  need  to  prove  the following
statement.

\begin{Proposition}
\label{pr:genus:0}
In  the   case   $g=0$   every   stratum  $\q(l_1,\dots,l_n)$  of
meromorphic quadratic differentials is nonempty and connected.
\end{Proposition}
\begin{proof}
For any divisor on ${\mathbb C}P^1$ with given multiplicities the
corresponding meromorphic  quadratic  differential  exists and is
unique up to a non-zero scalar. Thus, we have
$$
\q(l_1,...,l_n)/\C{\ast} \cong
\left( (\C{}P^1)^n\backslash diagonals\right) / \left( PSL(2,\C{})\times
(\text{finite symmetry group})
\right)
$$
Therefore  the  orbifold  $\q(l_1,\dots,l_n)$  is  nonempty   and
connected.
\end{proof}

Lemma~\ref{lm:Q:to:H},   the  observation   on   coincidence   of
dimensions of the  corresponding strata in~\eqref{eq:hyp:series},
together with Proposition~\ref{pr:genus:0}  justify the following
definition.

\begin{Definition}
By {\it hyperelliptic components} we call the following connected
components  of  the following strata of Abelian differentials  on
compact complex curves of genera $g\ge 2$:

The   connected   component   $\H^{hyp}(2g-2)$  of  the   stratum
$\H(2g-2)$ consisting of Abelian  differentials  on hyperelliptic
curves of genus $g$ corresponding to  the orbifold $\q(-1^{2g+1},
2g-3)$;

The  connected  component  $\H^{hyp}(g-1,g-1)$  of  $\H(g-1,g-1)$
corresponding to the orbifold $\q(-1^{2g+2},2g-2)$.
\end{Definition}

\begin{Remark}
Points of $\H^{hyp}(2g-2)$  (respectively of $\H^{hyp}(g-1,g-1)$)
are  Abelian  differentials on hyperelliptic curves of genus  $g$
which have a single  zero  of multiplicity $2g-2$ invariant under
the  hyperelliptic  involution (respectively a pair of zeroes  of
orders  $g-1$  symmetric  to  each  other  with  respect  to  the
hyperelliptic involution).

Note that if an Abelian differential on a hyperelliptic curve has
a  single zero  of  order $2g-2$ then  this  zero is  necessarily
invariant under  the  hyperelliptic  involution $\sigma$, because
$\sigma^\ast(\omega)=-\omega$   for   any  Abelian   differential
$\omega$. Therefore,  this  Abelian  differential  belongs to the
component $\H^{hyp}(2g-2)$. However, if  an  Abelian differential
$\omega$  has  two   zeroes  of  degrees  $g-1$,  there  are  two
possibilities:  the   zeroes   might   be   interchanged  by  the
hyperelliptic involution, and  they  might be invariant under the
hyperelliptic  involution.   In   the   first  case  the  Abelian
differential belongs to the component $\H^{hyp}(g-1,g-1)$,  while
in the second case it does not.
\end{Remark}

\subsection{Parity of a spin structure: a definition}

\begin{Definition}
A {\em spin structure} on a smooth compact complex curve $C$ is a
choice of a  half  of  the canonical class, i.e.  of  an  element
$\alpha\in Pic(C)$ such that $$2\alpha=K_C:=-c_1(T_C)$$ The  {\em
parity  of  the spin structure} is  the  residue modulo 2 of  the
dimension $$dim\,\Gamma(C,L)  =  dim\, H^0(C,L)$$ for line bundle
$L$ with $c_1(L)=\alpha$.
\end{Definition}

On a curve of genus  $g\ge  1$ there are $2^{2g}$ different  spin
structures  among  which  $2^{2g-1}  +  2^{g-1}$   are  even  and
$2^{2g-1}-2^{g-1}$  are  odd. It  follows  from  the  results  of
M.Atiyah~\cite{At} and  D.Mumford~\cite{Mum} that the parity of a
spin structure is invariant under continuous deformations.

Let  $\omega$  be   an   Abelian  differential  with  {\em  even}
multiplicities of zeroes, $k_i=2 l_i$ for all $i$, $i=1,\dots,n$.
The divisor of zeroes of $\omega$
$$
\text{Zeroes}(\omega)=2l_1 P_1+\dots+2l_n P_n
$$
represents the canonical class  $K_C$.  Thus, we have a canonical
spin structure  on $C$ defined by
$$ \alpha_\omega:= [l_1 P_1 + \dots +l_n P_n]\in Pic(C) $$
By continuity  the parity of this  spin structure is  constant on
each connected component of stratum $\H(2l_1,\dots,2l_n)$.

\begin{Definition}
We      say      that     a      connected      component      of
$\mathcal{H}(2l_1,\dots,2l_n)$ has  {\it  even}  or {\it odd spin
structure} depending on  whether  $\alpha_\omega$ is even or odd,
where $\omega$ belongs to the corresponding connected component.
\end{Definition}

In section~\ref{ss:top:def} we present  an  equivalent definition
of  the  parity  of   spin   structure  in  terms  of  elementary
differential topology.

\subsection{Main results}
\label{ss:main}

First of  all, we describe  connected components of strata in the
``stable  range''  when the genus of the  curve  is  sufficiently
large.

\begin{Theorem}
\label{th:ccor}
All connected components of any stratum  of Abelian differentials
on a curve of genus $g\ge 4$ are described by the following list:

The stratum $\mathcal{H}(2g-2)$ has  three  connected components:
the hyperelliptic one, $\mathcal{H}^{hyp}(2g-2)$,  and  two other
components:            $\mathcal{H}^{even}(2g-2)$             and
$\mathcal{H}^{odd}(2g-2)$  corresponding  to even  and  odd  spin
structures.

The stratum  $\mathcal{H}(2l,2l)$,  $l\ge  2$ has three connected
components:  the  hyperelliptic one,  $\mathcal{H}^{hyp}(2l,2l)$,
and   two   other  components:   $\mathcal{H}^{even}(2l,2l)$  and
$\mathcal{H}^{odd}(2l,2l)$.

All the  other strata of the form $\mathcal{H}(2l_1,\dots,2l_n)$,
where  all   $l_i\ge   1$,   have   two   connected   components:
$\mathcal{H}^{even}(2l_1,\dots,2l_n)$                         and
$\mathcal{H}^{odd}(2l_1,\dots,2l_n)$, corresponding to  even  and
odd spin structures.

The strata $\mathcal{H}(2l-1,2l-1)$, $l\ge 2$, have two connected
components;  one   of  them:  $\mathcal{H}^{hyp}(2l-1,2l-1)$   is
hyperelliptic;  the  other  $\mathcal{H}^{nonhyp}(2l-1,2l-1)$  is
not.

All the  other strata of Abelian  differentials on the  curves of
genera $g\ge 4$ are nonempty and connected.
\end{Theorem}

Finally we consider the list of connected components  in the case
of small genera $1\le g\le 3$, where some  components are missing
in comparison with the general case.

\begin{Theorem}
\label{th:g123}

The moduli space  of Abelian differentials  on a curve  of  genus
$g=2$    contains    two     strata:    $\mathcal{H}(1,1)$    and
$\mathcal{H}(2)$.  Each  of them is connected and coincides  with
its hyperelliptic component.

Each of  the  strata  $\mathcal{H}(2,2)$, $\mathcal{H}(4)$ of the
moduli space  of Abelian differentials  on a curve of genus $g=3$
has two  connected  components:  the  hyperelliptic  one, and one
having odd spin structure. The other  strata  are  connected  for
genus $g=3$.
\end{Theorem}

Parities  of   spin   structures  for  hyperelliptic  strata  are
calculated    in    the     Appendix~\ref{s:hypel},     Corollary
\ref{cr:hspin}.

Theorems~\ref{th:ccor}    and~\ref{th:g123}    were     announced
in~\cite{K}.

\subsection{Plan of the proof}

We possess two invariants of connected components: the components
could be either  hyperelliptic  or not, and in  the  case of even
multiplicities the associated spin structure could be either even
or  odd.  We show that these invariants  classify  the  connected
components. The maximal  number of connected components is 3, and
it is achieved for the  strata  $\H(2g-2)$ for $g\ge 4$. We  call
the stratum $ \H(2g-2)$ {\em minimal}.

Our plan of the proof is the following:

In section~\ref{s:spin} we give an alternative description of the
parity of the spin  structure  defined by an Abelian differential
having zeroes of  even degrees. For  a special class  of  Abelian
differentials   introduced    in   section~\ref{s:toolkit}   this
description  in  terms  of  differential topology will  make  the
computation of the parity of the spin structure especially easy.

The subset  of  points  $[(C,\omega)]$ whose horizontal foliation
has  only {\em  closed  leaves}, is dense  in  every stratum.  In
section~\ref{ss:diag} we consider Abelian  differentials  only of
this  type.  We  propose  a combinatorial way to  represent  such
Abelian  differentials   by  diagrams,  and  it  is  particularly
convenient for  the  minimal stratum. In section~\ref{ss:diag} we
establish a criterion for diagrams selecting  the ones associated
to Abelian  differentials.  We  call  corresponding diagrams {\it
realizable}. Also in section~\ref{ss:diag}  we  describe diagrams
corresponding to hyperelliptic
Abelian differentials.

We  complete  section~\ref{s:toolkit} by  introducing  a  surgery
(``bubbling a handle'') which  allows  us to construct an Abelian
differential  in  the  minimal  stratum  in  genus $g+1$ from  an
Abelian differential from  the minimal stratum in genus $g$. This
surgery can be applied to any Abelian differential; however, when
the horizontal  foliation  of  an  Abelian  differential has only
closed leaves, one  can  apply the surgery in  such  way that the
horizontal foliation  of  the resulting Abelian differential also
has only  closed leaves. In  this particular case the surgery can
be  described in  terms  of diagrams. Also  we  describe how  the
parity of the spin structure changes under the surgery.

In section~\ref{s:conn:comp} we prove the classification theorem.
First   we   prove   it   for    the    minimal    stratum.    In
section~\ref{ss:cc:minimal} we study  possible transformations of
realizable diagrams  representing  points  in the minimal stratum
preserving the  connected  component.  We  prove  by induction in
genus $g\ge 2$ that the classification of connected components of
the minimal stratum  $  \H(2g-2)$ is as in Theorems~\ref{th:ccor}
and~\ref{th:g123}. We have  to note that  a surgery used  in  the
step  of  induction  (``tearing  off  a  handle'')  is  based  on
combinatorial              Lemma~\ref{lm:m1adv}              from
appendix~\ref{ss:comb:rauzy:cl}    concerning    extended   Rauzy
classes.

In section~\ref{ss:local:stratification} we study the topology of
the adjacency  of strata, and  prove that the number of connected
components in every stratum adjacent  to  the  minimal stratum is
bounded  above  by  the  number  of  connected components of  the
minimal  stratum.  More precisely, we identify the  set  of  such
components with a quotient of the set $\pi_0(\H(2g-2))$.

In section~\ref{ss:merge}  we prove that any Abelian differential
which does not belong to the minimal stratum,  can be degenerated
to a  differential with less  zeroes. Thus, by induction we prove
that any  connected component of any  stratum is adjacent  to the
minimal  stratum,   opening   the   way   to   apply  results  of
section~\ref{ss:local:stratification}.

Using another class of transformations  of  diagrams  we prove in
section~\ref{ss:other:strata} that in certain cases two connected
components  of  any  stratum  adjacent  to  two  given  different
components of the minimal stratum coincide.

Using the previous results we prove the {\em upper} bound  on the
number of connected components of every  stratum.  On  the  other
hand, topological invariants plus a realization construction (see
the end  of  section~\ref{ss:other:strata})  give  a  {\em lower}
bound on the number of components.  These  two  bounds  coincide,
thus we obtain the main result.

Although we shall not do it explicitly in the present  paper, one
can easily  modify the proof for the case  of numbered zeroes and
obtain   essentially   the   same  classification  of   connected
components.

\section{Spin structure determined by an Abelian differential}
\label{s:spin}
In this  section we give an  alternative description of  the spin
structure  determined  by an Abelian differential with zeroes  of
even orders on a closed complex curve.


\subsection{Spin structure: topological definition}
\label{ss:top:def}

We  begin  by recalling the topological definition  of  the  spin
structure on  a  Riemann surface (see~\cite{Mil}, \cite{At}). Let
$M^2_g$ be  a Riemann surface of genus  $g$, and  let $P$ be  the
$S^1$-bundle  of   directions  of  non-zero  tangent  vectors  to
$M^2_g$. A {\it spin  structure}  on $M^2_g$ is a double-covering
$Q\to P$ whose restriction to each fiber of $P$ is  isomorphic to
the standard double covering $S^1 \stackrel{\ZZ}\rightarrow S^1$.

Since the structure group of the covering $Q\to P$ is just $\ZZ$,
the spin structures are in the one-to-one correspondence with the
$\ZZ$-valued  linear  functions  on $H_1(P;\ZZ)$, having  nonzero
value on the  cycle representing the  fiber $S^1$ of  $P$.  Thus,
spin structures are  classified by a coset of $H^1(M^2_g;\ZZ)$ in
$H^1(P;\ZZ)$.

In~\cite{J} D.Johnson  associates to every spin structure $\xi\in
H^1(P;\ZZ)$  on  a Riemann surface a $\ZZ$-valued quadratic  form
$\Omega_\xi$ on $H_1(M^2_g;\ZZ)$,  and  shows, that the parity of
the spin  structure  $\xi$  coincides  with  the Arf-invariant of
$\Omega_\xi$. We present briefly a  sketch  of  the  construction
from~\cite{J}.

First of all, there is a canonical lifting $c \mapsto \tilde{c}$,
$ H_1(M^2_g;\ZZ) \to H_1(P;\ZZ) $ (a map of sets) defined  in the
following  way.  Having  a  cycle $c\in H_1(M^2_g;\ZZ)$  one  can
represent it by a collection of  simple  closed  oriented  curves
$c=\sum_{i=1}^m[\alpha_i]$. Let $[\ha{\alpha_i}]$ be the cycle in
$H_1(P;\ZZ)$  represented  by  the  {\it  framed  curve}  in  $P$
consisting of positive  tangent  directions to $\alpha_i$. Let $z
\in H_1(P;\ZZ)$ be the homology  class  represented  by the fiber
$S^1$. The lifting is defined as
$$
c \mapsto \tilde{c} := \sum_{i=1}^m [\ha{\alpha_i}] + mz
$$
According to~\cite{J} the  map is well-defined. The map obeys the
following relation
$$
\widetilde{a+b}=\tilde{a}+\tilde{b}+(a\cdot b)z
$$
where $(a\cdot  b)$ is the intersection  index of cycles  $a$ and
$b$.  Notice that  the  lifting is {\it  not}  a homomorphism  of
groups.

A $\ZZ$-valued  quadratic  form $\Omega$ on $H_1(M^2_g;\ZZ)$ with
the  associated  bilinear  form  $(a,b)\mapsto a\cdot b$  is  any
function $\Omega: H_1(M^2_g;\ZZ) \to \ZZ$ such that
$$ \Omega(a+b) = \Omega(a)+\Omega(b)+ a\cdot b $$
Having a spin  structure $\xi\in H^1(P;\ZZ)$ one associates to it
the following quadratic form $\Omega_\xi$ on $H_1(M^2_g;\ZZ)$:
$$ \Omega_\xi(a) :\overset{def}{=}
\langle\xi,\tilde{a}\rangle $$
Given a symplectic basis $a_i,b_i  \in  H_1(M^2_g;\Z)$  the  {\it
Arf-invariant} of a quadratic form $\Omega_\xi$ is determined as
$$ \Phi(\Omega_\xi) :\overset{def}{=} \sum_{i=1}^g
\Omega_\xi(a_i)\Omega_\xi(b_i) (\mod 2) $$
It is  proved in~\cite{J} that the  parity of the  spin structure
$\xi$ coincides with the Arf-invariant of $\Omega_\xi$.

\subsection{Spin structure determined by an Abelian differential}

Consider an Abelian  differential  $\omega$ having zeroes of even
degrees $(2l_1,\dots,2l_n)$  on  a  Riemann  surface  $M^2_g$. It
determines  a   flat   structure   on   $M^2_g$   with  cone-type
singularities.  Recall,   that   this  flat  metric  has  trivial
holonomy.  In   particular,   outside   of   finite   number   of
singularities (corresponding to zeroes of  $\omega$)  we  have  a
well-defined  horizontal  direction.  Consider  a  smooth  simple
closed oriented curve  $\alpha$ on $M^2_g$ which does not contain
any zeroes of $\omega$. The flat structure allows us to determine
the index  $ind_\alpha\in \Z$ of  the field tangent to the curve;
$ind_\alpha$ coincides with the degree of the corresponding Gauss
map: the total change of the angle between the vector  tangent to
the curve, and the vector tangent to the  horizontal foliation is
equal to $2\pi\cdot ind_\alpha$.

The spin structure $\xi\in H^1(P;\ZZ)$ determined by $\omega$ has
the following property:
$$
\langle   \xi,   \widetilde{[\alpha]}   \rangle  +1  \equiv
\langle \xi , [\ha{\alpha}] \rangle = ind_\alpha (\mod 2)
$$

This property can be considered  as  a  topological definition of
the spin  structure  determined  by  an  Abelian differential. It
gives  also  the following effective way to  compute  the  parity
$\varphi(\omega)$ of  the  spin  structure  defined  by $\omega$:
choose oriented smooth paths $(\alpha_i, \beta_i)_{i=1,\dots, g}$
representing a symplectic basis of $H_1(M^2_g,\ZZ)$. Then
\begin{multline}
\varphi(\omega):=
\Phi(\Omega_\xi) =
\sum_{i=1}^g\Omega_\xi([\alpha_i])\cdot\Omega_\xi([\beta_i])\ (\mod 2) =
\sum_{i=1}^g\langle\xi,\widetilde{[\alpha_i]}\rangle
    \cdot\langle\xi,\widetilde{[\beta_i]}\rangle\ (\mod 2) =\\
=\sum_{i=1}^g(ind_{\alpha_i} +1)(ind_{\beta_i} +1)\ (\mod 2)
\label{eq:parity:of:spin:structure}
\end{multline}

In particular, using this definition it is easy  to calculate the
parity  of  the  spin  structure given any permutation  from  the
corresponding Rauzy class.

We complete this section with the following obvious statements.

\begin{Lemma}
\label{lm:ind0}
Let $\alpha$ be  a smooth simple closed oriented curve everywhere
transversal  to   the   horizontal   (vertical)  foliation.  Then
$ind_\alpha= 0$. Let  $\alpha$ be a  closed regular leaf  of  the
horizontal (vertical) foliation. Then $ind_\alpha= 0$.
\end{Lemma}

\begin{Lemma}
\label{lm:g1}
The spin structure  of an Abelian  differential on a  surface  of
genus one is always odd.
\end{Lemma}
\begin{proof}
An  Abelian  differential  $\omega$  on  a  surface of genus  one
defines  a  flat  metric  on  the  torus.  One  can  represent  a
symplectic basis of cycles on this flat torus by a pair of closed
geodesics  $\alpha,   \beta$.   By   Lemma~\ref{lm:ind0}  we  get
$ind_\alpha=ind_\beta=0$.                                   Thus,
formula~\ref{eq:parity:of:spin:structure}  gives  the   following
value for  the  parity  $\varphi(\omega)$  of  the spin-structure
defined by $\omega$
$$
\varphi(\omega) =(ind_\alpha +1)(ind_\beta +1)\ (\mod 2) = 1
$$
\end{proof}

\section{Preparation of a surgery toolkit}
\label{s:toolkit}

\subsection{Separatrix diagrams}
\label{ss:diag}

In  this  section  we   consider   a  special  class  of  Abelian
differentials.  Namely,  we   assume   that  all  leaves  of  the
horizontal foliation are either closed or connect critical points
(a  leaf joining  two  critical points is  called  a {\it  saddle
connection}  or  a  {\it  separatrix}). Later we will  be  saying
simply that the  horizontal foliation has only closed leaves. The
square  of  an  Abelian  differential having this property  is  a
particular  case  of   Jenkins--Strebel  quadratic  differential,
see~\cite{Strebel}.

\begin{Lemma}
\label{lm:closedleaves}
Abelian differentials whose horizontal (vertical) foliations have
only  closed  leaves  form  a dense subset in  arbitrary  stratum
$\H(k_1,\dots,k_n)$.
\end{Lemma}
\begin{proof}
We  prove  the statement for horizontal foliations; for  vertical
foliations it is completely analogous.  First  of  all, using the
period   mapping   one   concludes   immediately   that    points
$[(C,\omega)]$  with  rational periods  of $\omega_h:=Im(\omega)$
are dense  in arbitrary stratum. We  claim that for  these points
the horizontal foliation  (given by the kernel of $\omega_h$) has
only  closed  leaves.  The  reason  is  that  in  this  case  the
integration   of   $\omega_h$   gives   a   smooth   proper   map
$\pi:C\rightarrow \R{}  /  {\frac{1}  {N}}  \Z  \simeq S^1$ where
$N\in \N$ is a common denominator of periods  of $\omega_h$. Also
we have  the  equality  $\omega_h=\pi^\ast(dy)$ where $y$ denotes
the standard  coordinate on the real  line $\R{}$. All  leaves of
the  horizontal  foliations belong to fibers of $\pi$,  therefore
are either closed or connect critical points.
\end{proof}

We will  associate  with  each  Abelian differential $(C,\omega)$
whose horizontal foliation has only closed leaves a combinatorial
data called {\it separatrix diagram}.

We start with an informal explanation. Consider the  union of all
saddle connections  for  the  horizontal  foliation,  and add all
critical points (zeroes of $\omega$). We obtain a finite oriented
graph $\Gamma$. Orientation on the edges comes from the canonical
orientation of the horizontal foliation. Moreover, graph $\Gamma$
is drawn  on an oriented  surface, therefore it carries so called
{\it ribbon structure} (even  if  we forget about the orientation
of edges), i.e. on the star  of each vertex $v$ a cyclic order is
given, namely the clockwise order in which edges  are attached to
$v$. The direction of  edges  attached to $v$ alternates (between
directions toward  $v$ and from $v$)  as we follow  the clockwise
order.

It is  well known that any  finite ribbon graph  $\Gamma$ defines
canonically  (up  to an isotopy) an oriented surface  $S(\Gamma)$
with boundary. To  obtain this surface  we replace each  edge  of
$\Gamma$  by  a thin oriented strip (rectangle)  and  glue  these
strips  together  using  the  cyclic  order  in  each  vertex  of
$\Gamma$. In  our case surface $S(\Gamma)$  can be realized  as a
tubular $\varepsilon$-neighborhood (in the  sense  of transversal
measure) of the  union of all saddle connections for sufficiently
small $\varepsilon>0$.

The  orientation  of  edges   of   $\Gamma$  gives  rise  to  the
orientation of  the  boundary  of  $S(\Gamma)$.  Notice that this
orientation is {\it not} the same as the canonical orientation of
the boundary of  an  oriented surface. Thus, connected components
of the boundary of  $S(\Gamma)$  are decomposed into two classes:
positively   and   negatively   oriented  (positively  when   two
orientations of  the boundary components coincide and negatively,
when  they  are   different).   The  complement  to  the  tubular
$\varepsilon$-neighborhood of $\Gamma$ is a finite disjoint union
of  open  cylinders  foliated  by  oriented  circles. It gives  a
decomposition    of     the     set     of    boundary    circles
$\pi_0(\partial(S(\Gamma)))$  into  pairs  of  components  having
opposite signs of the orientation.

Now we are ready to give a formal definition:

\begin{Definition}
A  {\it separatrix  diagram}  (or simply a  {\it  diagram}) is  a
finite oriented ribbon graph $\Gamma$, and a decomposition of the
set of boundary components of $S(\Gamma)$ into pairs, such that
\begin{enumerate}
\item
the orientation of edges at any vertex is alternated with respect
to the cyclic order of edges at this vertex;
\item  there  is  one  positively  oriented  and  one  negatively
oriented boundary component in each pair.
\end{enumerate}
\end{Definition}

Notice that ribbon graphs which appear as a part of the structure
of a separatrix diagram are  very  special. Any vertex of such  a
graph has  even degree, and  the number of boundary components of
the associated surface with  boundary  is even. Notice also, that
in general the graph of a separatrix diagram is {\it not} planar.

Any separatrix  diagram  $(\Gamma,  pairing)$  defines  a  closed
oriented surface together with an embedding of $\Gamma$  (up to a
homeomorphism) into this  surface. Namely, we glue to the surface
with boundary  $S(\Gamma)$  standard oriented cylinders using the
given pairing.

In pictures representing diagrams we  encode  the  pairing on the
set of boundary components painting corresponding  domains in the
picture by some  colors (textures in the black-and-white text) in
such a way that every  color  appears exactly twice. We will  say
also that paired components have the {\it same color}.

\begin{figure}[htb]
%
\includegraphics{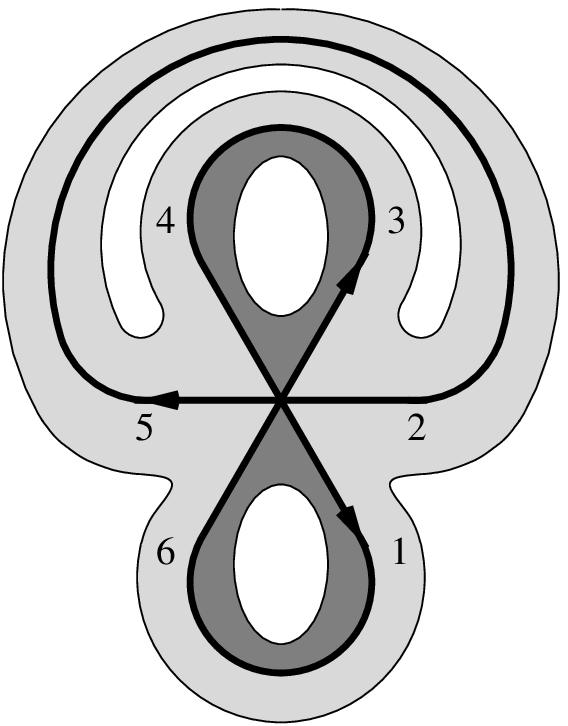}
\includegraphics{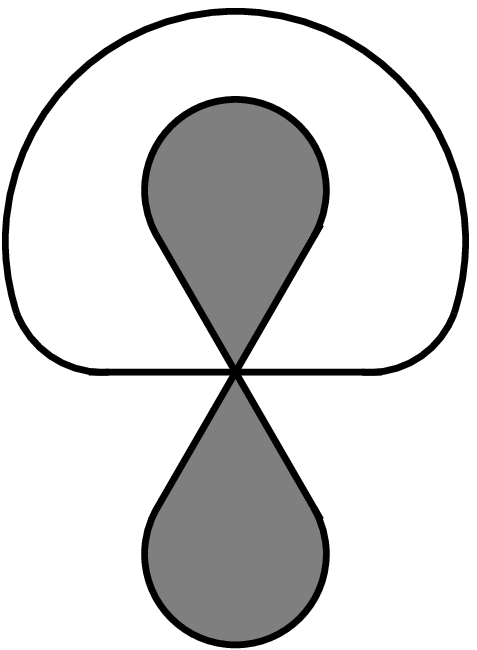}
\vspace{160bp}
\caption{
\label{pic:g2}
An example of  a separatrix diagram.  A detailed picture  on  the
left can be encoded by a schematic picture on the right.
}
\end{figure}

\begin{Example}
\label{ex:g2}
The ribbon  graph presented at Figure~\ref{pic:g2} corresponds to
the horizontal foliation  of an Abelian differential on a surface
of genus $g=2$.  The Abelian differential  has a single  zero  of
order 2. The ribbon graph has two pairs of boundary components.
\end{Example}

Any  separatrix   diagram   represents   an  orientable  measured
foliation with only closed  leaves  on a compact oriented surface
without boundary.  We say that  a diagram is {\it realizable} if,
moreover, this measured foliation can be chosen as the horizontal
foliation       of       some        Abelian        differential.
Lemma~\ref{lm:realiz:diag}   below    gives   a   criterion    of
realizability of a diagram.

Assign to each saddle connection a real variable standing for its
``length''. Now any boundary component is  also  endowed  with  a
``length'' obtained as sum of the ``lengths'' of all those saddle
connections which  belong to this component.  If we want  to glue
flat cylinders to the boundary components,  the  lengths  of  the
components in every pair should match each other.  Thus for every
two boundary components paired together  (i.e.  having  the  same
color) we  get a linear  equation: ``the length of the positively
oriented component equals  the  length of the negatively oriented
one''.

\begin{Lemma}
\label{lm:realiz:diag}
A diagram  is realizable if  and only if the corresponding system
of linear equations  on  ``lengths'' of saddle connections admits
strictly positive solution.
\end{Lemma}

The proof is obvious.

\begin{Example}
The diagram  presented  at  Figure~\ref{pic:g2}  has three saddle
connections, all of  them are loops. Let $p_{16}, p_{52}, p_{34}$
be their ``lengths''. There are two pairs of boundary components.
The corresponding system of linear equations is as follows:
$$
\left\{\begin{array}{l}
p_{34}=p_{16}\\
p_{16}+p_{52}=p_{34}+p_{52}
\end{array}\right.
$$
\end{Example}


Here  is  a  simple  but  important  result  which  together with
Lemma~\ref{lm:closedleaves}  shows  that  one  can  encode   (not
uniquely) connected components of strata by realizable separatrix
diagrams.

\begin{Lemma}
\label{lm:coinc:diag}
Let   the   horizontal  foliations   of   Abelian   differentials
$\omega_1,\omega_2$ have only closed leaves. If the corresponding
separatrix   diagrams   are   isomorphic,   then   both   Abelian
differentials belong to  the same connected component of the same
stratum of Abelian differentials.
\end{Lemma}
\begin{proof}
In  this  context  it  is  convenient  to  think  of  an  Abelian
differential as of a  flat  surface with cone type singularities,
with trivial holonomy and with a choice of a covariantly constant
horizontal direction.

A family of Abelian differentials sharing  the  same  diagram  is
parametrized  by  the  collection  of  ``horizontal''  parameters
representing the lengths of edges of the graph (i.e., the lengths
of saddle  connections)  and  by  the  collection of ``vertical''
parameters: heights of the cylinders,  and  twists  used to paste
them  in.  The   vertical   and  the  horizontal  parameters  are
independent. There are no constraints on vertical parameters: the
heights of  the  cylinders  are  arbitrary  positive numbers; the
twists are arbitrary angles. The horizontal  parameters belong to
a  simplicial  cone: they  are  presented  by  strictly  positive
solutions of a  system  of homogeneous linear equations described
in Lemma~\ref{lm:realiz:diag}.  Thus  the  space of parameters is
connected.
\end{proof}

\begin{Lemma}
\label{lm:arrows}
Diagram of  the  horizontal  foliation  of  Abelian  differential
$-\omega$  is  obtained  from  the  diagram   of  the  horizontal
foliation of  Abelian  differential  $\omega$  by  reversing  the
arrows (orientations of edges).
\end{Lemma}

The proof is obvious.

As a corollary, we obtain  a  necessary  and sufficient condition
for  diagrams   to   represent   a  {\it  hyperelliptic}  Abelian
differential from $\H(2g-2)$. First of  all,  such  a diagram has
one only vertex  of valence $4g-2$. Consider a small neighborhood
of the vertex  of  such graph; it is  represented  by $4g-2$ rays
joined at the vertex which are organized in a cyclic order. There
is a natural (local) involution  of  this  neighborhood, the {\it
central symmetry},  which fixes the  vertex and sends each ray to
the opposite one.

\begin{Lemma}
\label{lm:hypersymmetry}
For any diagram  with one vertex corresponding to a hyperelliptic
Abelian differential $(C,\omega)\in\H(2g-2)$ the central symmetry
extends to an involution of  the  ribbon  graph interchanging any
two paired boundary components.  Also  the number of cylinders in
the diagram is equal to one plus the number of two-element orbits
of  the  involution  on  the  set  of  the  edges  of  the  graph
(separatrix loops). Conversely,  any  diagram with one vertex and
properties as above is realizable and  represents a hyperelliptic
Abelian differential.
\end{Lemma}
\begin{proof}
Hyperelliptic involution  acts  as  a  central  symmetry near the
unique  zero  of   $\omega$,   also  it  transforms  $\omega$  to
$-\omega$. This implies  the symmetry of the graph underlying the
diagram.  Also  it  shows  that  the  decomposition  of  boundary
components into pairs is also invariant under the involution. Let
us prove that the  involution  preserves each pair. Suppose there
is a pair  of distinct cylinders  which are interchanged  by  the
involution. Change slightly the  ``height''  of one of them. This
corresponds  to   a   continuous   deformation  of  the  vertical
foliation, which  leaves  the horizontal foliation unchanged. The
deformed  Abelian  differential  is  supposed  to   stay  in  the
component    $\mathcal{H}^{hyp}(2g-2)$    which   leads    to   a
contradiction, since the involution does not exist anymore.

Let us establish now the  numerical  property  $n_c=n_2+1$  where
$n_c$ is the  number  of  cylinders and $n_2$ is  the  number  of
two-element orbits as in Lemma.  The  set of fixed points of  the
involution  consists of  a)  the vertex of  the  diagram, b)  the
middle point on every  involutive  separatrix loop, c) two points
in the interior of each cylinder. The total  number of separatrix
loops is equal  to $2g-1$, therefore  the number $n_1$  of  loops
invariant under the involution  is  equal to $2g-1-2 n_2$. Hence,
we have  $1+n_1+2n_c=2g+2(n_c-n_2)$  fixed  points.  On the other
hand, the number of fixed points of a hyperelliptic involution is
equal to $2g+2$ which implies that $n_c=n_2+1$.

Conversely, for a diagram with the properties listed in the Lemma
the  realizability  is  obvious  because  we  can assign to  each
separatrix  loop  the  same  length,  which  gives us a  positive
solution   of    the    system    of    linear   equations   from
Lemma~\ref{lm:realiz:diag}.  The  corresponding  surface  carries
canonically an involution with $2g+2$ fixed  points, therefore by
Hurwitz formula the quotient surface has genus zero and we are in
the hyperelliptic case.
\end{proof}

\begin{Remark}
\label{rm:arrows}
Consider a realizable separatrix diagram corresponding  to a {\it
connected} closed  surface,  and  forget  the  orientation of the
edges.  There are  exactly  two ways to  orient  again our  graph
(keeping the initial structure  of  the ribbon graph, and keeping
the initial distribution  of  the boundary components into pairs)
which lead to  a  realizable diagram:  the  initial way, and  the
opposite one. This is true  even  if the underlying graph of  the
diagram  is  not connected.  According  to  Lemma~\ref{lm:arrows}
these  two  orientations  correspond  to  Abelian   differentials
$\omega$ and  $-\omega$. Note that Abelian differentials $\omega$
and $-\omega$  belong  to  the same stratum $\cH(k_1,\dots,k_n)$;
moreover, they belong  to the same connected component since they
can  be  joined   inside  the  stratum  by  the  continuous  path
$e^{i\theta}\omega$, $\theta\in [0;\pi]$. Thus,  it  follows from
Lemmas~\ref{lm:coinc:diag}    and~\ref{lm:arrows}    that    both
orientations of  a  {\it  realizable}  diagram  represent Abelian
differentials  from  the  same  connected component of  the  same
stratum. Hence, if we care only  about  connected  components  of
strata then in pictures of separatrix diagrams we can omit arrows
(directions of edges).
\end{Remark}

\subsection{Bubbling handles}
\label{ss:forget}

In  this  section  we  describe  a  local  surgery  (``bubbling a
handle'') which modifies the surface in a small neighborhood of a
chosen  zero  of  the  Abelian  differential.  Here  it  will  be
convenient to use ``numbered'' versions  of  moduli  spaces  (see
Remark \ref{rm:number} from Introduction). Also here and later in
Section \ref{s:conn:comp} in order to alleviate notations we will
denote  a  point  $[(C,\omega)]$  of the moduli space  simply  by
$\omega$   and   will   write  slightly  incorrectly   $\omega\in
\H^{num}(k_1,\dots,k_n)$.

Topologically the surgery corresponds  to  adding a handle to the
surface. Metrically we choose a small disk centered at the chosen
conical singularity, then we  make  some geodesic cuts inside the
disk and paste in a small metric cylinder. Having started with an
Abelian differential
$\omega\in\mathcal{H}^{num}(k_1,\dots,k_{i-1},k_i,k_{i+1},\dots,k_n)$
we construct an Abelian differential
$\hat{\omega}\in\mathcal{H}^{num}(k_1,\dots,k_{i-1},k_i+2,k_{i+1},\dots,k_n)$,
where the  surface was modified in  the neighborhood of  the zero
$P_i$ of multiplicity  $k_i$. The surgery depends on one discrete
and on two complex parameters.

Revising this paper we decided to replace the  initial version of
the  surgery,  by   the   more  general  one  described  recently
in~\cite{Eskin:Masur:Zorich},  called there  the  ``figure  eight
construction''. Here  we present briefly this latter construction
consisting of two steps.

{\bf Breaking up a zero.} We  first describe how one can break up
a zero $P_i$ of multiplicity $k$ of an  Abelian differential into
two zeroes of multiplicities $k',  k''$,  where  $k'+k''=k$, by a
local surgery. In fact, we will need this construction also  in a
slightly more general case when parameter $k''$ is equal to zero.

Consider a metric disk  of  a small radius $\varepsilon$ centered
at the  point $P_i$, i.e. the set  of points  $Q$ of the  surface
such that Euclidean distance from $Q$ to the point $P_i$  is less
than or equal  to  $\varepsilon$. We suppose that $\varepsilon>0$
is chosen small enough,  so  that the $\varepsilon$-disk does not
contain any  other conical points  of the metric; we assume also,
that  the   disk  which  we   defined  in  the  metric  sense  is
homeomorphic to a topological disk. Then, metrically our disk has
a structure  of a regular cone  with a cone  angle $2\pi(k_i+1)$;
here $k_i$ is  the  multiplicity of the zero  $P_i$.  Now cut the
chosen disk  (cone) out of the surface. We  shall modify the flat
metric inside it preserving the metric at the  boundary, and then
paste the modified disk (cone) back into the surface.

\begin{figure}[ht]
%
\includegraphics{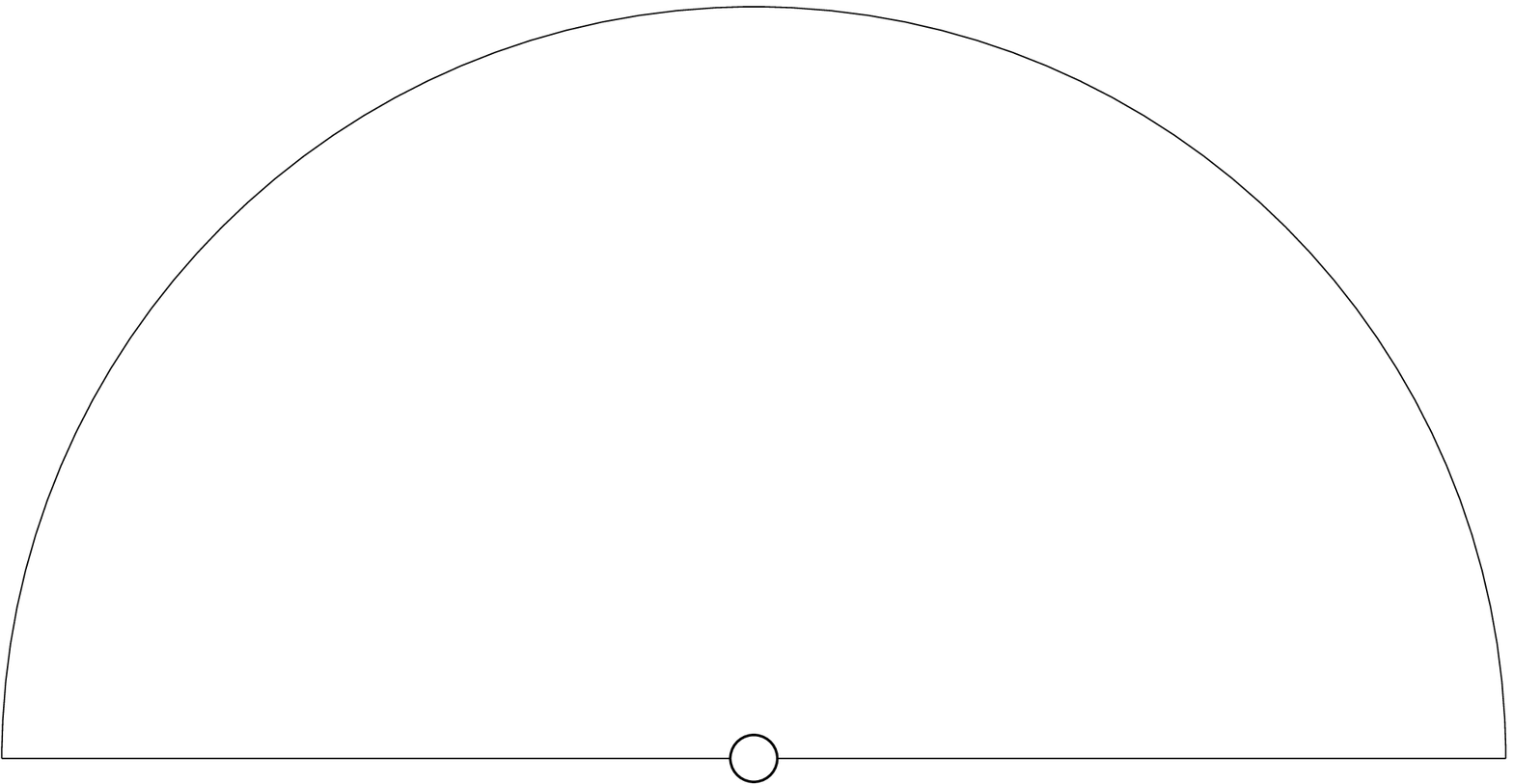}
\includegraphics{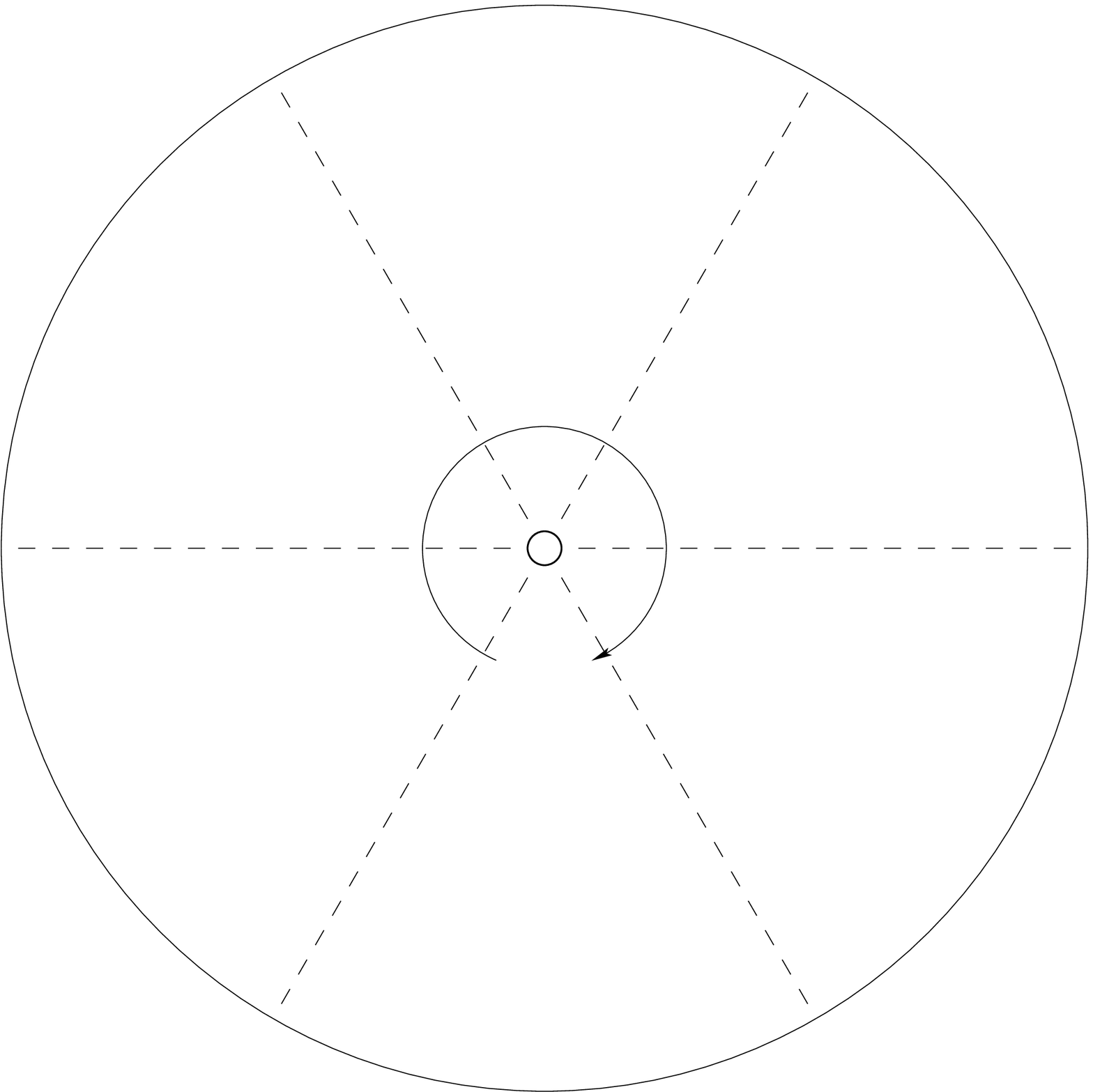}
\includegraphics{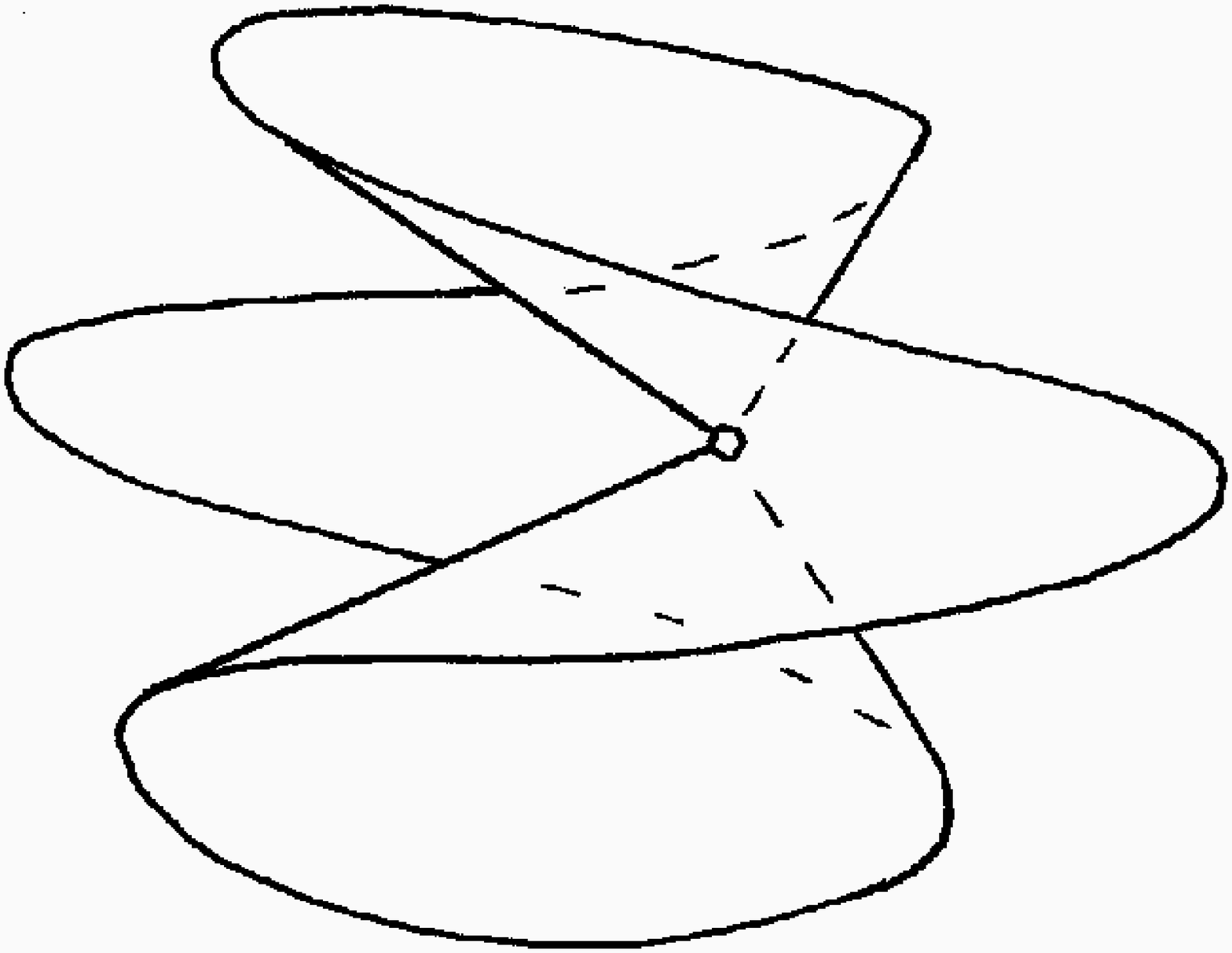}
%
%
\includegraphics{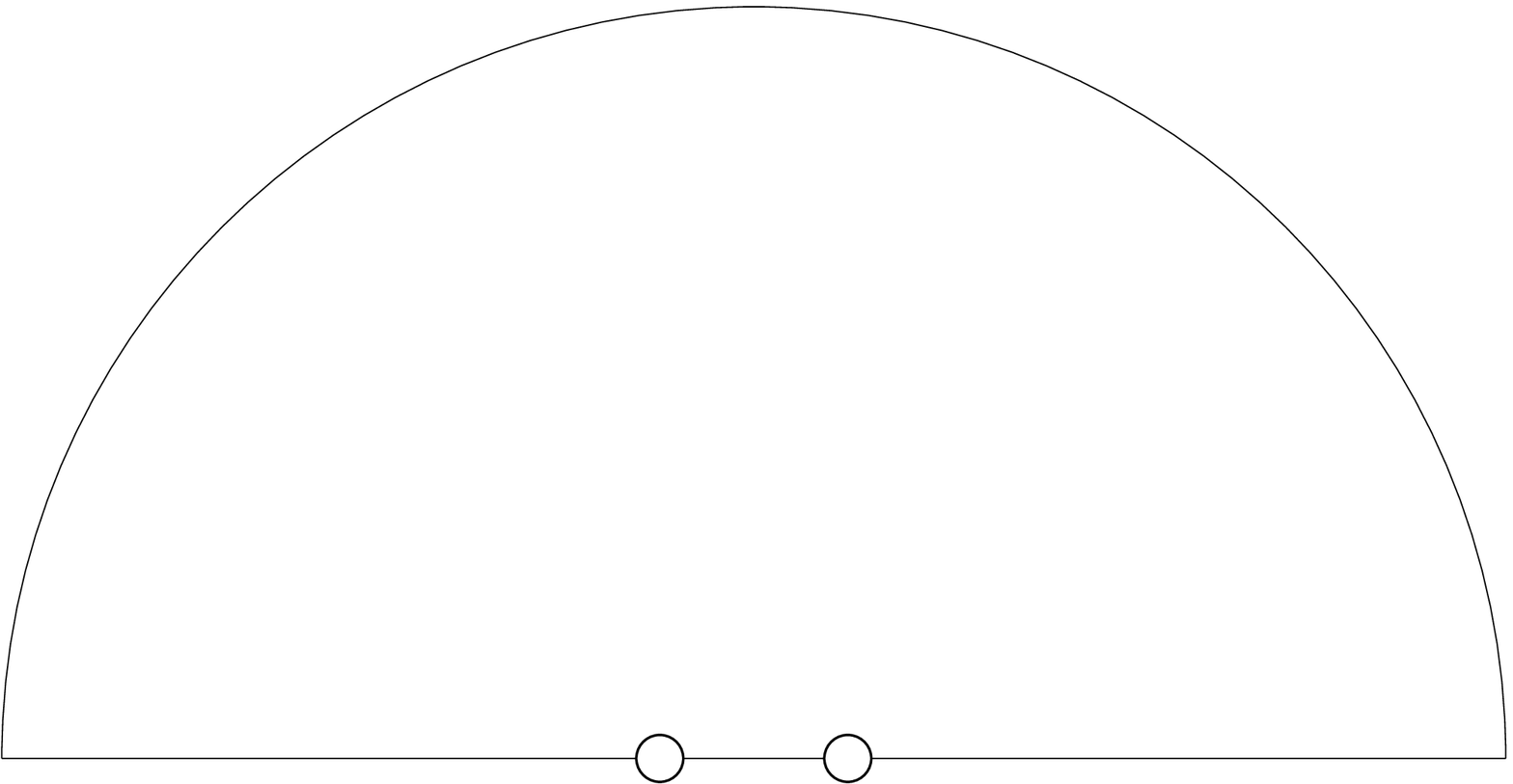}
\includegraphics{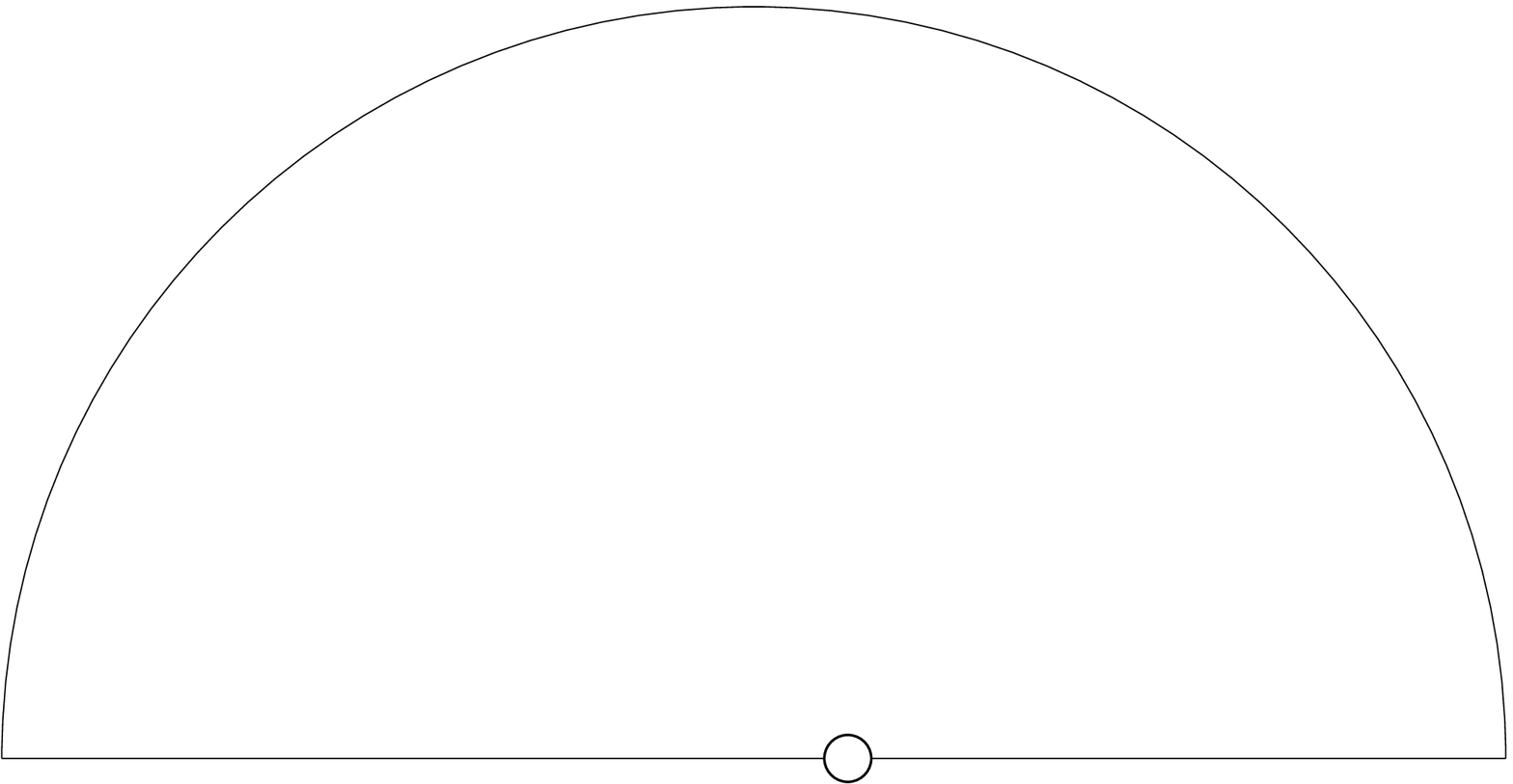}
\includegraphics{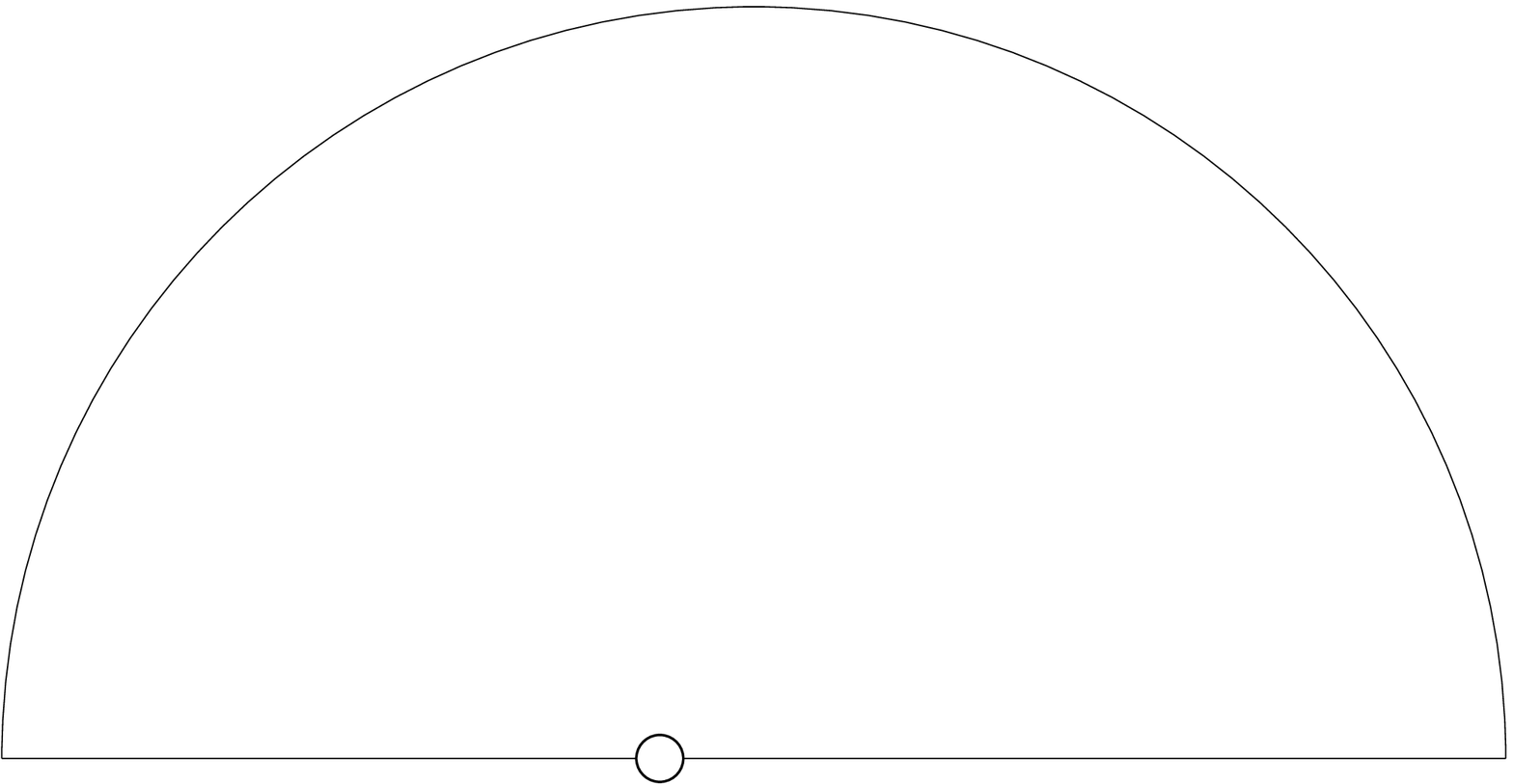}
\includegraphics{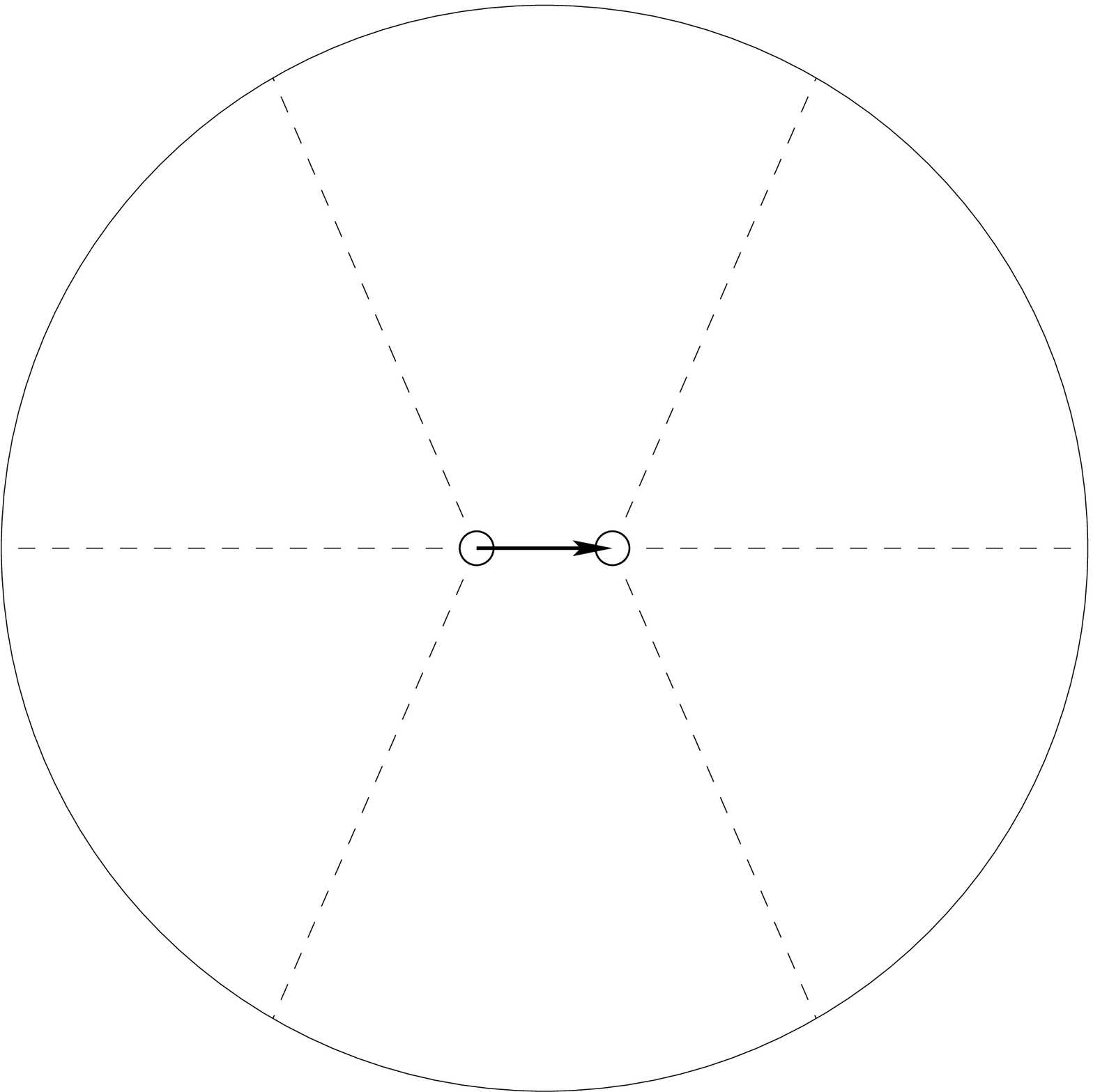}
\includegraphics{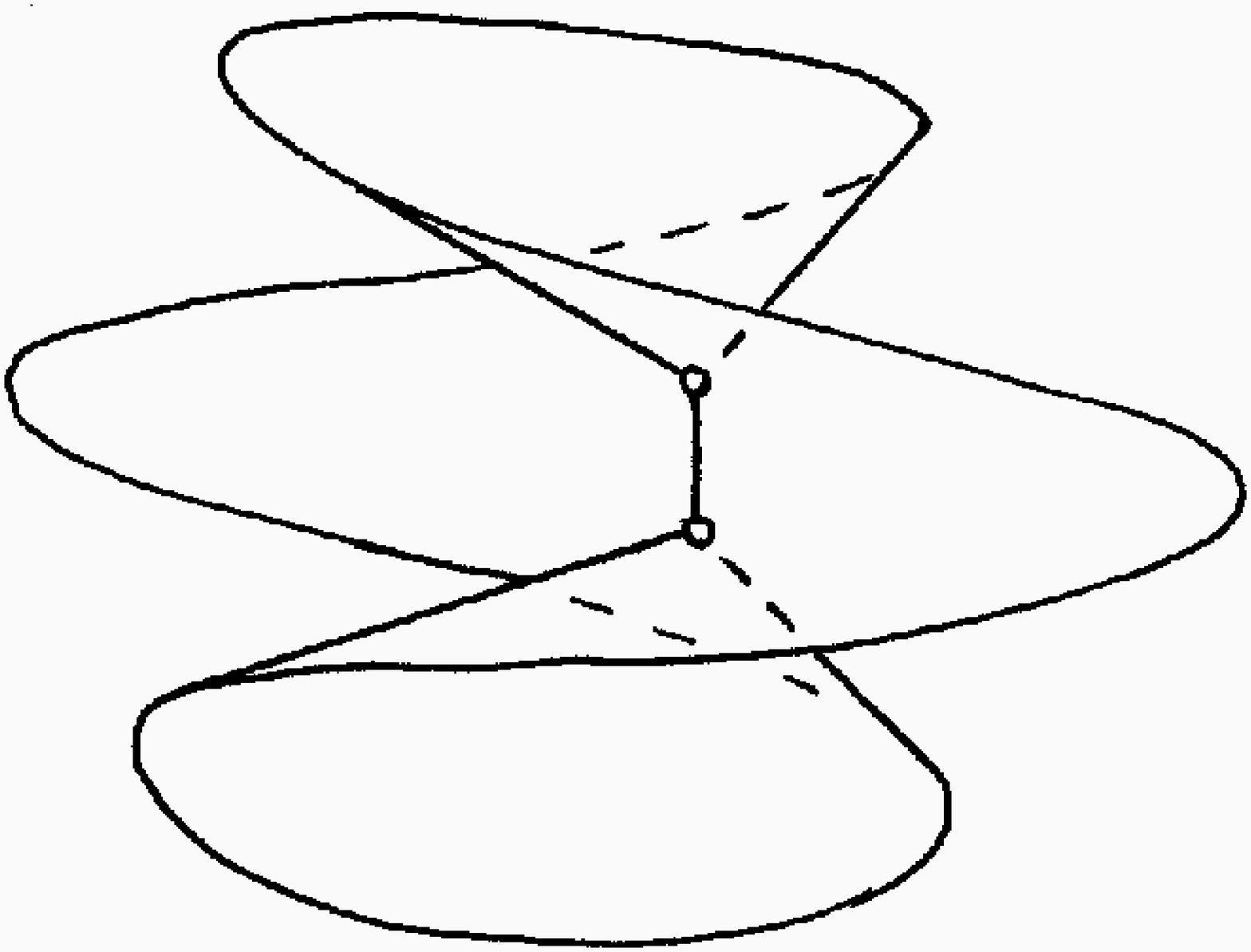}
%
%
\begin{picture}(0,0)(0,0)
\put(10,-10)
 {\begin{picture}(0,0)(0,0)
 \put(-154,-55){$\scriptstyle \varepsilon$}
 \put(-125,-55){$\scriptstyle \varepsilon$}
 \put(-38,-66){$\scriptstyle 6\pi$}
 \end{picture}}
\put(10,-30)
 {\begin{picture}(0,0)(0,0)
 \put(-145,-148){$\scriptstyle 2\delta$}
 \put(-160,-197){$\scriptstyle \varepsilon+\delta$}
 \put(-130,-197){$\scriptstyle \varepsilon-\delta$}
 \put(-160,-247){$\scriptstyle \varepsilon-\delta$}
 \put(-130,-247){$\scriptstyle \varepsilon+\delta$}
 \put(-39,-197){$\scriptstyle 2\delta$}
 \put(-65,-197){$\scriptstyle \varepsilon+\delta$}
 \put(-19,-197){$\scriptstyle \varepsilon+\delta$}
 \put(-64,-173){$\scriptstyle \varepsilon-\delta$}
 \put(-19,-173){$\scriptstyle \varepsilon-\delta$}
 \put(-65,-217){$\scriptstyle \varepsilon-\delta$}
 \put(-17,-217){$\scriptstyle \varepsilon-\delta$}
\end{picture}}
\end{picture}
\vspace{280bp} 
\caption{
\label{pic:breaking:up:a:zero} Breaking up a zero into two zeroes
(after~\cite{Eskin:Masur:Zorich}).}
\end{figure}

Our cone  can be glued  from $2(k_i+1)$ copies of standard metric
half-disks of the radius $\varepsilon$,  see  the  picture at the
top  of Figure~\ref{pic:breaking:up:a:zero}.  Choose  some  small
$\delta$, where  $0<\delta<\varepsilon$  and  change  the  way of
gluing  the  half-disks  as  indicated on the bottom  picture  of
Figure~\ref{pic:breaking:up:a:zero}. As patterns we still use the
standard  metric  half-disks,  but  we move slightly  the  marked
points on their diameters. Now we  use  two  special  half-disks;
they  have two  marked  points on the  diameter  at the  distance
$\delta$ from the center of the half disk. Each of  the remaining
$2k_i$  half-disks  has  a  single marked point at  the  distance
$\delta$ from the center of the half-disk. We are alternating the
half-disks with the marked point  moved  to the right and to  the
left from the  center.  The picture  shows  that all the  lengths
along identifications are matching; gluing the half-disks in this
latter way  we obtain a topological disk with  a flat metric; now
the flat metric has two cone-type  singularities  with  the  cone
angles $2\pi(k'+1)$  and  $2\pi(k''+1)$,  where $k'+k''=k_i$, and
$k', k''\in \Z_+$.  By convention we denote the multiplicities of
the newborn  zeroes in such way that $k'\ge  k''$. Here $2k'$ and
$2k''$ are the numbers of half-disks with one  marked point glued
in between the distinguished  pair  of half-disks with two marked
points.

By  technical  reasons  it  would be convenient to  include  into
consideration the  trivial case, when  $k''$ is equal to zero. In
this latter  case we, actually, do not change  the metric at all;
we just mark  a  point $P''$ at the  distance  $2\delta$ from the
point $P_i=P'$.

Note that a  small tubular neighborhood  of the boundary  of  the
initial  cone   is   isometric   to   the  corresponding  tubular
neighborhood of the boundary of the resulting object. Thus we can
paste it back into the surface. Pasting it back we can turn it by
any angle $\varphi$, where $0\le\varphi< 2\pi(k_i+1)$.

We described  how to break up a  zero of  multiplicity $k$ of  an
Abelian differential  into two zeroes of multiplicities $k',k''$,
where $k'+k''=k$, and  $k'\ge k''$. The construction is local; it
is parameterized by the two  free  real  parameters (actually, by
one complex parameter):  by  the small distance $2\delta$ between
the newborn  zeroes, and by  the direction $\varphi$ of the short
geodesic segment joining  the  two newborn zeroes. In particular,
as  a parameter  space  for this construction  one  can choose  a
punctured disk.

Now we can proceed with the second step of the construction.

\begin{figure}[hbt]
%
%
%
\includegraphics{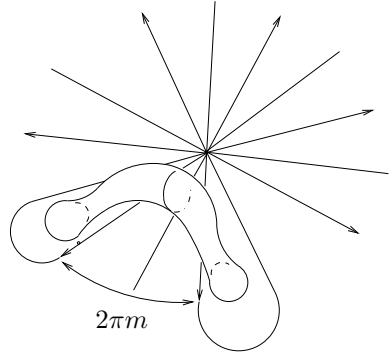}
\vspace{130bp}
\begin{picture}(0,0)(0,0)
\put(0,0) {\begin{picture}(0,0)(0,0)
\put(-52,6){$2\pi m$}
\end{picture}}
\end{picture}
\caption{
\label{pic:hand}
Two simple separatrix loops determining a handle.
}
\end{figure}

{\bf Bubbling a handle into a slit.}
Let us  slit the surface along the short  geodesic segment of the
length $2\delta$ joining the newborn zeroes $P', P''$  and let us
identify the endpoints of the slit. The resulting surface has two
boundary components joined together at  the  point  $P'=P''$.  By
construction  the  boundary components are geodesics in the  flat
metric  determined  by  $\omega$;  they  have   the  same  length
$2\delta$.  Take a  small  flat cylinder with  a  waist curve  of
length  $2\delta$  and  paste  it into our surface.  The  surface
$M^2_{g+1}$ is  constructed.  The  flat  structure on $M^2_{g+1}$
together with the choice of  the  horizontal  direction  uniquely
determine an Abelian differential $\hat\omega$ on $M^2_{g+1}$. By
construction  the  resulting  Abelian  differential  $\hat\omega$
belongs               to               the                stratum
$\mathcal{H}^{num}(k_1,\dots,k_i+2,\dots,k_n)$, where $\omega \in
\mathcal{H}^{num}(k_1,\dots,k_i,\dots,k_n)$,  and  $k_i$  is  the
multiplicity of the zero $P_i$ (the case $k_i=0$ is not excluded;
in  this  case  $P_i$  is  just  a  marked  point).  The  Abelian
differential $\hat\omega$ is obtained from $\omega$ by ``bubbling
a small handle'' at the zero $P_i$ (see Figure~\ref{pic:hand}).

This surgery is parameterized by the following list of
parameters:\newline
--- Discrete parameter $k'$, where  $k_i/2\le  k'\le  k_i$.  This
parameter  indicates   the   number   of   sectors   between  the
distinguished pair of sectors.  There  are $2k'+2$ sectors on the
one   side   and   $2k''+2$   on   the   other  side;  see   also
Figure~\ref{pic:hand}, where $m$ denotes $(k''+1)$;\newline
---  Pair  of   free   real  parameters  $\delta$  and  $\varphi$
responsible  for  the  breaking  up  a  zero;  in  the  resulting
construction they represent the length of the waist  curve of the
cylinder and direction in which  goes  the  corresponding  closed
geodesic. This pair of real parameters can be seen as one complex
parameter: the  period  of  the Abelian differential $\hat\omega$
along the waist curve of the new cylinder;\newline
--- Finally, we have  two  more free real parameters representing
the height  of the cylinder, and the twist  which we used pasting
it into the surface. They can be organized in a complex parameter
representing the period of $\hat\omega$ along the cycle following
the new handle.

Let us describe now the properties of this surgery.

\begin{Lemma}
\label{lm:rotate:handle:general}
Consider  Abelian  differentials $\hat\omega_1,  \hat\omega_2 \in
\mathcal{H}^{num}(k_1,\dots,k_i+2,\dots,k_n)$     obtained     by
``bubbling  a  handle'' at  the  same zero  $P_i$  of an  Abelian
differential                                           $\omega\in
\mathcal{H}^{num}(k_1,\dots,k_i,\dots,k_n)$. If the  angle  $2\pi
m$  between  the   sectors   corresponding  to  the  handle  (see
Figure~\ref{pic:hand}) is  the  same  for  $\hat\omega_1$ and for
$\hat\omega_2$,  then   there   exist   a   continuous   path  in
$\mathcal{H}^{num}(k_1,\dots,k_j+2,\dots,k_n)$            joining
$\hat\omega_1$ with $\hat\omega_2$; in particular, $\hat\omega_1,
\hat\omega_2$  belong  to  the  same connected component  of  the
stratum $\mathcal{H}^{num}(k_1,\dots,k_i+2,\dots,k_n)$.
\end{Lemma}
\begin{proof}
Fixing the  discrete parameter $k'$ (which  is equal to  $m-1$ if
$m>k_i/2-1$, and to $k_i-m+1$ otherwise) we describe ``bubbling a
handle'' at the zero $P_i$ of $\omega$ by two pairs of continuous
parameters   as   above.  Note  that  the  space  of   parameters
$\delta,\varphi$ describing breaking up the  zero  into  two is a
punctured  disk;  the  space  of  parameters  parameterizing  the
cylinder (the height of the cylinder and the twist used  to paste
it into the  surface)  is also a punctured  disk.  Thus the total
space of  parameters is a  direct product of two punctured disks,
which is obviously path-connected.
\end{proof}

\begin{Lemma}
\label{lm:lift:path}
Let       an       Abelian      differential       $\hat\omega\in
\mathcal{H}^{num}(k_1,\dots,k_i+2,\dots,k_n)$ be obtained from an
Abelian                                              differential
$\omega\in\mathcal{H}^{num}(k_1,\dots,k_i,\dots,k_n)$          by
``bubbling a  handle'' at some  zero $P_i$. Any path $\rho: [0;1]
\to \mathcal{H}^{num}(k_1,\dots,k_i,\dots,k_n)$ which  starts  at
$\omega$  can   be   lifted   to   a   path  $\hat\rho:  [0;1]\to
\mathcal{H}^{num}(k_1,\dots,k_i+2,\dots,k_n)$     starting     at
$\hat\omega$ by continuous ``bubbling a handle'' along $\rho$.
\end{Lemma}
\begin{proof}
Note that we  can  bubble a small handle  into  {\it any} Abelian
differential  of  the  given  stratum. This implies that  we  can
choose     an     appropriate     subset    in    the     stratum
$\mathcal{H}^{num}(k_1,\dots,k_i+2,\dots,k_n)$  which  would have
the  structure  of  a  (singular) fiber bundle over  the  stratum
$\mathcal{H}^{num}(k_1,\dots,k_i,\dots,k_n)$.  The  regular fiber
is  a  direct product of two  disks  punctured at the centers;  a
singular fiber  is a quotient of the direct  product of two disks
punctured at the centers by a finite group of symmetry.  Thus the
lifting described in the Lemma is the particular  case of lifting
of a path from the  base  to  the  total space of a fiber bundle.
\end{proof}

\begin{Remark}
\label{rm:forget:recall}
In other  words Lemma~\ref{lm:lift:path} means the following. Let
an Abelian  differential  $\hat\omega$  on  a  Riemann surface of
genus $g+1$ be constructed from an  Abelian differential $\omega$
on  a  surface of  genus  $g$  by  ``bubbling  a  small handle''.
Morally,   we   can   temporarily   ``forget''   (or   ``erase'')
corresponding handle;  modify  Abelian differential $\omega$ on a
surface of genus $g$ inside his stratum in a continuous  way, and
then  ``recall''  the  ``forgotten''  handle  for  the  resulting
Abelian   differential   $\omega'$   provided  that  the   metric
parameters of the handle are  sufficiently  small.  Then  Abelian
differential $\hat\omega'$ on a  surface  of genus $g+1$ which we
obtain  as a  result  of this construction  belongs  to the  same
connected component as Abelian differential $\hat\omega$.
\end{Remark}

Suppose  that  we   ``bubble  a  handle''  into  a  flat  surface
corresponding to  an  Abelian differential $\omega$ having zeroes
of  even  multiplicities.  The  resulting  Abelian   differential
$\hat\omega$  would  have the same property. In particular,  both
Abelian differentials determine spin structures on  corresponding
surfaces.   The   following   Lemma    compares    the   parities
$\varphi(\omega)$   and   $\varphi(\hat\omega)$  of   these  spin
structures.

\begin{Lemma}
\label{lm:change:of:parity}
Let       an       Abelian      differential       $\hat\omega\in
\mathcal{H}^{num}(2(l_1+1),2l_2,...,2l_n)$ on a  surface of genus
$g+1$  be   obtained  from  an  Abelian  differential  $\omega\in
\mathcal{H}^{num}(2l_1,2l_2,...,2l_n)$ on a surface of genus  $g$
by ``bubbling a handle''. Let $2\pi m$ be the angle of one of the
two    sectors     complementary     to     the    handle,    see
Figure~\ref{pic:hand}.  The  parities  of  the  spin   structures
determined by $\omega$ and by $\hat\omega$  are  related  in  the
following way:
$$ \varphi(\hat\omega) - \varphi(\omega) = m+1(\mod 2) $$
\end{Lemma}
\begin{proof}
Choose    a    collection    of     oriented     simple    curves
$(\alpha_i,\beta_i)_{i=1,\dots,g}$  on  the  initial  surface  of
genus $g$ representing a  symplectic  basis in the first homology
group.  Deforming,  if necessary, the paths inside their  isotopy
classes we  can make them  stay away from some small neighborhood
$U(P_1)$ of the zero $P_1$ under  consideration.  We  can  assume
that the surgery (``bubbling a  handle'')  was  performed  inside
this small neighborhood  $U$. In particular, the surgery does not
affect the initial collection of paths.

The initial basis  of  cycles can  be  completed to a  symplectic
basis on  the surface of  genus $g+1$ obtained after ``bubbling a
handle'' by adding two additional curves on the handle created in
the surgery.

One of  these new curves,  which we denote by $\alpha_{g+1}$, can
be represented  by a waist curve of the  cylinder which we pasted
in; moreover, we  can  choose  this waist curve to  be  a  closed
geodesic. By Lemma~\ref{lm:ind0} we get $ind_{\alpha_{g+1}}=0$.

The  second  curve, $\beta_{g+1}$,  can  be  chosen  as  follows.
Suppose for simplicity that the handle was bubbled  with the slit
made in the horizontal direction, and with a trivial twist, i.e.,
the leaf  of the vertical  foliation emitted from zero $P_1$ into
the new handle returns back to $P_1$. (Changing the twist we stay
in   the    same    connected    component    of    the   stratum
$\mathcal{H}^{num}(2(l_1+1),2l_2,\dots,2l_n)$; in particular,  we
do not change  the  parity of the spin  structure  defined by the
corresponding $\hat\omega$.) Take  a  circle centered at $P_1$ in
the initial surface. Making the radius of the circle sufficiently
small  we  get   the  circle  inside  the  neighborhood  $U$;  in
particular,  it  does  not  intersect any of the  initial  curves
$\alpha_i, \beta_i$, where $i\le g$. Choose an arc of this circle
joining  two  distinguished sectors  (see Figure~\ref{pic:hand}).
Then the endpoints of the arc might be joint by  a  piece of leaf
of the vertical foliation along  the  new  cylinder (recall, that
the  twist  of  the  handle   is   by   assumption  trivial).  By
construction the  resulting  closed path $\beta_{g+1}$ is smooth;
it  does  not  intersect  any  of   the   initial   curves;   and
$[\alpha_{g+1}]\cdot [\beta_{g+1}]=1$. Hence we  got  the desired
symplectic basis of cycles.

Direct calculation gives $ind_{\beta_{g+1}}=m$, since the tangent
vector to the  path  $\beta_{g+1}$ turns  by  the angle $2\pi  m$
along  the  arc,  and does not  turn  at  all  along the vertical
segment. Now  everything is ready  to compare the parities of the
spin  structures  of  $\omega$  and  $\hat\omega$;  here  we  use
formula~\ref{eq:parity:of:spin:structure}
\begin{multline*}
\varphi(\omega')=\sum_{i=1}^{g+1} (ind_{\alpha_i}+1)\cdot(ind_{\beta_i}+1) =\\
=\sum_{i=1}^{g} (ind_{\alpha_i}+1)\cdot(ind_{\beta_i}+1) +
(ind_{\alpha_{g+1}}+1)\cdot(ind_{\beta_{g+1}}+1) =\\
= \varphi(\omega) + (0+1)\cdot(m+1)
\end{multline*}
\end{proof}

Consider now ``bubbling a handle''  in  the  particular case when
the horizontal  foliation  of  the  initial  Abelian differential
$\omega$ has only closed leaves. If at the  intermediate stage of
``bubbling a handle''  we  break up  the  zero in the  horizontal
direction, then the horizontal foliation of the resulting Abelian
differential $\hat\omega$  obtained  after  ``bubbling a handle''
also has only  closed  leaves. It would be  convenient  for us to
reformulate the Lemmas above in this particular case  in terms of
the  diagrams   of   the   Abelian   differentials  $\omega$  and
$\hat\omega$.

\begin{Example}
\label{ex:bubble:into:torus}
Consider a flat torus and chose the horizontal direction on it in
such way that  the leaves of  the horizontal foliation  would  be
closed. ``Bubbling a handle'' in the horizontal direction at some
point of the torus we get  a surface of genus $2$ with horizontal
foliation  having  only  closed  leaves.  The   diagram  of  this
foliation is presented at Figure~\ref{pic:g2}. \end{Example}

We  call  a  separatrix  loop {\it simple} if  the  corresponding
outgoing  and   ingoing   separatrix   rays  are  neighbors  (see
Figure~\ref{pic:hand}). In  terms of diagrams of separatrix loops
``bubbling  a  handle'' in a horizontal direction corresponds  to
adding a pair of simple separatrix loops of the same color to the
initial diagram, see Example~\ref{ex:bubble:into:torus} above.

Note that when we ``bubble a  handle'' at a true zero of $\omega$
(and not  at a  regular point as in the  Example above) there are
several horizontal directions at a conical  point. In particular,
even when  we fix the discrete  parameter $k'$ there  are several
ways  to  ``bubble a handle'' in the  horizontal  direction.  The
first Lemma  says that fixing  the discrete parameter $k'$ we get
to  the  same  connected  component,  no  matter  which  of these
horizontal direction we choose.

\begin{Lemma}
\label{lm:rotate:handle:diagram}
Rotating a pair  of  simple separatrix  loops  of the same  color
corresponding to the handle just ``bubbled'' in such way that the
number of sectors between the pair  of  simple  separatrix  loops
stay   unchanged   we  obtain   diagrams   representing   Abelian
differentials  from   the   same   connected   component  of  the
corresponding stratum.
\end{Lemma}
\begin{proof}
The proof is a direct corollary of Lemma~\ref{lm:coinc:diag} and
Lemma~\ref{lm:rotate:handle:general}.
\end{proof}

\begin{Example}
The   realizable   diagrams  presented   at  Figure~\ref{pic:hg4}
correspond to  Abelian  differentials  from  the  same  connected
component.
\end{Example}

Now  let  us  reformulate Lemma~\ref{lm:change:of:parity} in  the
case when the leaves of horizontal foliations of  $\omega$ and of
$\hat\omega$ are closed. We assume  that  all  zeroes of $\omega$
and $\hat\omega$  have  even  degrees,  or  equivalently that the
corresponding separatrix  diagrams have vertices only of valences
$2(mod\ 4)$.

\begin{Lemma}
\label{lm:erase}
Let Abelian differentials $\omega$ on a surface of  genus $g$ and
$\hat\omega$  on  a  surface  of  genus   $g+1$  have  horizontal
foliations  with   only   closed   leaves.   Suppose   that   the
corresponding separatrix diagram of $\omega$ is obtained from the
diagram of $\hat\omega$ by erasing  a  pair  of simple separatrix
loops  corresponding to  the  same zero and  bounding  a pair  of
sectors of the same color (see  Figure~\ref{pic:hand}). Let $2\pi
m$  be  the angle of one  of  the two complementary sectors.  The
parities of the spin structures determined  by  $\omega$  and  by
$\hat\omega$ are related in the following way:
$$ \varphi(\hat\omega) - \varphi(\omega) = m+1(\mod 2) $$
\end{Lemma}
\begin{proof}
Bubbling an appropriate  handle in the flat surface determined by
$\omega$ we obtain  an Abelian differential with the same diagram
as $\hat\omega$. It follows  from  Lemma~\ref{lm:coinc:diag} that
all Abelian differentials corresponding to the  same diagram have
the same parity of the spin structure. Thus the Lemma is a direct
corollary of Lemma~\ref{lm:change:of:parity}.
\end{proof}

\section{Connected components of  the  strata}
\label{s:conn:comp}


\subsection{Connected   components   of   the   minimal   stratum
$\mathcal{H}(2g-2)$}
\label{ss:cc:minimal}

Here we will proceed  by induction in genus $g$ and will  use the
fact  that  any  connected  component  of   the  minimal  stratum
$\cH(2g)$  in  genus $g+1$ can be accessed  from  some  connected
component  of  the minimal stratum $\cH(2g-2)$ in  genus  $g$  by
``bubbling  a  handle''.  Surprisingly,  this  statement  is  not
trivial. Despite of many attempts we were unable to find a purely
geometric proof of the Lemma below; we use the arguments based on
combinatorial              Lemma~\ref{lm:m1adv}              from
Appendix~\ref{ss:comb:rauzy:cl}. The  difficulty which one  meets
here is as follows. In every connected component of $\cH(2g)$ one
can easily find Abelian differentials having a cylinder filled by
regular closed leaves of the horizontal  foliation. Moreover, one
can get examples  when  such  cylinder is bounded by  a  pair  of
simple separatrix loops.  However, we should warn the reader that
many  of  these Abelian differentials are {\it  not}  results  of
``bubbling a handle'', see~\cite{Eskin:Masur:Zorich} for details.

\begin{Lemma}
\label{lm:accessible:by:bubbling}
Any connected component of  the  minimal stratum $\cH(2g)$ can be
accessed from some connected component  of  the  minimal  stratum
$\cH(2g-2)$  by  ``bubbling  a  handle''  into   a  flat  surface
corresponding to an Abelian differential $\omega\in\cH(2g-2)$.
\end{Lemma}
\begin{proof}
Consider   the    extended   Rauzy   class    $\mathfrak{R}_{ex}$
corresponding  to   the   connected   component  of  the  stratum
$\mathcal{H}(2g)$ under consideration. Choose a permutation  $\pi
\in \mathfrak{R}_{ex}$ as in  Lemma~\ref{lm:m1adv}  from Appendix
A.  Consider  an Abelian differential $\hat\omega$ obtained as  a
suspension over  an  interval  exchange  transformation  with the
permutation   $\pi$,   with   integer   $\lambda_i$,   and   with
$\lambda_1=\lambda_m$.  Since  all  $\lambda_i$ are integer,  the
vertical foliation of  $\hat\omega$ has only closed leaves. It is
easy  to  see that the vertical  foliation  has a pair of  simple
separatrix loops: that is an ingoing and outgoing separatrices in
each loop are neighbors, and this pair of simple loops determines
a  handle  (cf. Figure~\ref{pic:hand});  subintervals  $I_1$  and
$I_m$ belong to this handle.

Consider  now  an  Abelian  differential $\omega$ obtained  as  a
suspension  over   an   interval   exchange  transformation  with
permutation $\pi'$  as  in  Lemma~\ref{lm:m1adv},  and  with  the
vector  of  lengths  $\lambda_2,  \dots,  \lambda_{m-1}$,   where
$\lambda_i$ are same as above. By the choice  of the permutations
$\pi$  and   $\pi'$   we   get  $\omega\in\cH(2g-2)$.  Since  all
$\lambda_i$ are integer, the vertical foliation  of $\omega$ also
has only closed leaves. Moreover,  it  is easy to check that  the
separatrix diagram of  the  vertical foliation of $\hat\omega$ is
obtained from the separatrix diagram of the vertical foliation of
$\omega$ by ``bubbling a handle''. The statement of the Lemma now
follows from Lemma~\ref{lm:coinc:diag}.
\end{proof}

From now on till the end of this subsection we will work in terms
of separatrix  diagrams  representing  Abelian differentials from
stratum $\H(2g-2)$.

\begin{figure}[htb]
%
%
\includegraphics{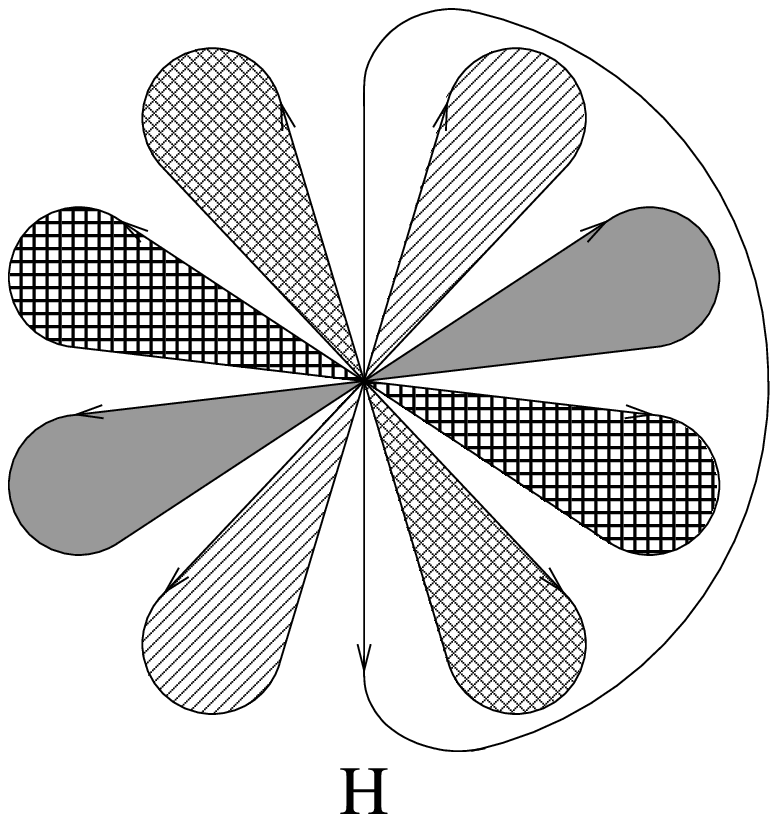}
\includegraphics{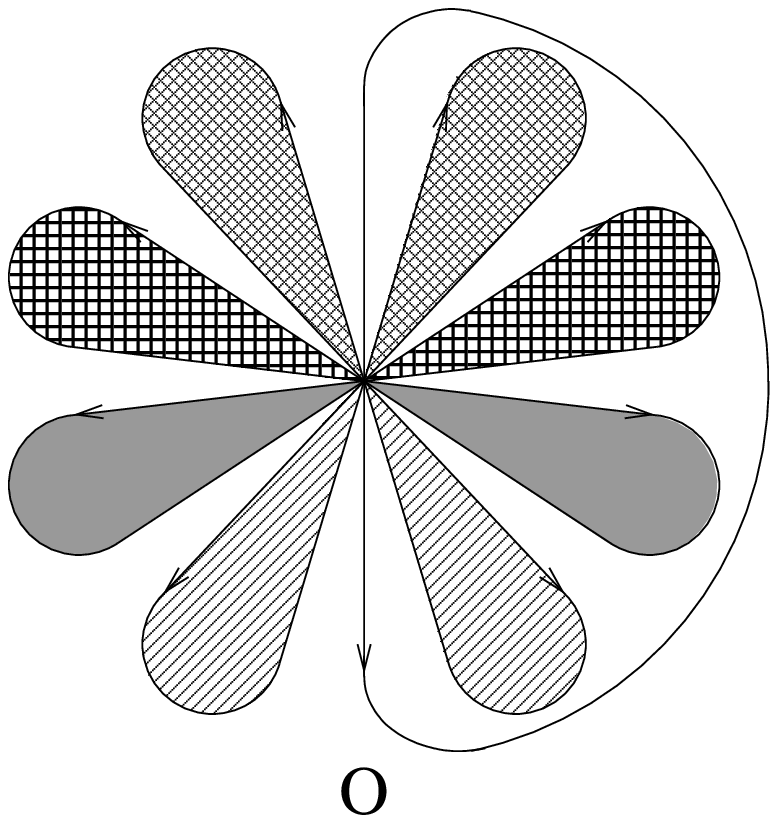}
\includegraphics{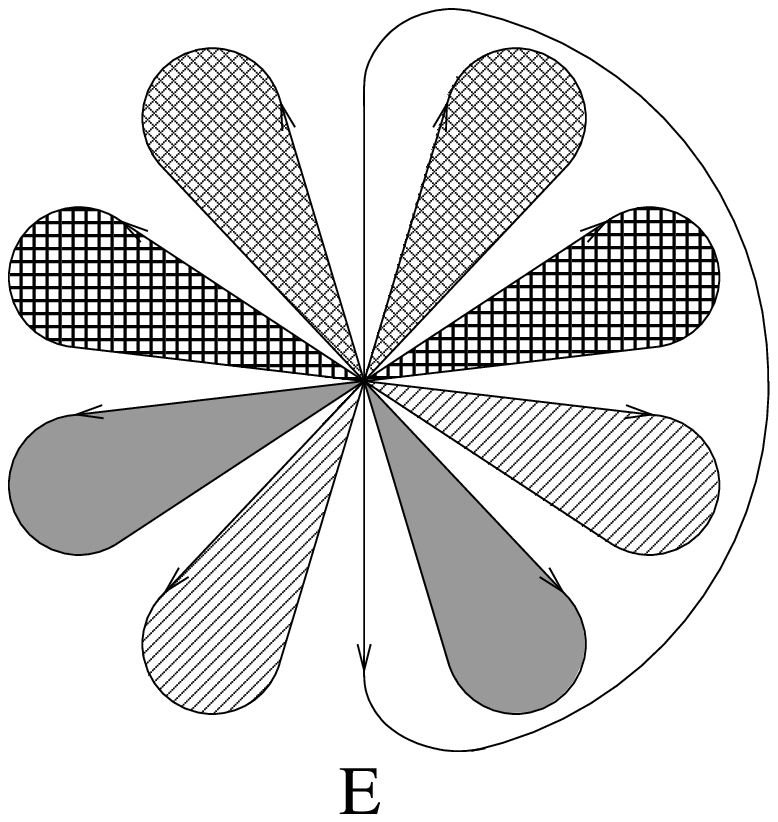}
\vspace{130bp}
\caption{
\label{pic:g5}
The diagrams represent the  following  components on a surface of
genus $5$:  H)  hyperelliptic  component;  O) odd spin structure;
E)~even spin structure.
}
\end{figure}

Consider the following  three diagrams depending on genus $g$ (cf
Figure~\ref{pic:g5} for  the  case  $g=5$). Each diagram contains
$2g-1$  separatrix  loops  ($g\ge  2$)  corresponding  to  $4g-2$
separatrix  rays  $r_1,\dots,r_{4g-2}$.  We  enumerate  the  rays
counterclockwise starting from the one pointing to the ``south''.
Join the following ordered pairs  of  rays by arcs; the order  of
the rays  in each pair  defines the orientation of the separatrix
loop. Join $r_1$  and $r_{2g}$; join $r_{2i+1}$ and $r_{2i}$, for
$i=1,\dots,g-1$;    join    $r_{2i-1}$     and    $r_{2i}$    for
$i=g+1,\dots,2g-1$. The diagrams  differ only by the way we paint
them.

--- Case H. ($g\ge 2$)
Paint with the same  colors  pairs of domains centrally symmetric
to each other.

--- Case O. ($g\ge 2$)
Paint  with  the  same  colors  pairs  of domains symmetric  with
respect to the vertical axis.

--- Case E.  ($g\ge  3$)
Paint with the same colors the sectors corresponding to the loops
$r_{3} r_{2}$ and $r_{4g-5} r_{4g-4}$. Paint with the same colors
the  sectors  corresponding  to  the  loops   $r_{5}  r_{4}$  and
$r_{4g-3} r_{4g-2}$. Paint with the same colors the rest pairs of
domains symmetric to each other with respect to the vertical axis
(cf Figure~\ref{pic:g5}).

\begin{Lemma}
\label{lm:realize3}
Every  diagram of  the  type  H,  O, or  E is  realizable  by the
horizontal foliation of an Abelian differential  from the stratum
$\H(2g-2)$.
\end{Lemma}
\begin{proof}
In  our  case the system of  linear  equations on the lengths  of
separatrix  loops (see  Lemma~\ref{lm:realiz:diag})  is  trivial:
simple separatrix loops bounding the  domains  of  the same color
have equal length. Thus it obviously has positive solutions.
\end{proof}

\begin{Lemma}
\label{lm:3types}
Let  an   Abelian   differential  $\omega$  have  the  horizontal
foliation represented by one of the diagrams $H, O, E$.

---    If     the     diagram     is     diagram     $H$,    then
$\omega\in\mathcal{H}^{hyp}(2g-2)$.

--- If the  diagram  is diagram $O$, then  $\omega$  has odd spin
structure.  If  $g=2$, then the diagram $O$  coincides  with  the
diagram $H$,  and $\omega$ is  hyperelliptic; for $g>2$ it is not
hyperelliptic.

---  If  the diagram is diagram  E,  then $\omega$ has even  spin
structure.  If  $g=3$  then  the diagram $E$ coincides  with  the
diagram $H$,  and $\omega$ is  hyperelliptic; for $g>3$ it is not
hyperelliptic.
\end{Lemma}
\begin{proof}
Diagram  $H$  is  always  centrally  symmetric  ; diagram $O$  is
centrally symmetric only for $g=2$ when  it  coincides  with  the
diagram $H$; diagram $E$ is  centrally  symmetric  only for $g=3$
when  it  coincides  with  the  diagram  $H$.  Thus  according to
Lemma~\ref{lm:hypersymmetry} these are the only cases when we get
a hyperelliptic Abelian differential.

The  parity  of  the  spin  structure  determined  by  an Abelian
differential with  a horizontal foliation  of the type $O$ or $E$
is   computed  inductively   using   Lemma~\ref{lm:erase},   with
Lemma~\ref{lm:g1} serving  as  the  base  of  induction. Strictly
speaking, we can not apply Lemma~\ref{lm:erase} to the bubbling a
handle to the torus  ($g=1$  case) because we introduced bubbling
only at zeroes of positive multiplicity. It is easy to check that
our arguments  work as well for the bubbling  at a regular point,
i.e. at  a ``zero of multiplicity  zero''. The diagram  for $g=1$
case consists of one vertex and one loop.
\end{proof}

\begin{Proposition}
\label{pr:1of3}
Any connected  component  of  the stratum $\mathcal{H}(2g-2)$ for
$g\ge  2$  contains  an  Abelian differential with  a  horizontal
foliation represented by one of the diagrams $H, O, E$.
\end{Proposition}

\begin{Corollary}
\label{cr:3c}
For $g=2$ the stratum $\mathcal{H}(2)$ is connected; it coincides
with  the   hyperelliptic   component.   For  $g=3$  the  stratum
$\mathcal{H}(4)$  has  exactly  two  connected  components:   the
hyperelliptic  one  $\mathcal{H}^{hyp}(4)$, and  one  having  odd
parity of the spin structure $\mathcal{H}^{odd}(4)$. For $g\ge 4$
the  stratum  $\mathcal{H}(2g-2)$  has  exactly  three  different
connected         components:          $\mathcal{H}^{hyp}(2g-2)$,
$\mathcal{H}^{odd}(2g-2)$, and $\mathcal{H}^{even}(2g-2)$.
\end{Corollary}

\begin{proof}
The   Corollary   follows   immediately   from   combination   of
Proposition~\ref{pr:1of3}    with    Lemma~\ref{lm:3types}    and
Lemma~\ref{lm:coinc:diag}.
\end{proof}

The  rest  part  of  this  section is  devoted  to  the  proof of
Proposition~\ref{pr:1of3}.

\begin{proof}[Proof of Proposition~\ref{pr:1of3}]
The  diagrams  which  can  be  obtained  one from  the  other  by
reversing  the  arrows are equivalent in our considerations,  see
Remark~\ref{rm:arrows}. Throughout  this  proof  we mostly do not
distinguish such diagrams.

First     note     that     for    the    connected     component
$\mathcal{H}^{hyp}(2g-2)$  the  statement  of the proposition  is
obvious:  by  Lemmas~\ref{lm:realize3}  and  \ref{lm:3types}  the
diagram  H  is  realizable  by  the  horizontal  foliation  of  a
hyperelliptic Abelian differential.

Every  Riemann  surface  of  genus $g=2$ is  hyperelliptic  which
implies  that   any   Abelian   differential   in   the   stratum
$\mathcal{H}(2)$                is                 hyperelliptic,
$\mathcal{H}(2)=\mathcal{H}^{hyp}(2)$ and hence  $\mathcal{H}(2)$
is connected. Thus for $g=2$ the Proposition is proved.

Assume that  Proposition~\ref{pr:1of3}  is  proved for all genera
smaller than or equal to $g$, where $g\ge 2$. Let us prove it for
genus $g+1$. To make a step  of induction we have to decrease the
genus  of  the  surface  by one.  In  order  to  do this  we  use
Lemma~\ref{lm:accessible:by:bubbling}   saying    that   in   any
connected component of $\mathcal{H}(2g)$, $g\ge 2$,  one can find
an Abelian  differential  $\hat\omega$ obtained from some Abelian
differential $\omega\in\cH(2g-2)$ by ``bubbling a handle''.

``Forget''       the        corresponding       handle       (see
Remark~\ref{rm:forget:recall}).  By  the induction  assumption we
can deform  continuously  the  corresponding Abelian differential
$\omega$ on a surface of genus  $g$ to fit one of the diagrams H,
O, or E.  Now  we can ``bubble the  forgotten  handle'' along the
path     in     the      stratum     $\mathcal{H}(2g-2)$,     see
Lemma~\ref{lm:lift:path}.  Proposition~\ref{pr:1of3}  now follows
from the Lemma below.
\end{proof}

\begin{Lemma}
\label{lm:HOE:mod}
Consider an Abelian  differential  $\hat\omega\in\mathcal{H}(2g)$
obtained  by ``bubbling  a  handle'' at the  zero  of an  Abelian
differential $\omega\in\mathcal{H}(2g-2)$ having  the  horizontal
foliation of one of the types $H, O, E$ in genus $g$. There exist
a continuous  path in $\mathcal{H}(2g)$ joining $\hat\omega$ with
an Abelian differential having the horizontal foliation of one of
the types $H, O, E$ in genus $g+1$.
\end{Lemma}
\begin{proof}
Note that by Lemma~\ref{lm:rotate:handle:diagram}  we  may always
assume that the  pair of simple separatrix loops representing the
handle just ``bubbled'' is symmetric with respect to the vertical
axis.

If the diagram obtained after ``bubbling a handle'' is already of
one of  the types  $H, O, E$, the statement  of the Lemma becomes
trivial.

The pair of  simple separatrix loops representing the handle just
``bubbled''  is  colored in black on all  the  figures.  Speaking
about a pair of simple separatrix  loops we always mean a pair of
simple separatrix loops  representing  a handle, and thus colored
by the same color at the diagram.

The  idea  of the  proof  is  the  following.  Our  diagrams have
numerous  pairs  of  simple  separatrix loops of the  same  color
representing handles. We choose an  appropriate  pair  of  simple
separatrix loops of the same color and temporarily ``forget'' it,
see             Remark~\ref{rm:forget:recall}.              Using
Lemma~\ref{lm:realiz:diag} we  check that the modified diagram is
realizable. Using the induction assumption we deform continuously
the corresponding Abelian differential in genus $g$ to one having
one  of   the  diagrams  $H,   O,  E$.  Then  we  ``recall''  the
``forgotten'' handle.  By Lemma~\ref{lm:lift:path} the  resulting
diagram in genus  $g+1$  represents an Abelian differential which
can be obtained from the initial one $\hat\omega$ by a continuous
deformation inside $\mathcal{H}(2g)$.

We  start  with the  general  case  assuming  that  $g\ge  4$; we
consider the small genera $g=2,3$ separately.

{\bf Case h)} Let the initial Abelian differential $\hat{\omega}$
in genus $g+1$ correspond  to the diagram of the type $H$  with a
``handle bubbled into it''. If the resulting diagram is centrally
symmetric,       it       satisfies        conditions        from
Lemma~\ref{lm:hypersymmetry} and therefore  Abelian  differential
$\hat\omega$  belongs  to  the hyperelliptic connected  component
$\mathcal{H}^{hyp}(2g)$.  Hence,  we can join it by a  continuous
path  with  an Abelian differential corresponding to the  diagram
$H$ in genus $g+1$.

Suppose now that the diagram obtained after ``bubbling a handle''
in the  diagram $H$ is not  centrally symmetric. For  the initial
genus  $g\ge  4$  we obtain a  diagram  in  genus  $g+1$ having a
centrally  symmetric   pair  of  {\it  simple}  separatrix  loops
different from the pair just  ``bubbled''.  We  can always choose
this new pair of  simple separatrix loops in such a way  that the
diagram obtained after  ``forgetting'' this new pair would not be
centrally symmetric. It is easy to see that the resulting diagram
is realizable. Thus, by induction  assumption  we  can deform the
resulting Abelian differential inside $\mathcal{H}(2g-2)$ to  one
corresponding to  one of the  diagrams $O$ or $E$. ``Recall'' the
``forgotten'' handle. The resulting diagram is  obtained from one
of the diagrams of the type  $E$ or $O$ by ``bubbling a handle''.
We have reduced this case to one of the cases o) or e).

{\bf Case  o)} A diagram obtained  by ``bubbling a  handle'' into
the diagram of type $O$ in  genus $g\ge 4$ is either again of the
type   $O$,   or   it   is   of    the    type    presented    at
Figure~\ref{pic:case:o}.

\begin{figure}[htb]
%
\includegraphics{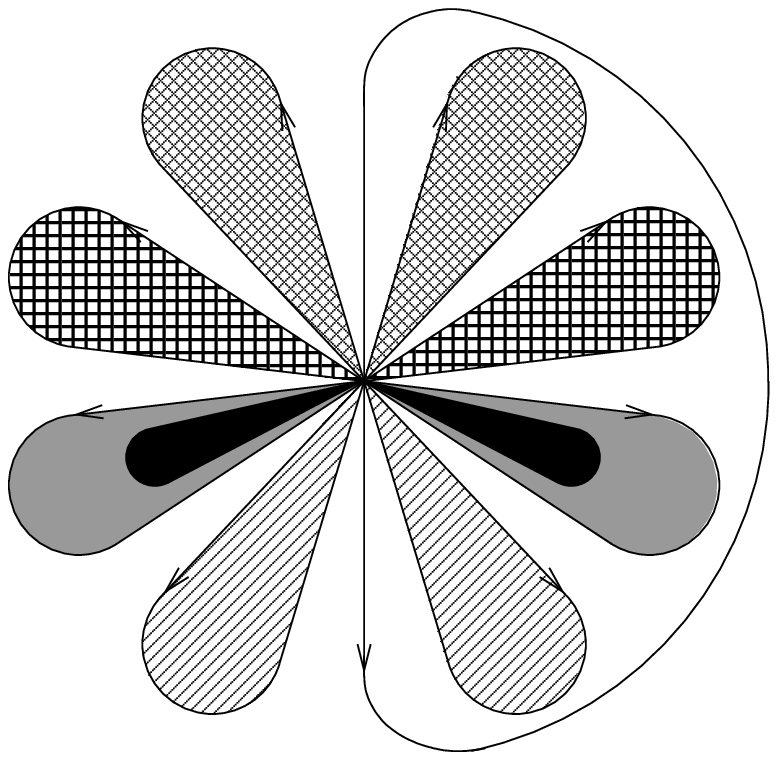}
\vspace{100bp}
\caption{
\label{pic:case:o}
Case o).
}
\end{figure}

In the latter case there is a symmetric pair of simple separatrix
loops next  to the vertical axis. By reversing  the arrows on the
diagram, if necessary, we may assume that this new pair of simple
separatrix loops is  next to the  top vertical ray  (one  without
arrow). Let us mark the top vertical ray, and ``forget'' this new
pair  of  simple  separatrix  loops.  The  resulting  diagram  is
obviously  realizable.  By   Lemma~\ref{lm:change:of:parity}   it
represents an  Abelian  differential $\omega'$ having even parity
of the spin structure. For initial genus $g\ge 4$ it would not be
centrally-symmetric. Thus,  by  the  induction  assumption we can
deform   $\omega'$   by   a    continuous    deformation   inside
$\mathcal{H}(2g-2)$  to  an   Abelian  differential  representing
diagram $E$. ``Recall'' the ``forgotten'' handle  near the marked
ray. It is represented by a pair of simple separatrix  loops next
to the top vertical separatrix ray symmetric with  respect to the
vertical diameter. The diagram thus  obtained  is  diagram $E$ in
genus $g+1$.

{\bf  Case  e)} Consider a diagram obtained  after  ``bubbling  a
handle'' into diagram $E$ in genus  $g\ge 4$. If it is already of
the type $E$, we have nothing to modify.

If  the  new  handle  was  ``bubbled''  inside the  pair  of  top
symmetric sectors  (see,  for  example, Figure~\ref{pic:e:g5}) we
may  turn  the pair  of  black petals  (keeping  fixed the  angle
between them), say, placing them inside  the  bottom  two  petals
(see  Lemma~\ref{lm:rotate:handle:diagram}).  Thus we  may assume
that the top symmetric petals of the diagram  stay unchanged upon
``bubbling a handle'', see, e.g. Figure~\ref{pic:e}.

\begin{figure}[ht!]
%
\includegraphics{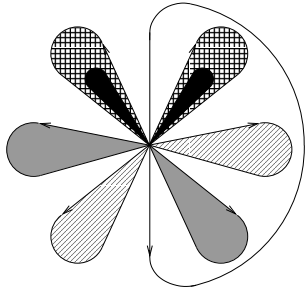}
\vspace{100bp}
\caption{
\label{pic:e:g5}
A handle ``bubbled'' into diagram $E$ in genus $g=4$.
}
\end{figure}

``Forgetting''  the  top  pair  of symmetric petals we  obtain  a
diagram which is obviously realizable by  an Abelian differential
$\omega'\in\mathcal{H}(2g-2)$,     and     which      is      not
centrally-symmetric,  except  only one case when $g=4$ and  black
petals are inserted between the bottom pairs of  petals. Thus, if
we are not in this exceptional case, by  the induction assumption
we  can  deform $\omega'$  by  a  continuous  deformation  inside
$\mathcal{H}(2g-2)$ to an Abelian  differential  representing one
of  the  diagrams  $O$  or $E$. ``Recalling''  the  ``forgotten''
handle we obtain  an Abelian differential representing one of the
diagrams $O$ or $E$ in genus $g+1$.

In the exceptional case when  the  initial genus $g$ is equal  to
$4$ and black  petals are inserted  between the bottom  pairs  of
petals  of  diagram E ``recall'' the ``forgotten''  top  pair  of
symmetric petals at  the initial place.  Turn the pair  of  black
petals (keeping fixed the  angle  between them) by seven sectors.
The black petals are again symmetric with respect to the vertical
axis and we obtain a diagram of the type $E$ in genus $g=5$.

To complete the proof of the  Lemma we need to consider the small
genera $g=2,3$.

\begin{figure}[htb!]
%
\includegraphics{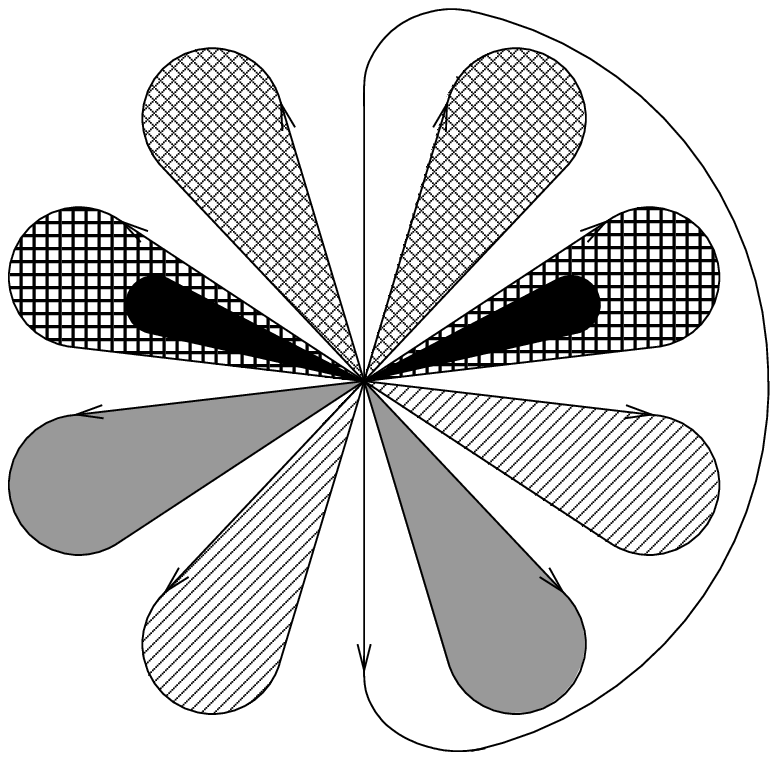}
\vspace{100bp}
\caption{
\label{pic:e}
Case e).
}
\end{figure}

{\bf Small genera)} For initial genus $g=2$ the  diagram $E$ does
not exist,  and the diagrams $H$ and $O$  coincide. Thus, the new
diagram is obtained by ``bubbling a handle'' into  the diagram of
the type $H$ in genus  $g=2$.  The diagram obtained is either  of
the type $O$ in genus $g=3$, or it is centrally-symmetric. In the
latter  case by  Lemma~\ref{lm:hypersymmetry}  the  corresponding
Abelian  differential  belongs to  the  hyperelliptic  component.
Hence,  it  can  be joined by  a  continuous  path  to an Abelian
differential corresponding to the  diagram  of the type $H$. This
completes consideration of genus $g=2$.

Consider now  diagrams in genus $g=3$. The case  when a handle is
``bubbled'' in  a diagram  of the type $O$ in  genus $g=3$ can be
treated the  same way  as in the general case  o), except that we
can  now  obtain  a  diagram which is centrally  symmetric,  thus
reducing it to the case of `bubbling a handle'' into  the diagram
of the type $H$.

\begin{figure}[htb]
%
%
\includegraphics{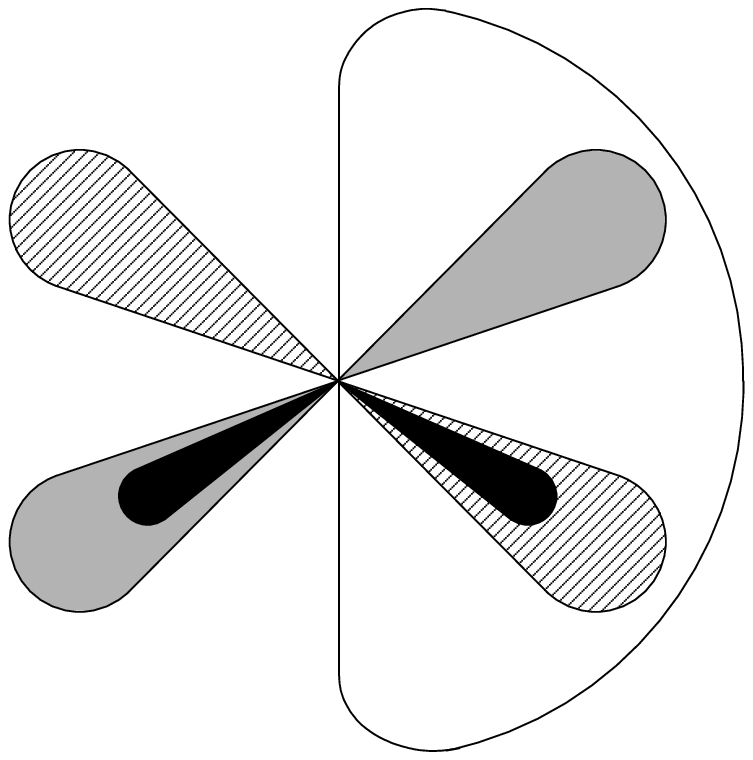}
\includegraphics{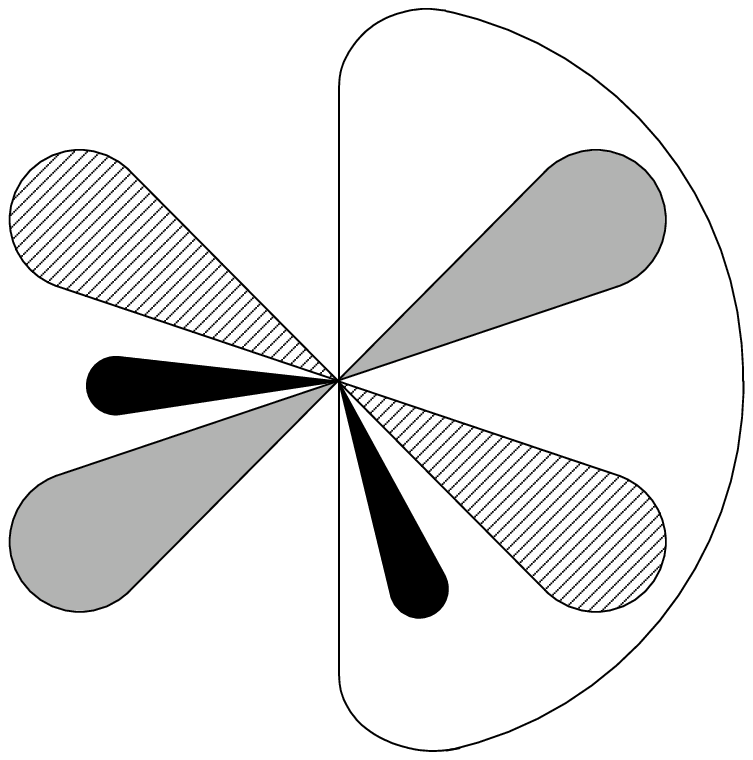}
\includegraphics{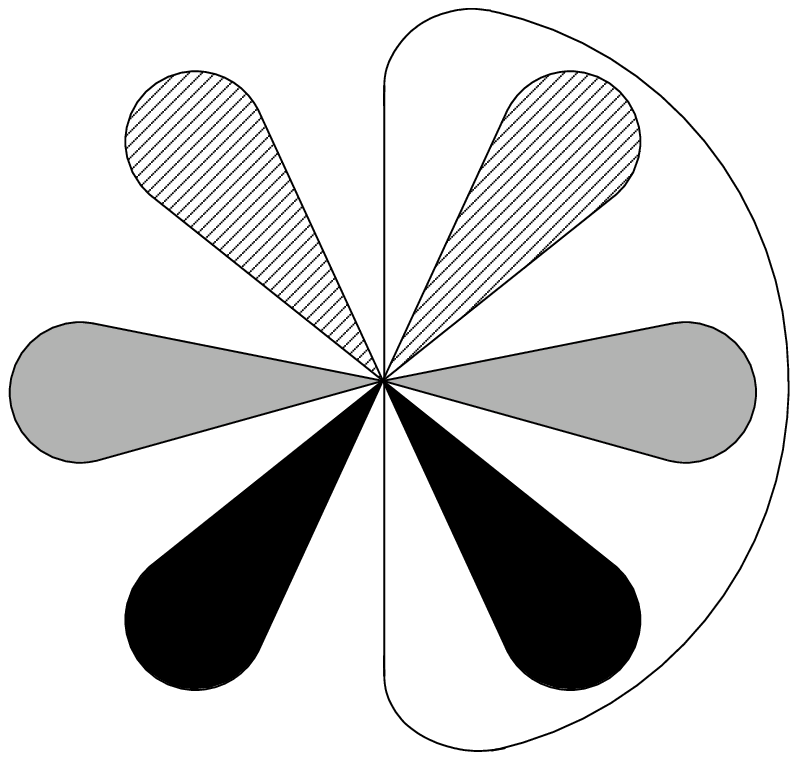}
\vspace{110bp} 
\caption{
\label{pic:hg4}
Modifying diagrams by rotating a pair of petals of the same
color.
}
\end{figure}

The diagrams  $H$ and $E$ coincide  in genus $g=3$.  ``Bubbling a
handle'' in  the diagram $H$ in genus $g=3$  we obtain either one
of  the  diagrams  $H, E$ in genus $g=4$ (up to the change of the
orientation of  the foliation, see Remark~\ref{rm:arrows}) or the
diagram presented at Figure~\ref{pic:hg4} on the left.

In the  latter case we can rotate  the pair  of black petals  one
sector   clockwise   to   obtain   the   diagram   presented   at
Figure~\ref{pic:hg4}        in        the       middle;        by
Lemma~\ref{lm:rotate:handle:diagram}   this   realizable  diagram
represents the same connected component  of  the  same stratum as
the initial  one.  Using Lemma~\ref{lm:realiz:diag} one can check
that erasing {\it any} pair of  petals of the same color from the
resulting diagram  (the  middle  one on Figure~\ref{pic:hg4}) one
gets a realizable  diagram. In particular,  we may think  of  the
pair of petals colored in light grey (the petals corresponding to
North-East and South-West directions) as  of  the  pair of petals
``just       bubbled''.       Hence,      we       may      apply
Lemma~\ref{lm:rotate:handle:diagram}  to  this  pair  of  petals.
Rotating two  sectors clockwise this pair  of petals of  the same
color  (the  ones  corresponding  to  North-East  and  South-West
directions) we  modify the middle diagram on Figure~\ref{pic:hg4}
to the diagram of the type $O$ in genus $g=4$  (the  right one on
Figure~\ref{pic:hg4}).

Lemma~\ref{lm:HOE:mod} is proved.
\end{proof}

\subsection{Stratification of $\cH_g$ near a given stratum}
\label{ss:local:stratification}

Let  $(C,\omega)$  be   a  complex  curve  $C$  with  an  Abelian
differential       $\omega$       representing      a       point
$x=[(C,\omega)]\in\cH_g$ of the moduli space $\cH_g$ of {\it all}
Abelian  differentials,  $g\ge  2$.  A germ $U$ of  the  orbifold
$\cH_g$ at the point $x$ is the quotient  $\tilde U/\Gamma$ where
$\tilde U$ is a germ  of  a complex manifold of dimension  $4g-3$
and  $\Gamma=Aut(C,\omega)$  is a finite group acting on  $\tilde
U$,  preserving   the   base   point   $\tilde   x\in  \tilde  U$
corresponding     to     $x$.    (For     a     generic     point
$x=[(C,\omega)]\in\cH_g$ the group of  automorphisms  $\Gamma$ is
trivial.)  Our  goal here  is  to describe  the  germ $\tilde  U$
together with the stratification of $\tilde  U$ by multiplicities
of zeroes induced from $\cH_g$.

By definition $\tilde U$  is  a universal analytic deformation of
the pair $(C,\omega)$,  i.e. it is  an analytic family  of  pairs
$(C_y,\omega_y)_{y\in \tilde U}$ together with an  identification
$i:(C_{\tilde  x},\omega_{\tilde  x})\simeq  (C,\omega)$  of  the
distinguished element $(C_{\tilde x},\omega_{\tilde  x})$  of the
family with $(C,\omega)$, such  that  any germ of deformations of
$(C,\omega)$ is induced canonically from $\tilde U$.

Let $P_1,\dots,P_n$ be all zeroes of $\omega$ (enumerated in some
order), and $k_1,\dots,k_n$ be their multiplicities.  Let us also
choose local coordinates $z_1,\dots,z_n$ near the  zeroes in such
way that $z_i(P_i)=0$ and $\omega=z_i^{k_i} d z_i$. The choice of
local  coordinates  $z_i$  is  canonical up to  a  transformation
$z_i\mapsto  \xi_i  z_i$  where  $\xi_i$  is  a  root  of  unity,
$\xi_i^{k_i+1}=1$. A deformation $(C_y,\omega_y)_{y\in \tilde U}$
defines for each $i,\,\,\,1\le i\le n$, a deformation of the germ
$z_i^{k_i} d z_i$ of a holomorphic 1-form on a complex curve.

It  is  an  easy   and   well-known  corollary  of  the  standard
deformation theory of singularities of  functions  that  for  the
germs $z^k  d z$ of  1-forms there exists a universal deformation
over  a   germ   of  $(k-1)$-dimensional  manifold.  Namely,  the
following Proposition holds.

\begin{Proposition}
\label{pr:deformation}
Let  $\pi:(\cC,c)\rightarrow(  \cB,b)$ be a map between germs  of
analytic spaces with based points with non-singular fibers of dimension
$1$, and let $\omega\in\Gamma(\cC,T^*_{\cC/\cB})$  be  a $1$-form
along the fibers  of  $\pi$ not  equal  identically to zero.  Let
$z_b$   be   a  local  coordinate  on  $\pi^{-1}(b)$  such   that
$\omega|_{\pi^{-1}(b)}$ is equal  to  $z_b^k dz_b$ for some $k\ge
0$. Assume that $k\ge 1$. Then there exist a unique collection of
$k-1$ functions $a_2,\dots a_k$ on $\cB$ vanishing at  $b$, and a
holomorphic function $z$ on $\mathcal  C$  extending  $z_b$  such
that $$\omega=(z^k+a_2 z^{k-2}+\dots+a_k) dz$$
\end{Proposition}
\begin{proof}
Germs of 1-forms in one variable can be identified by integration
with germs of  functions modulo constants. Now apply the standard
fact: the universal deformation of the germ $z^{k+1}$ is given by
the   formula   $z^{k+1}+a_2'    z^{k-1}+\dots+a'_{k+1}$    where
$(a_i')_{i=2,\dots,k+1}$  are  parameters.  It  remains  to   let
$a_i:=a_i'\cdot(k+1-i)/(k+1)$ for $i=2,\dots,k$.
\end{proof}

We denote by $\cP_k$ the germ in the space $\C{k-1}$ endowed with
coordinates  $a_2,\dots,a_k$  parameterizing $1$-forms  $(z^k+a_2
z^{k-2}+\dots+a_k) dz$  near $z=0\in \C{}$. The above Proposition
says that  any  deformation of a germ of a zero of order $k$ of a
$1$-form  is  induced  canonically  from $\cP_k$. We also  get  a
canonical local coordinate $z$ on  deformed  germs  of curves. In
notations  of  the   Proposition   we  call  a  point  $z_{b'}\in
\pi^{-1}(b')$  given  by  the  equation  $z(z_{b'})=0$  the  {\it
holomorphic  center  of masses} of zeroes (near  the  base  point
$c$). The reason is  that  for a polynomial form $\omega=(z^k+a_2
z^{k-2}+\dots+a_k) dz$ the arithmetic mean of  zeroes of $\omega$
coincides  with  $0\in  \C{}$.  (Note  that  the  notion  of  the
holomorphic center of masses is {\it  not}  invariant  under  the
$GL(2,\R{})_+$-action for the case $k\ge 3$.)

Now  we are  ready  to construct local  coordinates  on the  germ
$\tilde U$ associated  with a global  curve $C$ endowed  with  an
Abelian  differential  $\omega$.  Using the notations  introduced
above, we define a map
\begin{equation}
\label{eq:germ} \Phi:\tilde{U}\rightarrow \prod_{i=1}^n
\cP_{k_i}\times H^1(C,\{P_1,\dots,P_n\};\C{}).
\end{equation}
The components of  this map have the following description. First
of all, for each $i,\,\,\,1\le i\le n$, we  construct a canonical
map      $\tilde      U\rightarrow      \cP_{k_i}$,      applying
Proposition~\ref{pr:deformation} to a neighborhood  of  the point
$P_i$. Secondly,  for deformed curves $(C',\omega')$ with Abelian
differentials we  have  {\it canonical} local holomorphic centers
of  masses   $P_1',\dots,   P_n'\in   C'$.   We   associate  with
$(C',\omega')$  an   element  in  $H^1(C,\{P_1,\dots,P_n\};\C{})$
(close       to       $[\omega]$)       using       $[\omega']\in
H^1(C',\{P'_1,\dots,P'_n\};\C{})$  and  an identification  of the
cohomology      spaces      $$H^1(C,\{P_1,\dots,P_n\};\C{})\simeq
H^1(C',\{P'_1,\dots,P'_n\};\C{})$$  given  by any  continuous map
$$(C,\{P_1,\dots,P_n\})\rightarrow(C',\{P'_1,\dots,P'_n\})$$
close to the identity map (in other words, using the  holonomy of
the Gauss---Manin connection).

An easy calculation with the tangent spaces shows  that $\Phi$ is
a  local  isomorphism. Thus, we constructed, in  a  sense,  local
coordinates  in  a  neighborhood  of  any  point of $\cH_g$.  The
stratification  of  $\tilde  U$  given by the  multiplicities  of
zeroes is obvious. Namely, we should count the  numbers of zeroes
of   given   multiplicities  in   deformed   polynomial   Abelian
differentials. Also, the transversal  slice  in $\tilde U$ to the
stratum containing the base  point  $\tilde x$ is identified with
the product of germs $\cP_{k_i}$.

Using this description of the local structure of  $\cH_g$ we draw
the main conclusion for our classification program:

\begin{Corollary}
\label{cor:local:conn}
For  any  stratum  $\cH(k_1,\dots,k_n)$  of $\cH_g$ and  for  any
connected  component  $S$ of $\H(2g-2)$ there exists exactly  one
connected component $S'$ of $\cH(k_1,\dots,k_n)$ adjacent to $S$,
i.e. such that $S$ is contained in the closure $\overline{S'}$ of
$S'$ in $\cH_g$.
\end{Corollary}
\begin{proof}
Let us prove  that  for any point $x=[(C,\omega)]\in\cH(2g-2)$ of
the  minimal  stratum  $\cH(2g-2)\subset\cH_g$  one  can  find  a
sufficiently  small  neighborhood $U(x)\subset\cH_g$  of  $x$  in
$\cH_g$   such    that    the   intersection   of   $U(x)$   with
$\cH(k_1,\dots,k_n)$ is nonempty and  connected.  Obviously Lemma
follows from this statement.

Applying       formula~\ref{eq:germ}       to       a       point
$x=[(C,\omega)]\in\cH(2g-2)$ of  the minimal stratum  $\cH(2g-2)$
we  establish   a   local   diffeomorphism   between   the   germ
$\tilde{U}(x)$ and $\cP_{2g-2}\times  H^1(C,\{P_1\};\C{})$.  Here
$P_1\in C$ is the single zero of order $2g-2$. Our  statement now
follows  from   the  fact  that   the  germ  of  any  stratum  in
$\cP_{2g-2}$ is nonempty and connected.
\end{proof}

\begin{Corollary}
\label{cor:2comp}
For any parameters $(k_1,\dots, k_n)$, where  $\sum k_i=2g-2$ and
$n\ge    2$,     any     component    $S$    of    the    stratum
$\cH(k_1+k_2,k_3,\dots,k_n)$ of $\cH_g$ and for any two connected
components $S_1, S_2$ of $\H(2g-2)$ to  which  $S$  is  adjacent,
there     exists     a    connected     component     $S'$     of
$\cH(k_1,k_2,\dots,k_n)$ which  is  also  adjacent  to  $S_1$ and
$S_2$.
\end{Corollary}
\begin{proof}
By assumption  we  have  $\overline{S}\supset  S_1\cup  S_2$.  It
follows  from  our  picture  of   the   local   behavior  of  the
stratification  that  there exists a connected component $S'$  of
$\cH(k_1,k_2,\dots,k_n)$      adjacent      to     $S$,      i.e.
$\overline{S'}\supset S$. Therefore we have $\overline{S'}\supset
\overline{S}\supset S_1\cup S_2$.
\end{proof}

\subsection{Merging zeroes and adjacency to the minimal stratum}
\label{ss:merge}

In this section we prove

\begin{Proposition}[Merging zeroes]
\label{pr:col:zer:2}
For any  given  parameters  $k_1,\dots,k_n$,  where  $n\ge 2$ and
$\sum_i k_i=  2g-2$, the closure  of any connected component of a
stratum   $\mathcal{H}(k_1,k_2,\dots,k_n)\subset\cH_g$   contains
some      connected      component       of      the      stratum
$\mathcal{H}(k_1+k_2,k_3,\dots,k_n)$.
\end{Proposition}

\begin{Remark}
Note that in general the similar statement is no longer  true for
the  strata  of  quadratic  differentials.  Say,  one of the  two
connected components  of  the  stratum  of  meromorphic quadratic
differentials having  a single simple  pole, and a single zero of
degree $9$ is  adjacent to the stratum of quadratic differentials
with a single  zero  of degree $8$, and  the  other component ---
not.
\end{Remark}

We start  the proof with  the following technical result (see the
similar  result   in~\cite{M79}  for  the  principal  stratum  of
quadratic differentials).

\begin{Lemma}
\label{lm:1cyl}
Every  connected  component of every stratum contains an  Abelian
differential whose  horizontal  foliation  has only closed leaves
and the corresponding diagram has only one cylinder.
\end{Lemma}
\begin{proof}
First of all, let us  choose  an Abelian differential in a  given
component of a stratum whose  {\it  vertical}  foliation has only
closed   leaves,   see   Lemma~\ref{lm:closedleaves}.   Deforming
slightly  this  differential  preserving  the  structure  of  the
vertical measured foliation we can  assume  that  the  horizontal
foliation is  uniquely  ergodic,  in  particular,  minimal.  This
follows immediately  from results of H.~Masur, see\cite{M82}, and
W.~A.~Veech, see~\cite{V82}.

Let us pick a point $P$  on our surface which is not connected by
a leaf of the horizontal foliation to a critical point.  The leaf
of the horizontal  foliation starting at $P$ is everywhere dense.
For any $\epsilon>0$ we can now follow the leaf until  it returns
for the first time to the distance $\epsilon$ from the point $P$.
For sufficiently small  $\epsilon$ if we connect the endpoints of
our piece  of horizontal leaf  by the geodesic interval of length
$\epsilon$, and then perturb slightly the  obtained closed curve,
we  obtain  a  smooth  closed curve $\gamma$ transversal  to  the
vertical foliation (see~\cite{Z} for details).

If we  choose curve $\gamma$ long  enough, it will  intersect all
cylinders of the  vertical  foliation and, moreover, all vertical
saddle connections. Our goal now is to modify  the flat structure
on our surface preserving the  vertical  measured  foliation  and
making the horizontal foliation satisfy  the  properties  in  the
Lemma.

Let us remove  all vertical saddle connections and curve $\gamma$
from  our  surface; let us forget the  horizontal  foliation.  We
obtain a finite  collection  of curvilinear $4$-gons endowed with
vertical  measured  foliations. Now let us construct on  $4$-gons
new  horizontal  measured  foliations,  in  such   way  that  the
horizontal parts of the boundary of every $4$-gon would be leaves
of the new foliation. Namely, we say that the vertical  length of
each $4$-gon is equal to one, and if a critical point is situated
on the vertical part  of the boundary of a $4$-gon, it  should be
exactly in the middle. The last condition can be always fulfilled
because by  our assumption we  cannot have more than one critical
point on any vertical part of the boundary of any $4$-gon.

We endow  each $4$-gon with the  canonical flat metric  (and with
two  measured  foliations) and  then  glue  the  rectangles  thus
obtained together.  On the new  surface both the vertical and the
horizontal  foliations  have  only  closed  leaves.  The  diagram
corresponding  to  the vertical  foliation  of  the  new  surface
coincides with the diagram for the  initial Abelian differential,
thus  we  land  to  the  same  connected component. Finally,  the
horizontal foliation of the new  surface  has  only one cylinder,
the union of two strips  of  width  $1/2$  on the left and on the
right of the closed line $\gamma$.
\end{proof}

\begin{Remark}
Actually,  we  have  proved  a stronger statement:  {\it  Abelian
differentials  satisfying  conditions of  Lemma~\ref{lm:1cyl} are
dense in every connected component of any stratum.}
\end{Remark}
   %
Here is the proof of this statement.
We did  not change the vertical foliation  at all,  so we do  not
deform the closed 1-form $\omega_v$.  We  modify  the  horizontal
foliation in  two steps: at the first step  we perturb the closed
1-form $\omega_h$  to  some  $\omega'_h$  to  make the horizontal
foliation  uniquely  ergodic.  At  the  second  step  we  replace
$\omega'_h$  by   a   closed   1-form   $\omega''_h$,
whose cohomology class,
by
construction, is  Poincar\'e-dual  to  the  cycle $[\gamma]$; the
final  horizontal   foliation   is   the   kernel   foliation  of
$\omega''_h$.  Note  that  it  follows  from  ergodicity  of  the
intermediate  foliation  corresponding to  $\omega'_h$  and  from
definition of  $\omega''_h$  that  choosing  the  curve  $\gamma$
sufficiently long we get
$$
\int_\rho \omega'_h\approx \frac{1}{|\gamma|}\cdot
(\text{number of intersections of $\gamma$ with $\rho$})
\approx \frac{1}{|\gamma|}\int_\rho \omega''_h
$$
for any path $\rho$ transversal  to  $\omega'_h$  (here we assume
that the total area of the surface measured in the flat metric is
normalized   to   one).   By   construction   the   integral   of
$\frac{1}{|\gamma|}\cdot\omega''_h$  along  a piece  of  leaf  of
$\omega'$ is  close to zero (provided this piece  of leaf is much
shorter than $\gamma$). Thus, choosing $\gamma$ sufficiently long
we can  make $\omega'_h$ and  $\frac{1}{|\gamma|}\cdot\omega''_h$
arbitrarily  close.  Hence, the  resulting  Abelian  differential
determined     by      a      pair      of     closed     1-forms
$(\omega_v,\frac{1}{|\gamma|}\cdot\omega''_h)$  is  close to  the
initial      Abelian      differential      corresponding      to
$(\omega_v,\omega_h)$.
   %

Now we are ready to prove Proposition~\ref{pr:col:zer:2}.

\begin{proof}[Proof of Proposition~\ref{pr:col:zer:2}]
Given   a    connected    component    $S'$    of    a    stratum
$\mathcal{H}(k_1,...,k_n)$,   choose   an  Abelian   differential
$\omega$ in  this  component  as in Lemma~\ref{lm:1cyl}. Consider
the  diagram   of   its  horizontal  foliation.  Consider  saddle
connections of this diagram. It is easy to see that {\it any}
choice of strictly  positive  lengths of these saddle connections
gives  a  solution of  the  system  of  linear  equations  (as in
Lemma~\ref{lm:realiz:diag}).

Since  the  union  of  all nonsingular leaves of  the  horizontal
foliation forms a single cylinder,  the  underlying  graph of the
diagram is connected. In particular, every saddle is connected to
at least one {\it other} saddle by a separatrix.

Consider a diagram  obtained  by shrinking this saddle connection
to a point. The diagram is obviously realizable; it represents an
Abelian        differential        from        the        stratum
$\mathcal{H}(k_1+k_{j_1},k_2,k_3,\dots,\widehat{k_{j_1}},\dots,k_n)$.
(For the  moment we cannot control the index  $j_1$; we just know
that   $j_1\neq   1$.)  We  may  shrink  the  saddle   connection
continuously without  changing  other  parameters of the diagram.
Thus  we  get  a  continuous
path with interior in
the  chosen  component  of
$\mathcal{H}(k_1,\dots,k_n)$
and with one of the endpoints belonging to
$\mathcal{H}(k_1+k_{j_1},k_2,k_3,\dots,\widehat{k_{j_1}},\dots,k_n)$.
We have  proved that the  connected component $S'$ is adjacent to
some connected  component of a stratum  with a smaller  number of
zeroes. Repeating  inductively  this  procedure  we conclude that
$S'$ is adjacent to some connected component $S_1\in\cH(2g-2)$ of
the minimal stratum.

By   Corollary~\ref{cor:local:conn}   there  exist   a  connected
component  $S  \subset  \cH(k_1+k_2,k_3,\dots,k_n)$  adjacent  to
$S_1$.  By  Corollary~\ref{cor:2comp}  there  exist  a  connected
component $S''\subset\cH(k_1,\dots,k_n)$ adjacent to  $S$  and to
$S_1$.        Since        both       connected        components
$S',S''\in\cH(k_1,\dots,k_n)$ are adjacent  to  $S_1\in\cH(2g-2)$
Corollary~\ref{cor:local:conn}   implies    that   $S'=S''$.   By
construction    $S''$    is    adjacent     to     the    stratum
$\cH(k_1+k_2,k_3,\dots,k_n)$   which   completes  the   proof  of
Proposition~\ref{pr:col:zer:2}.
\end{proof}

It would be convenient to formulate an intermediate result of our
proof as a separate Corollary.

\begin{Corollary}
\label{cr:collaps:zeroes}
The    closure     of    any    component    of    any    stratum
$\mathcal{H}(k_1,\dots,k_n)$  contains  a connected  component of
the stratum $\mathcal{H}(2g-2)$, where $2g-2=k_1+\dots+k_n$.
\end{Corollary}

   %

\subsection{Connected components of general strata}
\label{ss:other:strata}

First  note  that  any  stratum  $\mathcal{H}(k_1,\dots,k_n)$  is
nonempty for any collection of positive integers $k_i$, such that
the sum  of all $k_i$ is even, see~\cite{MS}.  Another way to see
this is to perturb Abelian differentials from $\mathcal{H}(2g-2)$
(see  Corollary~\ref{cor:local:conn})  which we  have constructed
directly, see Lemma~\ref{lm:realize3}.

Moreover,    for    any   positive    integers   $l_1,\dots,l_n$,
$l_1+\dots+l_n=g-1$,   perturbing   an    Abelian    differential
$\omega\in\mathcal{H}(2g-2)$  we  obtain an  Abelian differential
$\omega'\in\mathcal{H}(2l_1,\dots,2l_n)$, having  the same parity
of  the  spin  structure  as  the  initial  Abelian  differential
$\omega\in\mathcal{H}(2g-2)$, as follows from  invariance  of the
parity  of  the  spin-structure  under  continuous  deformations,
see~\cite{At}, \cite{Mum}. Thus, using our direct construction of
hyperelliptic     components     $\mathcal{H}^{hyp}(2g-2)$    and
$\mathcal{H}^{hyp}(g-1,g-1)$,     and     perturbing      Abelian
differentials      from      the       connected       components
$\mathcal{H}^{odd}(2g-2)$, $\mathcal{H}^{even}(2g-2)$  we can get
all    the   components    listed    in    Theorems~\ref{th:ccor}
and~\ref{th:g123}.  To  complete the proofs of these Theorems  we
have to prove that all the components listed in the  Theorems are
connected, and that there are no other components.

\begin{Lemma}
\label{lm:atmost3c}
Any  stratum  $\mathcal{H}(k_1,\dots,k_n)$  has  at  most   three
connected components.
\end{Lemma}
\begin{proof}
The   statement   follows   immediately   from   combination   of
Corollary~\ref{cr:3c},    Corollary~\ref{cor:local:conn},     and
Corollary~\ref{cr:collaps:zeroes}.
\end{proof}

A  component  of  $\mathcal{H}(2g-2)$  uniquely  determines   the
embodying  component  of  $\mathcal{H}(k_1,\dots,k_n)$,  but  the
embodying component may contain in the closure two, or even three
components of $\mathcal{H}(2g-2)$, see Propositions~\ref{pr:k1k2}
and~\ref{pr:even2odd} below.

\begin{Proposition}
\label{pr:k1k2}
For any genus $g\ge 3$  and  any  $k$,  $1\le k <g-1$, there is a
continuous    path    $\gamma:[0;1]    \to   \H_g$   such    that
$\gamma(]0,1[)\subset      \mathcal{H}(k,2g-k-2)$,       endpoint
$\gamma(1)$   belongs   to   the   hyperelliptic   component   of
$\mathcal{H}(2g-2)$, and  endpoint  $\gamma(0)$ belongs to one of
two nonhyperelliptic components of $\mathcal{H}(2g-2)$.
\end{Proposition}
\begin{proof}
The path  is  presented at Figure~\ref{pic:h2nonh}. Every diagram
is easily seen to be  realizable.  Note that we may preserve  the
heights  and the  widths  (measured in our  flat  metric) of  all
cylinders along the path; we  just  change  the identification of
the boundary components. This implies the continuity of the path.

The  bottom  diagram   is   centrally  symmetric  and  obeys  the
conditions of  Lemma~\ref{lm:hypersymmetry}. Thus it  corresponds
to  an  Abelian  differential from $\mathcal{H}^{hyp}(2g-2)$.  By
assumption of the Proposition  $k>0$  and $g-k-1>0$. Thus the top
diagram is not centrally symmetric (see Figure~\ref{pic:h2nonh}).
Hence it corresponds to a nonhyperelliptic component.
\end{proof}

\begin{Proposition}
\label{pr:even2odd}
For any  genus $g\ge 4$ and any $k$, $1\le k \le g/2$, there is a
continuous    path    $\gamma:[0;1]    \to   \H_g$   such    that
$\gamma(]0,1[)\subset\mathcal{H}(2k-1,2(g-k)-1)$,   one   of  the
endpoints lies  in  the  component $\mathcal{H}^{even}(2g-2)$ and
another endpoint lies in the component $\mathcal{H}^{odd}(2g-2)$.
\end{Proposition}
\begin{proof}
The path is presented at Figure~\ref{pic:e2odd}. Again it is easy
to  see  that all the diagrams  are  realizable, and that we  may
preserve the  heights and the widths  of all cylinders  along the
path. It  is easy  to see that neither top  nor bottom diagram is
centrally  symmetric  ,  thus  they  represent   nonhyperelliptic
components of $\mathcal{H}(2g-2)$.

\begin{figure}
   %
\includegraphics{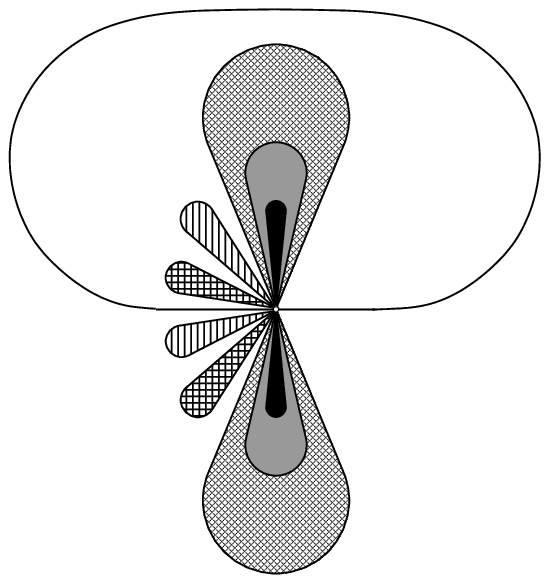}
\includegraphics{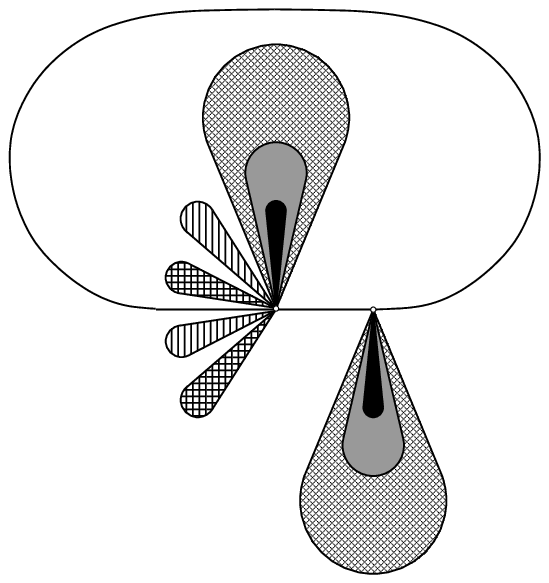}
\includegraphics{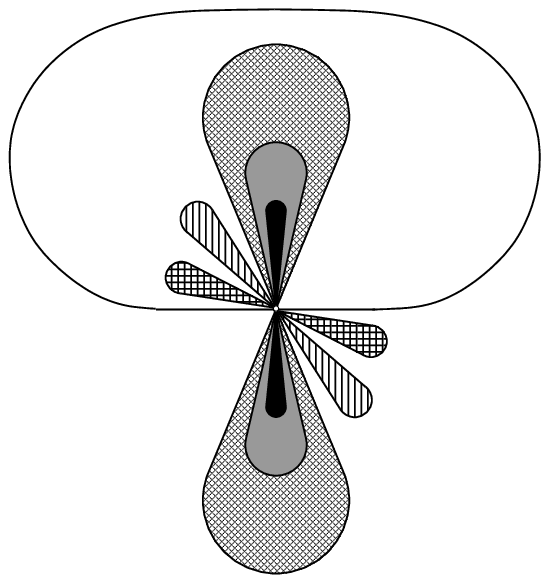}
\vspace{546.93bp} 
\begin{picture}(0,0)(0,0) 
\put(0,0){\begin{picture}(0,0)(0,0)
\put(-105,245){$g-k-1$ loop $\bigg\{ \bigg.$}
\put(60,215){$k$ loops}
\end{picture}}
\end{picture}

\caption{
\label{pic:h2nonh}
A  path  in  $\mathcal{H}(k,2g-2-k)$  joining  hyperelliptic  and
nonhyperelliptic components of $\mathcal{H}(2g-2)$.
}
\end{figure}

\begin{figure}
\includegraphics{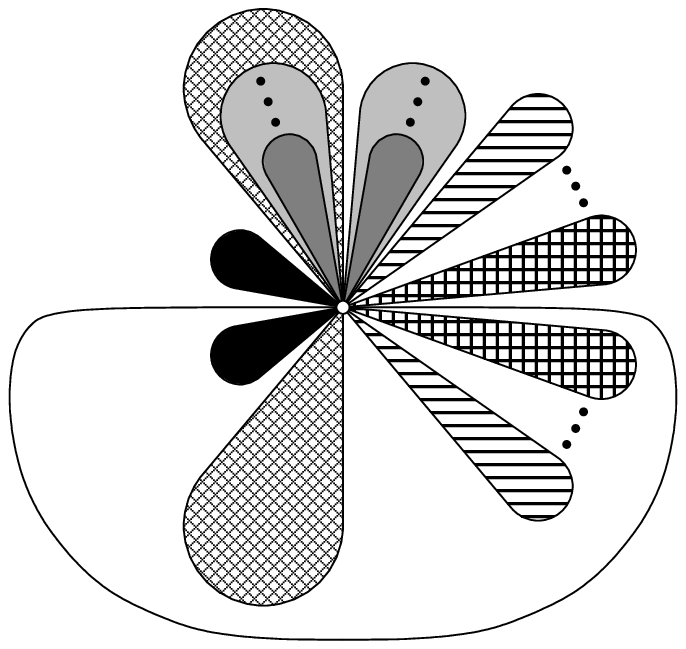}
\includegraphics{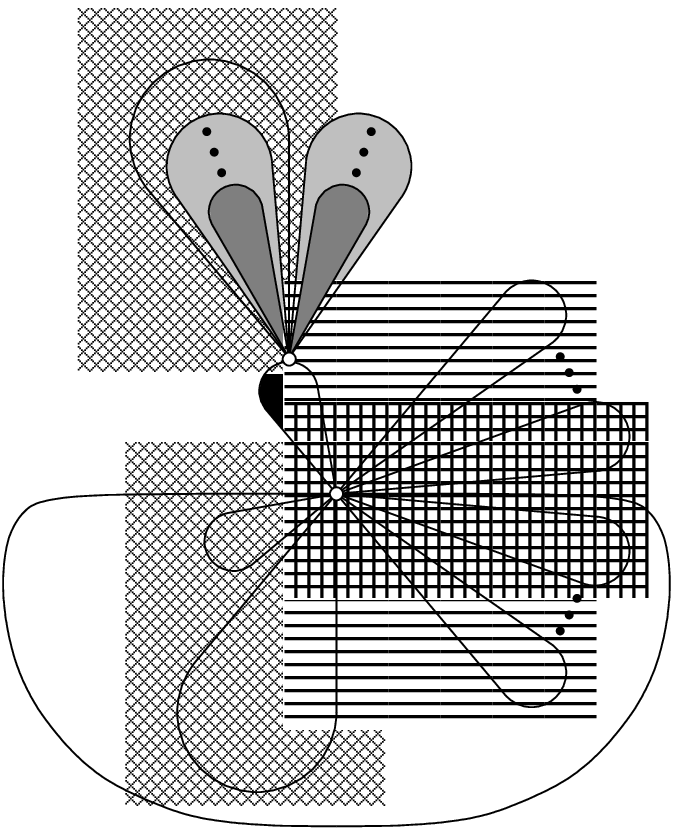}
\includegraphics{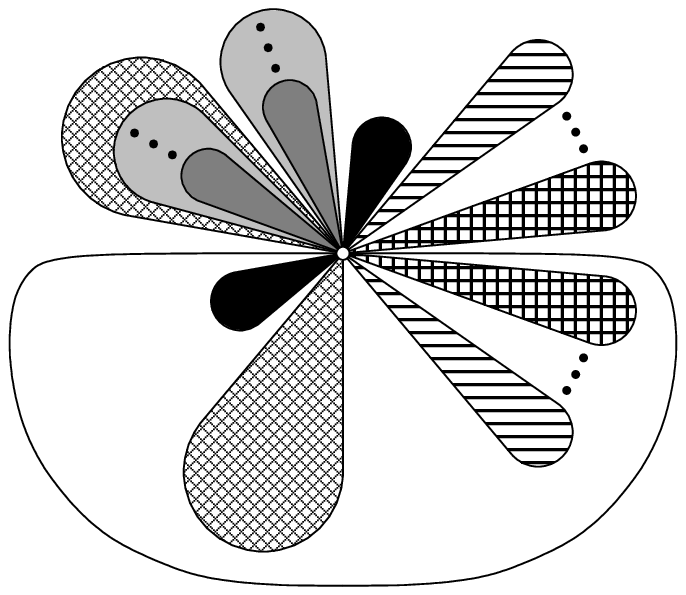}
\vspace{546.93bp} 
\begin{picture}(0,0)(0,0) 
\put(0,0){\begin{picture}(0,0)(0,0) 
\put(30,339){$\bigg. \bigg\} k-1$ loop}
\put(85,290){$\Bigg. \Bigg\} g-k-2$ loops}
\end{picture}}
\end{picture}
\caption{
\label{pic:e2odd}
A      path      in       $\mathcal{H}(2k-1,2(g-k)-1)$ joining
$\mathcal{H}^{even}(2k-2)$ and $\mathcal{H}^{odd}(2k-2)$.
}
\end{figure}


Let  us  prove  that  the  parities   of   the   spin   structure
corresponding to the top  and  the bottom diagrams are different.
Constructing the diagrams we may ``bubble'' the handle painted in
black at the very last step. Thus we  may ``erase'' corresponding
pair of simple loops  both on the top and on the  bottom diagram.
After  erasing  this pair  of  simple loops  we  obtain the  same
diagram on  top and at the bottom. Let  $\varphi_0$ be the parity
of the spin structure corresponding to the diagram thus obtained.
By  Lemma~\ref{lm:erase}   the   parity  of  the  spin  structure
corresponding to  the  initial  top  diagram  equals $\varphi_0$,
while  the  parity  of  the spin structure corresponding  to  the
initial bottom diagram equals $\varphi_0+1$.
\end{proof}

Now  we  are ready  to  finish the  proof  of the  classification
Theorems~\ref{th:ccor} and~\ref{th:g123}. Recall  that we already
have a surjection  from the set  of connected components  of  the
minimal stratum $\H(2g-2)$  to the set of connected components of
any other  stratum, and  also two invariants (to be  or not to be
hyperelliptic, or to  have even or odd spin structure) separating
connecting components.

Let us start with the  case  $g\ge 4$. For stratum $\H(2g-2)$  we
have achieved  already  the classification with three components,
see  Corollary~\ref{cr:3c}.  Similarly, we treat the case of  the
stratum $\H(2l,2l)$,  $l\ge 2$. For all  the other strata  of the
form  $\H(2l_1,\dots,  2l_n)$  we  will have only  two  connected
components distinguished by the parity  of  spin  structure.  The
reason  is  that the  component  adjacent  to  the  hyperelliptic
component  in  $\H(2g-2)$ is also adjacent to a  nonhyperelliptic
component,   as  follows   from   Proposition~\ref{pr:k1k2}   and
Corollary~\ref{cor:2comp}.

For strata  $\H(2k-1,2k-1)$ with $k\ge 2$  we will have  only two
components distinguished  now  by  hyperellipticity,  now  we use
Proposition~\ref{pr:even2odd} and Corollary~\ref{cor:2comp}.

For any other stratum we will have at least one of multiplicities
which   is  odd   and   not  equal  to   $g-1$.   In  this   case
Propositions~\ref{pr:k1k2}  and  \ref{pr:even2odd} together  with
Corollary~\ref{cor:2comp} finish  the  job,  showing that we have
only    one    connected    component.     Thus,     we    proved
Theorem~\ref{th:ccor}.

In  the  case $g=3$ the upper  bound  of the number of  connected
components is $2$. We treat cases $\H(4)$ and  $\H(2,2)$ as above
and get two components. In all  other cases we will have at least
one  of  multiplicities  equal  to  $1$   and   here   we   apply
Proposition~\ref{pr:k1k2}. In the  case  $g=2$ the upper bound is
already  equal  to  $1$  and  we  conclude that  all  strata  are
connected.\qed

\appendix

\section{Rauzy classes and zippered rectangles}
\label{s:Rcl}
\label{a:rauzy:cl}

\subsection{Interval exchange transformations}
\label{ss:iet}

In this  section we recall  the notions of {\it interval exchange
transformation}, of  {\it Rauzy class}, see~\cite{Rauzy}, and the
construction  of  a   complex   curve  endowed  with  an  Abelian
differential   by   means   of  {\it  ``zippered   rectangles''},
see~\cite{V82}.

Consider  an  interval  $I\subset\R{}$,  and  cut   it  into  $m$
subintervals of  lengths  $\lambda_1,\dots,  \lambda_m$. Now glue
the subintervals together in another  order,  according  to  some
permutation   $\pi\in   \mathfrak{S}_m$    and   preserving   the
orientation. We again obtain an interval $I$ of  the same length,
and hence we  get  a  mapping $T:I\to I$, which  is  called  {\it
interval exchange  transformation}.  Our  mapping  is a piecewise
isometry, and it preserves the orientation  and Lebesgue measure.
It is singular  at  the points  of  cuts, unless two  consecutive
intervals separated  by a point  of cut are mapped to consecutive
intervals in the image.

\begin{Remark}
\label{rm:pi:and:pi:inv}
Note, that  actually there are  two ways to glue the subintervals
``according  to  permutation $\pi$''.  We  may  send  the  $k$-th
interval to  the place $\pi(k)$, or we may  have the intervals in
the image appear in  the  order $\pi(1),\dots,\pi(m)$. We use the
first  way;  under  this  choice  the  second way corresponds  to
permutation $\pi^{-1}$.
\end{Remark}

Given an interval exchange transformation $T$  corresponding to a
pair       $(\lambda,\pi)$,       $\lambda\in\R{m}_+$,       $\pi
\in\mathfrak{S}_m$,               set                $\beta_0=0$,
$\beta_i=\sum_{j=1}^i\lambda_j$,        and        $I_i         =
[\beta_{i-1},\beta_i[$.           Define           skew-symmetric
$m\hspace{-3pt}\times\hspace{-3pt}m$-matrix   $\Omega(\pi)$    as
follows:
\begin{equation}
\label{eq:Omega}
  \Omega_{ij}(\pi)=\left\{
\begin{array}{rl}
      1 & \text{if $i<j$ and $\pi(i)>\pi(j)$}\\
     -1 & \text{if $i>j$ and $\pi(i)<\pi(j)$}\\
      0 & \text{otherwise}
\end{array}
\right.
\end{equation}
Consider a translation  vector  $\delta=\Omega(\pi)\cdot\lambda$.
Our interval exchange transformation $T$ is defined as follows:
$$ T(x)=x+\delta_i,\text{\hspace{1cm} for $x\in I_i$, $1\le i\le
m$} $$

\subsection{Extended Rauzy classes}
\label{ss:ext:rauzy:cl}

Consider      an      Abelian       differential       $\omega\in
\mathcal{H}(k_1,...,k_n)$  on   a  surface  of  genus  $g\ge  2$.
Consider   corresponding  vertical   (or   horizontal)   measured
foliation on  the  Riemann  surface.  For  generic $\omega$ every
nonsingular leaf  of the foliation  is dense on the surface. Take
an interval $I$ transversal to  the  foliation.  Our foliation is
oriented, so it defines the Poincar\'e map (the first return map)
$I\to  I$.  It  is  easy  to see that the map $T$ is an  interval
exchange transformation.  The  number of intervals under exchange
is $2g+n-1$, $2g+n$, or $2g+n+1$ depending on the  choice of $I$.
(Morally,  one  has  to  place the endpoints of  the  transversal
interval on  the critical leaves of  the foliation to  obtain the
minimal possible  number  of  subintervals.)  In  particular  the
choice of  transversal  interval  $I$ determines some permutation
$\pi$. Consider  the  set  $\mathfrak{R}_{ex}$  of  all  possible
permutations $\pi\in\mathfrak{S}_{2g+n-1}$ which  can be obtained
by choosing different transversal intervals $I$. It was proved by
W.A.Veech in~\cite{V82} that the set $\mathfrak{R}_{ex}$ does not
depend on the choice of a generic representative  $\omega$ in any
connected  component  of  $\mathcal{H}(k_1,\dots,k_n)$.  The  set
$\mathfrak{R}_{ex}$  is  called  {\it   extended   Rauzy  class},
see~\cite{Rauzy},~\cite{V82}.

Conversely, given  an interval exchange transformation $T:I\to I$
one  can  construct  a   complex   curve  $C_g$  and  an  Abelian
differential $\omega$ on it, such that the Poincar\'e map induced
by the vertical foliation on the appropriate embedded subinterval
would give  the  initial interval exchange transformation. Though
the choice of the pair $(C_g,\omega)$ is not  unique, topology of
$(C_g,\omega)$ (genus, degrees $k_1, \dots,  k_n$  of  zeroes  of
$\omega$,    and    even    the     connected     component    of
$\mathcal{H}(k_1,\dots,k_n)$)  are  uniquely  determined  by  the
permutation $\pi$. We review the construction  of {\it suspension
over     an     interval      exchange     transformation}     in
Appendix~\ref{ss:ziprec}, more details can be found in~\cite{M82}
or in~\cite{V82}.

In the section below we present a direct combinatorial definition
of the extended Rauzy class, see~\cite{Rauzy}, \cite{V82}.

\subsection{Combinatorics of Rauzy classes}
\label{ss:comb:rauzy:cl}

Note,     that      if     for     some     $k<m$     we     have
$\pi\{1,\dots,k\}=\{1,\dots,k\}$, then the corresponding interval
exchange transformation $T$ decomposes into two interval exchange
transformations. We consider only the class $\mathfrak{S}^0_m$ of
{\it irreducible} permutations --- those which  have no invariant
subsets of the form $\{1,\dots,k\}$, where $1\le k < m$.

Permutation $\pi$  is called {\it  degenerate} if it obeys one of
the   following   conditions   (see  3.1--3.3  in~\cite{M82}   or
equivalent conditions 5.1--5.5 in~\cite{V82}):

\noindent for some $1\le j < m$,
\begin{equation*}
\begin{array}{rcl}
\pi(j)   &=& m \\
\pi(j+1) &=& 1 \\
\pi(1)   &=& \pi(m) +1
\end{array}
\end{equation*}
for some $1\le j < m$,
\begin{equation*}
\begin{array}{rcl}
\pi(j+1) &=& 1 \\
\pi(1)   &=& \pi(j) +1
\end{array}
\end{equation*}
for some $1\le j < m$,
\begin{equation*}
\begin{array}{rcl}
\pi(j+1) &=& \pi(m)+1\\
\pi(j)   &=& m
\end{array}
\end{equation*}
Otherwise permutation $\pi$ is called {\it nondegenerate}.

We denote  by $\tau_k\in\mathfrak{S}_m$, $1\le k<m$ the following
permutation:
\begin{equation*}
\begin{array}{lll}
\tau_k &=& (1,2,\dots,k,k+2,\dots,m,k+1)
\;\;\; 1\le k < m-1 \nonumber\\
\tau_{m-1} &=& (1,2,\dots,m)=id \nonumber
\end{array}
\end{equation*}
Permutation  $\tau_k$  cyclically moves one step forward all  the
elements occurring after the element $k$.

Consider two  maps $a,b:\mathfrak{S}_m^0\to \mathfrak{S}_m^0$  on
the set of irreducible permutations (see~\cite{Rauzy}):
\begin{equation*}
\begin{array}{rcl}
a(\pi) &=& \pi\cdot\tau_{\pi^{-1}(m)}^{-1}\\
b(\pi) &=& \tau_{\pi(m)}\cdot\pi
\end{array}
\end{equation*}
where one should consider product of  permutations as composition
of    operators     ---     from     right    to    left.    Say,
$b(2,3,1)=(1,3,2)\cdot(2,3,1)=(3,2,1)$.    We    may     consider
permutation as a pair of orderings of a finite set:  a ``domain''
ordering and an  ``image''  ordering. Operator $b$ corresponds to
the modification of the  image  ordering by cyclically moving one
step forward those letters occurring after the image  of the last
letter in  the domain, i.e.,  after the letter $m$. Operation $a$
corresponds to the modification of the ordering of  the domain by
cyclically moving one step forward those  letters occurring after
one going to the last place, i.e., after $\pi^{-1}(m)$.

Note, that
\begin{equation*}
\left(a(\pi)\right)^{-1}=b(\pi^{-1})
\end{equation*}

In components the maps $a,b$ are represented as follows:
\begin{equation*}
\begin{array}{rcl}
a(\pi)(j) &=&
   \begin{cases}
   \pi(j)   & j\le \pi^{-1}(m)\\
   \pi(m)   & j= \pi^{-1}(m)+1\\
   \pi(j-1) & \text{other $j$}
   \end{cases}
\\ &&\\
b(\pi)(j) &=&
   \begin{cases}
   \pi(j)   & \pi(j)\le \pi(m)\\
   \pi(j)+1 & \pi(m)<\pi(j)<m\\
   \pi(m)+1 & \pi(j)=m
   \end{cases}
\end{array}
\end{equation*}

\begin{Definition}
The {\it  Rauzy  class}  $\mathfrak{R}(\pi)$  of  an  irreducible
permutation  $\pi$  is  the  subset  of   those  permutations  in
$\mathfrak{S}_m^0$  which  can  be  obtained from $\pi$  by  some
composition of maps $a$ and $b$.
\end{Definition}

Consider the permutation $\pi_0=(m,m-1,\dots,2,1)$, and the map
$$
\operatorname{Ad}_{\pi_0}: \pi \mapsto \pi_0^{-1} \pi \pi_0 =
\pi_0 \pi \pi_0
$$
Note that the map $\operatorname{Ad}_{\pi_0}$ maps an irreducible
permutation to an irreducible one.

\begin{Definition}
The {\it  extended  Rauzy  class}  $\mathfrak{R}_{ex}(\pi)$ of an
irreducible permutation $\pi$ is the subset of permutations which
can be  obtained from $\pi$ by some composition  of the maps $a$,
$b$, and $\operatorname{Ad}_{\pi_0}$.
\end{Definition}

\begin{Remark}
A   Rauzy   class   $\mathfrak{R}(\pi)$  (extended  Rauzy   class
$\mathfrak{R}_{ex}(\pi)$)  of  a nondegenerate  permutation $\pi$
contains only nondegenerate permutations.
\end{Remark}

\begin{Thm}[W.A.Veech, \cite{V82}]
The extended Rauzy  classes  of nondegenerate permutations are in
the one-to-one  correspondence  with  the connected components of
the strata in the moduli  spaces  of  Abelian  differentials.
\end{Thm}

Using classification of the strata obtained in the current paper,
article~\cite{Z:represent} presents an explicit construction of a
representative of any extended Rauzy class.

\begin{Lm}[G.Rauzy, \cite{Rauzy}]
Any Rauzy class $\mathfrak{R}$ contains at  least one permutation
$\pi$ with the property
$$
\pi(m)=1 \qquad \pi(1)=m
$$
\end{Lm}

For the convenience of the reader  we give a sketch of the proof.
We want  to  fulfill  constraints $\pi(m)=1$ and $\pi^{-1}(m)=1$.
Suppose that it is not the case. Let us compare  numbers $\pi(m)$
and $\pi^{-1}(m)$.  If the smallest  of them is greater than $1$,
then applying one of operations $a$ or $b$ several times  one can
make  another  number  strictly  smaller. If the  smallest  among
$\pi(m)$ and $\pi^{-1}(m)$ is equal to $1$, then  applying one of
operations  $a$ or  $b$  several times one  can  make {\it  both}
numbers $\pi(m)$ and $\pi^{-1}(m)$ equal to $1$. \qed

We need the following modification of this Lemma.

\begin{Lemma}
\label{lm:m1adv}
Any  extended  Rauzy  class $\mathfrak{R}_{ex}$ of  nondegenerate
permutations contains at least one  permutation  $\pi$  with  the
following two properties
$$
\pi(m)=1 \qquad \pi(1)=m
$$
The permutation
$$
\pi':=
\begin{pmatrix}
\pi(2) \ttt{2} \pi(3) \ttt{2}  \cdots \ttt{2} \pi(m-2) \ttt{2} \pi(m-1) \\
    2  \ttt{2}     3  \ttt{2}  \cdots \ttt{2}     m-2  \ttt{2}     m-1
\end{pmatrix}
$$
obtained  as  a  restriction  of   $\pi$   to   the  ordered  set
$\{2,3,\dots,m-1\}$ is irreducible.
\end{Lemma}
\begin{proof}
Consider a permutation $\pi$ as  in  the  previous Lemma. Suppose
that  the  restriction  $\pi'$  of  $\pi$  to the ordered  subset
$\{2,3,\dots,m-1\}$  is  reducible. Choose  the  maximal  integer
$a<m-1$  such   that   $\pi'$   leaves  the  set  $\{2,\dots,a\}$
invariant. In other words  chose  the rightmost position where we
can break permutation $\pi'$ into two nonempty permutations.

Consider  the  following   ordered  subsets:
$$
\begin{array}{rcl}
A   & \! := &\! \{2,\dots,a\} \\
B_1 & \! := &\! \{a+1,\dots,\pi(m-1)-1\}\\
B_2 & \! := &\! \{\pi(m-1),\dots,m-1\}
\end{array}
$$
where  $B_1$  is an empty set when  $\pi(m-1)=a+1$.  Replace  the
initial permutation $\pi$ by  the  following one contained in the
same extended Rauzy class:
\begin{equation}
\label{eq:pi2}
\begin{pmatrix}
m & 1 &     & \pi(A) &   &  \pi(B_1) &   &   & \pi(B_2) &     & \\
  &   & B_2 &        & 1 &           & A & | &          & B_1 & m
\end{pmatrix}
\end{equation}
This  permutation   is   obtained   from   permutation  $\pi$  by
composition  of  the following  two  operations.  We  first  make
modification from the right by cyclically moving one step forward
the elements of the top line occurring after the letter $m$. Then
we  make  modification   from   the  left  by  cyclically  moving
$card(B_2)$  steps  forward  the  elements  of  the  bottom  line
occurring before letter $m$.

After reenumeration of the elements in the standard  order we see
that in this standard enumeration  our  new  permutation  $\pi_2$
again has the property $\pi_2(1)=m$ and $\pi_2(m)=1$. Restriction
$\pi_2'$ of this new permutation to  the subset $\{2,\dots,m-1\}$
may be again reducible. We are going to prove that the restricted
permutation may  split only to the right of  the marked place. In
other  words   we  are  going  to   prove  that  if   the  subset
$\{2,\dots,a_2\}$ is invariant under  $\pi'_2$  then $a_2 \ge a +
\Card{B_2} > a$.

Since the  initial  permutation  $\pi$  is  nondegenerate we have
$\pi(m-1)\neq  m-1$  (to  see  this  let  $j=m-1$  in  the second
condition on degenerate permutations at  the  beginning  of  this
section). Thus $\card{B_2}>1$. Hence the letter $1$ cannot be the
second letter in  the bottom line of~\eqref{eq:pi2}. Thus, if the
splitting occurs,  the  leftmost  invariant  subset contains more
than one element. Looking at  the  top  line of~\eqref{eq:pi2} we
see  that  this  means  that the leftmost invariant  subset  must
contain at  least one element  of $\pi(A)$. Looking at the bottom
line we see that the leftmost invariant subset  contains at least
one  element  of $B_2$. Note that  the  set $A$ considered as  an
unordered set was chosen to  be  invariant  under the permutation
$\pi$. Thus $\pi(A)$ does not intersect  with  $B_2$.  Hence  the
splitting may occur only to the right of the word $\pi(A)$ in the
top line. Thus the leftmost invariant subset must contain all the
elements of the  unordered set $\pi(A)=A$. Thus the splitting may
mapped only to  the  right of the marked  position  at the bottom
line.

Repeating  inductively   this  procedure  we  finally  obtain  an
irreducible restricted permutation.
\end{proof}

\subsection{Zippered rectangles (after W.A.Veech)}
\label{ss:ziprec}
Having an interval exchange transformation  $T:I\to  I$  one  can
``suspend'' a smooth closed complex  curve  $C_g$  and an Abelian
differential $\omega$  over $T$. Here  we present the idea of the
``suspension'';  one  can find all the details  in  the  original
paper of W.A.Veech~\cite{V82}.

\begin{figure}[ht]
%
%
%
\includegraphics{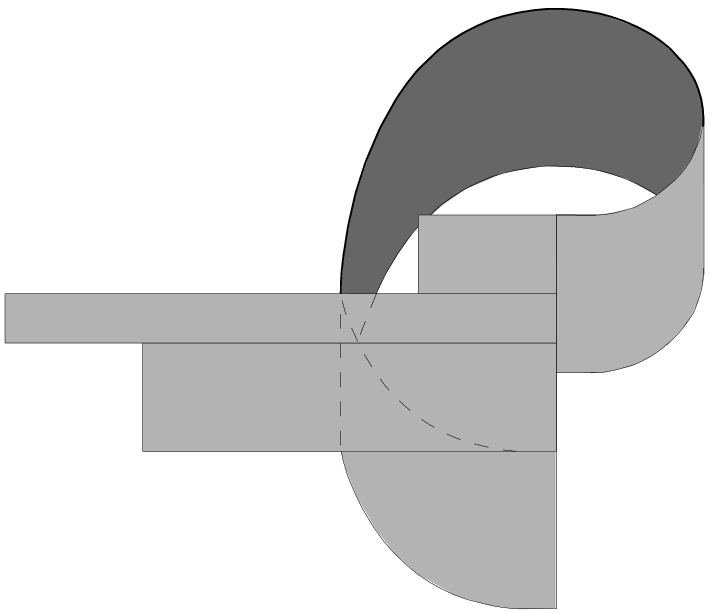}
\includegraphics{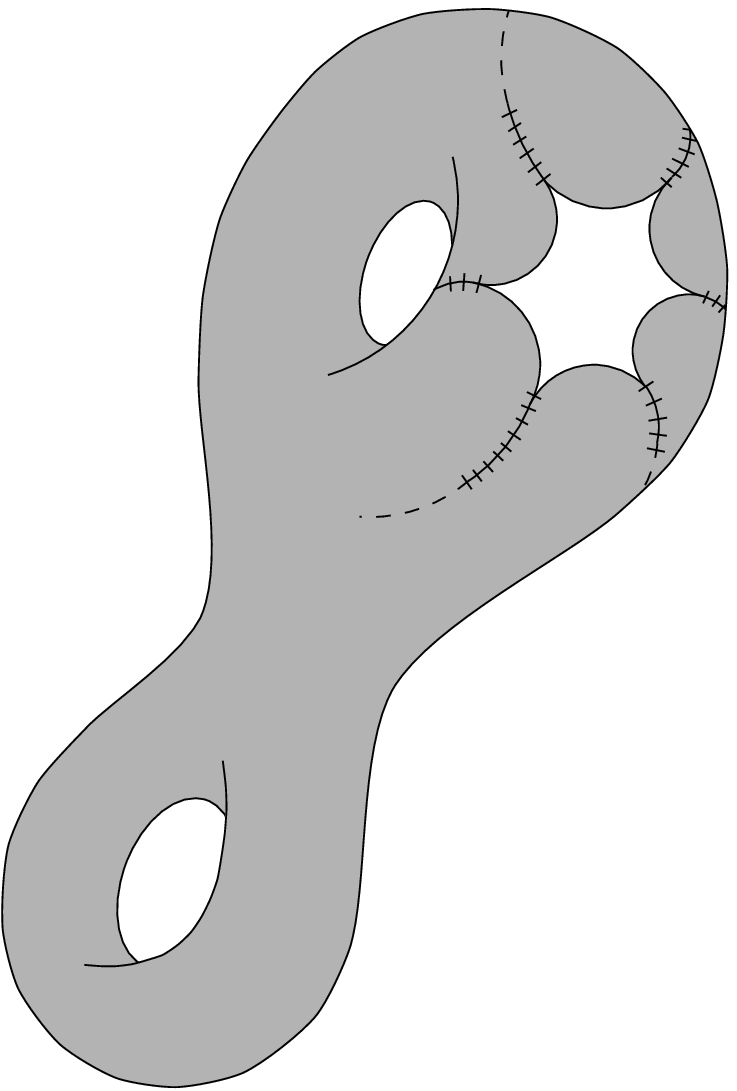}
\vspace{165bp} 
\caption{
\label{pic:zr}
Suspension over  the  interval  exchange  transformation with the
permutation $\pi=\{4,3,2,1\}$ produces a surface of  genus 2 with
an Abelian differential having single zero of order 2.
}
\end{figure}

\begin{figure}[ht]
%
%
\includegraphics{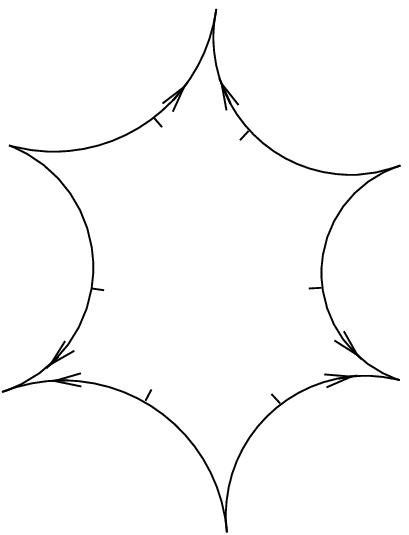}
\vspace{130bp}
\begin{picture}(0,0)(0,0)
\put(0,0){
\begin{picture}(0,0)(0,0)
\put(7,80){$h_2$}
\put(9,24){$h_3$}
\put(52,8){$h_4+h_1$}
\put(108,24){$h_2$}
\put(113,80){$h_3$}
\put(52,103){$h_1+h_4$}
\put(97,-5){$i$}
\put(145,56){$ii$}
\put(99,122){$iii$}
\put(28,122){$iv$}
\put(-15,54){$v$}
\put(29,-5){$vi$}
\put(-120,45){\shortstack{1\\4}}
\put(-105,45){\shortstack{2\\3}}
\put(-90,45){\shortstack{3\\2}}
\put(-75,45){\shortstack{4\\1}}
\multiput(-113.1,57)(15,0){3}{$\vee$}
\multiput(-113.1,43)(15,0){3}{$\wedge$}
\put(-113,65){\footnotesize\it i}
\put(-112,35){\footnotesize\it vi}
\put(-98,65){\footnotesize\it v}
\put(-97,35){\footnotesize\it ii}
\put(-83,65){\footnotesize\it iii}
\put(-84,35){\footnotesize\it iv}
\end{picture}}
\end{picture}
\caption{
\label{pic:hole}
The lengths  of the sides of the  hole for  a suspension over  an
interval   exchange    transformation   with   the    permutation
$\pi=\{4,3,2,1\}$.
}
\end{figure}

Place  the  interval $I$ horizontally in the plane  $\R{2}=\C{}$.
Place a rectangle $R_i$ over each subinterval $I_i\subset I$; the
rectangle $R_i$ has the width $\lambda_i=|I_i|$ and some altitude
$h_i$. Later on we shall pose some restrictions on the altitudes.
Glue the  top horizontal side  of rectangle $R_i$ to the interval
$T(I_i)$ at the  base. There are still no identifications between
the vertical sides of the rectangles, so we get a Riemann surface
with several ``holes''; each boundary component is a union of the
vertical       sides       of      the       rectangles      (see
Figures~\ref{pic:zr}, \ref{pic:hole}). Now  start ``zipping'' the
holes  (see  Figure~\ref{pic:zr}). If the altitudes $h_i$ of  the
rectangles, and the  altitudes $a_i$ till which we ``zipper'' the
rectangles   obey   some  linear   equations   and   inequalities
(see~\cite{V82}), then we  manage to eliminate all the holes. The
Riemann surface thus constructed has natural  flat structure with
cone-type singularities;  the  complex structure, coming from the
initial complex structure  on  the plane $\C{}=\R{2}$, extends to
the conical points. The Abelian differential  $\omega$ is locally
represented  as  $dz$, where $z$ is the  standard  coordinate  in
$\C{}$.

As we already mentioned the altitudes $h_i$, and  $a_i$ obey some
linear  relations  (cf.  Figure~\ref{pic:hole});  it  is   proved
in~\cite{V82}  that  the family of solutions is always  nonempty.
This family  has dimension $m=2g+k-1=\dim  H^1(C_g,\{\text{zeroes
of   }\omega\})$,   which  coincides  with  the  number  $m$   of
subintervals under exchange, $\pi\in\mathfrak{S}_m$.

\section
{Abelian differentials on hyperelliptic curves}
\label{s:hypel}

Let  $\omega$  be an  Abelian  differential  on  a  hyperelliptic
complex  curve  such that  all  zeroes of  $\omega$  are of  even
degrees. Let the canonical divisor $K(\omega)$ of $\omega$ be
\begin{equation}
\label{eq:komega}
K(\omega)= 2(k_1 P_{i_1}+ \dots +  k_p P_{i_p})+
2\left( l_1(P^+_1 + P^-_1) + \dots + l_q(P^+_q + P^-_q)\right)
\end{equation}
where
$$
\sum_{i=1}^p k_i + 2\sum_{j=1}^q l_j = g-1
$$
By  $P_{i_n}$  we  denote  the points which are  invariant  under
hyperelliptic involution; by $P^\pm_j$ we  denote  the  pairs  of
points symmetrical  to each other under hyperelliptic involution.
We assume that all the indicated points are distinct.

\begin{Proposition}
\label{pr:hspin}
The  parity  of  the   spin   structure  defined  by  an  Abelian
differential on a hyperelliptic  curve  is given by the following
formula:
$$ \varphi(\omega) \equiv \dim \left| \frac{1}{2} K(\omega)
\right|+1 (\mod 2) = \sum_{i=1}^p \left[ \frac{k_i}{2}
\right] + \sum_{j=1}^q l_q +1 (\mod 2) $$
\end{Proposition}
\begin{proof}
In  our   case  a  base  of   solutions  of  the   linear  system
$\cfrac{1}{2} K(\omega)$ can be constructed explicitly.
\end{proof}

\begin{Corollary}
\label{cr:hspin}
Parity  of  the  spin  structure   determined   by   an   Abelian
differential      from      the      hyperelliptic      component
$\mathcal{H}^{hyp}(2g-2)$ equals
$$ \varphi(\mathcal{H}^{hyp}(2g-2)) \equiv
\left[\frac{g+1}{2}\right] (\mod 2) $$
Parity  of  the  spin  structure of the  hyperelliptic  component
$\mathcal{H}^{hyp}(g-1,g-1)$, for odd genera $g$ equals
$$ \varphi\left(\mathcal{H}^{hyp}(g-1,g-1)\right) \equiv
\left(\frac{g+1}{2}\right) (\mod 2) \quad for\ \ odd\ \ g $$
\end{Corollary}

\centerline{\sc{Acknowledgments}}

The authors thank M.~Duchin, M.~Farber, I.~Itenberg, M.~Kazarian,
E.~Lanneau, and  the  referee  for  helpful  comments. The second
author is grateful to MPI f\"ur Mathematik at Bonn, to FIM of ETH
at Z\"urich, and to IHES for  hospitality  while  preparation  of
this paper, as well as to CNRS Projects 5376, 7726 for support of
collaboration between University of Rennes and Moscow Independent
University.


\end{document}